\documentclass[11pt]{article}
\usepackage{amsmath, amssymb, amscd}

\newtheorem{thm}{Theorem}
\newtheorem{prop}[thm]{Proposition}
\newtheorem{lem}[thm]{Lemma}

\newtheorem{cor}[thm]{Corollary}
\newtheorem{rem}[thm]{Remark}
\newtheorem{df}[thm]{Definition}

\renewcommand{\epsilon}{\varepsilon}
\renewcommand{\phi}{\varphi}
\renewcommand{\deg}{\widetilde{\operatorname{deg}}}
\renewcommand{\P}{\operatorname{P}}
\renewcommand{\S}{\operatorname{S}}

\newcommand{\BB}{\mathbb}

\newcommand{\pf}{\noindent {\it Proof. }}
\newcommand{\qed}{\nopagebreak $\qquad$ $\square$ \vskip5pt}
\newcommand{\separate}{\vskip5pt}
\newcommand{\supp}{\operatorname{supp}}
\newcommand{\im}{\operatorname{Im}}
\newcommand{\re}{\operatorname{Re}}
\newcommand{\tr}{\operatorname{Tr}}

\newcommand{\B}{\overline}
\newcommand{\HC}{\BB H_{\BB C}}
\newcommand{\HR}{\BB H_{\BB R}}

\newcommand{\sgn}{\operatorname{sign}}
\newcommand{\sgnx}{\operatorname{sign}_3}

\usepackage[OT2,T1]{fontenc}
\newcommand\textcyr[1]{{\fontencoding{OT2}\fontfamily{wncyr}\selectfont #1}}
\newcommand{\Zh}{\textit{\textcyr{Zh}}}

\begin{document}

\title{\bf Split Quaternionic Analysis and Separation of the Series for
$SL(2,\BB R)$ and $SL(2,\BB C)/SL(2,\BB R)$}
\author{Igor Frenkel and Matvei Libine}
\maketitle

\begin{abstract}
We extend our previous study of quaternionic analysis based on
representation theory to the case of split quaternions $\HR$.
The special role of the unit sphere in the classical quaternions
$\BB H$ -- identified with the group $SU(2)$ -- is now played by the group
$SL(2,\BB R)$ realized by the unit quaternions in $\HR$.
As in the previous work, we use an analogue of the Cayley transform to
relate the analysis on $SL(2,\BB R)$ to the analysis on the imaginary
Lobachevski space $SL(2,\BB C)/SL(2,\BB R)$ identified with the one-sheeted
hyperboloid in the Minkowski space $\BB M$.
We study the counterparts of Cauchy-Fueter and Poisson formulas on $\HR$ and
$\BB M$ and show that they solve the problem of separation of the discrete
and continuous series.
The continuous series component on $\HR$ gives rise to the minimal
representation of the conformal group $SL(4,\BB R)$, while the discrete series
on $\BB M$ provides its $K$-types realized in a natural polynomial basis.
We also obtain a surprising formula for the Plancherel measure of $SL(2,\BB R)$
in terms of the Poisson-type integral on the split quaternions $\HR$.
Finally, we show that the massless singular functions of four-dimensional
quantum field theory are nothing but the kernels of projectors onto the
discrete and continuous series on the imaginary Lobachevski space
$SL(2,\BB C)/SL(2,\BB R)$.
Our results once again reveal the central role of the Minkowski space
in quaternionic and split quaternionic analysis as well as a deep connection
between split quaternionic analysis and the four-dimensional quantum field
theory.
\end{abstract}


\section{Introduction}

Six years after William Rowan Hamilton's fundamental discovery of quaternions,
in 1849 James Cockle introduced a related algebra of ``coquaternions'',
which in modern language is called the algebra of split quaternions
$\HR$. As an algebra, $\HR$ is isomorphic to the real $2\times2$ matrices
and has signature $(2,2)$ with respect to the determinant viewed as a
quadratic form. Unlike Hamilton's quaternions, not every non-zero element in
$\HR$ is invertible. However, there are strong parallels between the geometry
and analysis on the two algebras.

In our previous paper \cite{FL} we began studying analysis on the algebra
of quaternions $\BB H$ from the point of view of representation theory of
the quaternionic conformal group $SL(2,\BB H)$.
It is well known that the spaces of harmonic, left- and right-regular
functions admit natural actions of $SL(2,\BB H)$ by fractional linear
transformations.
We regarded $\BB H$ as a real form of the space of complex quaternions $\HC$,
and the questions of unitarity of these representations led us
to another real form in $\HC$, namely the Minkowski space $\BB M$.
The conformal group of $\BB M$ is naturally identified with
$SU(2,2)$ and acts by unitary transformations on the spaces of harmonic,
left- and right-regular functions on $\BB M$.
By harmonic functions on $\BB M$ we mean the solutions of the wave equation;
similarly, left- and right-regular functions on $\BB M$ are solutions of
the left and right Weyl equations, which together form the massless Dirac
equations. To relate these two real forms of $\HC$ we considered the Cayley
transform which maps $\BB M$ into $U(2)$ and the unit two-sheeted hyperboloid
$H^3 \subset \BB M$ into the unit sphere $S^3 \subset \BB H$.
The quaternionic picture, however, had an important advantage over $\BB M$
-- it gave a natural realization of the $K$-types of the unitary
representations of $SU(2,2)$ as the spaces of polynomials in quaternionic
coordinates. While the group $SU(2,2)$ does not act naturally on $\BB H$,
its maximal compact subgroup $K= S(U(2) \times U(2))$ does act there and one
can identify the basis of $K$-types with the matrix coefficients of
the group $U(2)$ and its subgroup $SU(2)$.

In this work we consider analysis on another real form of
$\HC$, namely the split quaternions $\HR$, and study the action of
the conformal group $SL(2,\HR) \simeq SL(4,\BB R)$.
Just as the quaternionic analysis on $\BB H$ is based on analysis on the
unit sphere identified with the group $SU(2)$, the split quaternionic
analysis on $\HR$ has at its foundation analysis on the unit hyperboloid
identified with the group $SL(2,\BB R) \simeq SU(1,1)$.
It is a classical result that the analysis on $SL(2,\BB R)$ has the discrete
and continuous components corresponding to the two types of irreducible
unitary representations of this group.
All results in this paper will come in two flavors -- for left-/right-regular
functions and solutions of the appropriate wave equation.
For simplicity, in this introduction we only announce the results for
left-regular functions.

One of the central results of our paper is a solution of the problem of
separation of the series on $SU(1,1)$ using the main formulas of
quaternionic analysis on $\HR$, namely the counterparts of the Cauchy-Fueter
and Poisson integrals.
Again, the relationship between the analysis on $\HR$ and $\BB M$ is an
important tool of our study.
We use an analogue of the Cayley transform previously studied in \cite{KouO}
which now maps $\BB M$ into $U(1,1)$ and the unit one-sheeted hyperboloid
$H^3_{\BB M} \subset \BB M$ into the unit hyperboloid $H^3_{\BB R} \subset \HR$
(Proposition \ref{Cayley}).
There are natural identifications $H^3_{\BB M} \simeq SL(2,\BB C)/SU(1,1)$
and $H^3_{\BB R} \simeq SU(1,1)$, and this Cayley transform interchanges the
discrete and continuous components on $H^3_{\BB M}$ and $H^3_{\BB R}$.
The study of harmonic analysis on $SL(2,\BB C)/SU(1,1)$ -- also known as the
imaginary Lobachevski space -- goes back to Gelfand and his school \cite{GGV}
and is based on the geometry of horospheres.
These methods were used for the problem of separation of the series on
$SU(1,1)$ in more recent works by Gindikin \cite{Gi1,Gi2}.
(See also references therein for related work on this problem.)
In our approach we realize functions on $SL(2,\BB C)/SU(1,1)$ and $SU(1,1)$
using harmonic extensions to the flat spaces $\BB M$ and $\HR$.
This idea was previously explored in greater generality by Strichartz \cite{St}.
In our work the methods and formulas of quaternionic analysis are naturally
applied to solve the problems of harmonic analysis on $SU(1,1)$ and
$SL(2,\BB C)/SU(1,1)$.

We study first the discrete series component of split quaternionic analysis,
which goes in a strong parallel with the case of (classical) quaternions
following the results and constructions of our paper \cite{FL}.
It turns out that the discrete components of the harmonic\footnote
{By harmonic functions on $\HR$ we mean the solutions of the ultrahyperbolic
wave equation, which is in parallel with the Minkowski case $\BB M$.},
left- and right-regular functions on $\HR$ yield the same unitary
representations of the group $SU(2,2)$ as in the (classical) quaternionic case.
The unitary structure and the identification of the representations become
transparent after the Cayley transform relating $\BB M$ with $U(1,1)$
and mapping the tube domains $\BB T^{\mp}$ in $\HC$
(see Section 3.5 in \cite{FL}) into certain complex semigroups $\Gamma^{\mp}$
with the Shilov boundary $U(1,1)$. These are Ol'shanskii semigroups
which were first introduced in \cite{Vin}, \cite{Ol}.

One can obtain again polynomial realizations of the discrete components of
the harmonic and regular functions on $\HR$ by restricting the group $SU(2,2)$
to another subgroup $K_{\BB R} = S(U(1,1) \times U(1,1))$ and identifying the
polynomial basis of $K_{\BB R}$-representations with the matrix coefficients
of the discrete component of the group $U(1,1)$ and its subgroup $SU(1,1)$.
In particular, the two irreducible representations ${\cal V}^{\mp}$ of
$SU(2,2)$ realized in the space of left-regular functions
on $\HR$ are decomposed with respect to the shifted degree operator
$\deg = \operatorname{deg}+1$ as follows:
\begin{equation}  \label{V(n)}
{\cal V}^- = \bigoplus_{n \in \BB Z} {\cal V}^-(n), \qquad
{\cal V}^+ = \bigoplus_{n \in \BB Z} {\cal V}^+(n)
\end{equation}
where ${\cal V}^-(n)$ (respectively ${\cal V}^+(n)$), $n\ne 0$, are the
irreducible representations from the holomorphic (respectively antiholomorphic)
discrete series, and the case $n=0$ corresponds to the limit
of the holomorphic (respectively antiholomorphic) discrete series.
Then the Fueter formula yields projections $\P^{\mp}$ of the space of
left-regular functions $f$ onto the holomorphic or antiholomorphic discrete
series components depending on the domain of the variable $W \in \Gamma^{\mp}$ in
\begin{equation}  \label{H_R-proj-intro}
-\frac1{2\pi^2}
\int_{Z \in H^3_{\BB R}} \frac {(Z-W)^{-1}}{\det(Z-W)} \cdot *dZ \cdot f(Z)
= \bigl(\P^{\mp}f_{>0} - \P^{\mp}f_{<0}\bigr)(W),
\end{equation}
where $*dZ$ is a certain naturally defined quaternionic-valued differential
3-form on $\HR$ and $\P^{\mp}f_{>0}$ (respectively $\P^{\mp}f_{<0}$) denotes the sum
of the components in (\ref{V(n)}) of positive (respectively negative) shifted
degree (Theorem \ref{Hilbert-transform}).
Note that the components $\P^{\mp}f_{>0}$ and $\P^{\mp}f_{<0}$ enter with opposite
signs and the limits of the discrete series would cause the integral to be
divergent. (In the case of Poisson formula the limits of the discrete series
are projected out.)
Since $H^3_{\BB R} \simeq SU(1,1)$ lies inside the Shilov boundary of
$\Gamma^{\mp}$, the values of the functions on $H^3_{\BB R}$ can be recovered
by taking limits.
The Fueter formula for the right-regular functions and the Poisson formula
for harmonic functions have similar structures.
The results and formulas for the middle series of representations of $SU(2,2)$
obtained in \cite{FL} have even more exact analogues for the discrete
component of split quaternionic analysis.
In particular, the limits of the discrete series do not occur in the
restrictions of these representations to $K_{\BB R}$.


Fundamentally new features of split quaternionic analysis appear
when we study the continuous series component. Now the space of harmonic
functions on $\HR$ of the continuous series component gives rise to a single
irreducible representation of the conformal group
$SL(4,\BB R)/\{\pm 1\} \simeq SO(3,3)$, known as the minimal
representation\footnote
{Strictly speaking, $SO(3,3)$ does not have a minimal representation,
so by ``minimal representation of $SO(3,3)$'' we mean the representation
$(\varpi^{p,q}, V^{p,q})$ of $O(p,q)$ in the notations of \cite{KobO} for $p=q=3$.
When $p+q$ is an even number greater than or equal $8$,
one gets a genuine minimal representation.}.
It was studied for an arbitrary signature in
\cite{KobO} (also see references therein for the previous work on this subject).
Again, the representation theory of $SL(4,\BB R)$ is crucial for studying
the continuous component of split quaternionic analysis.
On the other hand, the latter illuminates various aspects
of the representation theory.
In particular, applying the Cayley transform to the space of harmonic
functions on $\HR$ spanning the minimal representation, we realize the minimal
representation in a space ${\cal D}^0_{\BB M}$ with basis consisting of certain
harmonic polynomials on $\BB M$.
This basis appears naturally when we restrict the minimal representation of
$SL(4,\BB R)$ to its subgroup $K_{\BB M} = \BB R^{>0} \times SL(2,\BB C)$.
Decomposing ${\cal D}^0_{\BB M}$ with respect to the shifted degree operator
$\deg$ on $\BB M$ yields
$$
{\cal D}^0_{\BB M} = \bigoplus_{n \in \BB Z} {\cal D}^0_{\BB M}(n),
$$
where each ${\cal D}^0_{\BB M}(n)$, $n \in \BB Z$, is an irreducible
representation of $SL(2,\BB C)$ and can be decomposed further into the
irreducible polynomial subspaces with respect to the compact subgroup
$SU(2) \subset SL(2,\BB C)$.
This description also gives an explicit realization of the $K$-types of
the minimal representation  relative to the maximal compact subgroup
$K'=SO(4) \subset SL(4,\BB R)$.
The Poisson and Cauchy-Fueter integrals provide us projections onto the
discrete component on $\BB M$, but now the procedure is more subtle than in
the case of discrete series projectors in $\HR$ such as (\ref{H_R-proj-intro}).
The projection is obtained as a boundary value of these integrals in the
sense of hyperfunctions, namely as a limit of the difference of two integrals.
Thus for the Fueter integral we obtain the following result:
\begin{equation}  \label{M-discr-intro-1}
\lim_{\epsilon \to 0^+} \bigl( \Phi(W,i\epsilon) - \Phi(W,-i\epsilon) \bigr)
= \bigl(\P^{discr(\BB M)} f_{>0} - \P^{discr(\BB M)} f_{<0}\bigr)(W),
\end{equation}
where
\begin{equation}  \label{M-discr-intro-2}
\Phi(W, i\epsilon) = \frac i{2\pi^2} \int_{Y \in H^3_{\BB M}}
\frac{(Y-W)^+}{(\det(Y-W) +i\epsilon)^2} \cdot *dY \cdot f(Y)
\end{equation}
and $\P^{discr(\BB M)}$ is the projection of the space of left-regular functions on
$\BB M$ onto the discrete series component (Theorem \ref{M-discr-proj-reg}).
Note that, as in the case of the discrete series projectors on $\HR$
(\ref{H_R-proj-intro}), $\P^{discr(\BB M)} f_{>0}$ and $\P^{discr(\BB M)} f_{<0}$
enter with opposite signs and the functions $f$ with $\deg f=0$ would cause
the integrals to be divergent. (In the case of Poisson formula the functions
in ${\cal D}^0_{\BB M}(0)$ are projected out.)
Thus we can view ${\cal D}^0_{\BB M}(0)$ as the limit of the discrete series
on $H^3_{\BB M} \simeq SL(2,\BB C)/SU(1,1)$.

Finally, we return to the $\HR$-setting and study the continuous series
component by applying the Cayley transform to the discrete series component on
$\BB M$. Quite surprisingly, the Poisson and Cauchy-Fueter integrals do not
give the projectors as in (\ref{M-discr-intro-1}), (\ref{M-discr-intro-2}),
but instead become diagonal operators commuting only with the subgroup
$S(GL(2,\BB R) \times GL(2,\BB R))$ of the conformal group $SL(4,\BB R)$,
and the diagonal density is precisely the inverse of the Plancherel density
of $SL(2,\BB R)$! Recall that in the parameterization of the continuous
series by $\chi=(l,\epsilon)$, where $l=-\frac12+i\lambda$, $\lambda \in \BB R$
and $\epsilon \in \{0,\frac12\}$, the Plancherel density is given by
$$
Pl(\chi) = \begin{cases}
\lambda \tanh(\pi\lambda) & \text{if $\epsilon=0$;}  \\
\lambda \coth(\pi\lambda) & \text{if $\epsilon=\frac12$.}
\end{cases}
$$
(See, for example, \cite{Va}.)
Then we have the following identity
\begin{equation}  \label{H-cont-intro-1}
\lim_{\epsilon \to 0^+} \bigl(
\tilde\Phi(W,i\epsilon) - \tilde\Phi(W,-i\epsilon) \bigr)
= \frac {\phi_{\chi}(W)}{Pl(\chi)}, \qquad W \in SU(1,1),
\end{equation}
where
\begin{equation}  \label{H-cont-intro-2}
\tilde\Phi(W, i\epsilon) = \frac1{2\pi^2} \int_{X \in SU(1,1)}
\frac{\phi_{\chi}(X)}{\det(X-W) +i\epsilon} \,dX^3
\end{equation}
and $\phi_{\chi}(X)$ belongs to the linear span of matrix coefficients of the
irreducible representation corresponding to $\chi=(l,\epsilon)$
(Theorem \ref{cont_ser_proj-2}).
We remark that for the discrete series of $SL(2,\BB R)$ the Plancherel density
is $-(2l+1)$ and we could rewrite the Poisson integral for the discrete
series projector in a form similar to
(\ref{H-cont-intro-1})-(\ref{H-cont-intro-2}) using the kernels
$\bigl( \det(X-W) \pm i\epsilon \bigr)^{-1}$.
Thus quaternionic analysis over the split quaternions $\HR$ naturally
reveals the important salient features of harmonic analysis on $SL(2,\BB R)$.

In \cite{FL2} we differentiate the family of operators $\operatorname{Pl}_R$
introduced in Theorem \ref{cont_ser_proj-2} with respect to $R$ and show that
the effect of $\frac{d}{dR} \operatorname{Pl}_R$ on the discrete and continuous
series components can be easily computed using the Schr\"odinger model for
the minimal representation of $O(3,3)$ and the results of Kobayashi-Mano from
\cite{KobM}, particularly their computation of the integral expression for the
operator ${\cal F}_C$.
This provides an independent verification of the coefficients involved
in Theorems \ref{cont_ser_proj}, \ref{cont_ser_proj-2}.

To summarize our studies of quaternionic analysis on $\BB H$ in \cite{FL}
and split quaternionic analysis on $\HR$ in this work we would like to
point out once again the special role of the Minkowski space $\BB M$ related
to the two algebras of quaternions by the two types of Cayley transform.
Thus analysis on $\BB H$ and $\HR$ is equivalent to analysis on the
two-sheeted and one-sheeted hyperboloids in $\BB M$ respectively.
Moreover, doing analysis on $\BB M$ in many ways facilitates and clarifies
analysis on $\BB H$ and $\HR$.
The relation between $\BB H$, $\HR$ and $\BB M$ also reveals the hierarchy
of the symmetry groups of these three spaces, which can be presented as follows:
$$
\begin{array}{ccccccccc}
SO(5,1) & \: & \: & \: & SO(4,2) & \: & \: & \: & SO(3,3) \\
\: & \diagdown & \: & \diagup & \: & \diagdown & \: & \diagup & \: \\
\: & \: & SO(4,1) & \: & \: & \: & SO(3,2) & \: & \:  \\
\: & \diagup & \: & \diagdown & \: & \diagup & \: & \diagdown & \: \\
SO(4) & \: & \: & \: & SO(3,1) & \: & \: & \: & SO(2,2) \\
\text{\Large$\wr$} & \: & \: & \: & \text{\Large$\wr$} & \: & \: & \: &
\text{\Large$\wr$} \\
\BB H \simeq \BB R^4 & \: & \: & \: & \BB M \simeq \BB R^{3,1} & \: & \: & \: &
\HR \simeq \BB R^{2,2}
\end{array}
$$
(In this diagram some groups are replaced with locally isomorphic ones.)
Clearly, the groups at the bottom of the diagram are the metric-preserving
transformations of the three real forms of $\HC$ of our interest,
while the groups at the top row are the corresponding conformal groups.
The groups in the middle row appear as subgroups of the conformal groups
preserving the hyperboloids. Thus the Minkowski space literally plays the
central role in quaternionic analysis on $\HC$ and its real forms
$\BB H$ and $\HR$.

In the same way as various results in quaternionic analysis on $\BB H$
can be generalized to Clifford analysis on arbitrary Euclidean spaces (the
book \cite{GSp} contains a comprehensive bibliography on these subjects),
split quaternionic analysis on $\HR$ as well as analysis on the Minkowski
space $\BB M$ can be further extended to Clifford analysis on real vector
spaces of arbitrary signature.
However, quaternionic analysis on $\BB H$, $\HR$ and $\BB M$ is singled out
because of its relation to the harmonic analysis on the simplest groups
$SU(2)$, $SL(2,\BB R)$ and homogeneous spaces $SL(2,\BB C)/SU(2)$,
$SL(2,\BB C)/SL(2,\BB R)$ as well as uniqueness of the space of quaternions
$\BB H$.
Moreover, we believe that the special role of the four-dimensional case will
reveal itself at deeper levels in the future developments of the subject
when the relation with exotic smoothness in four-dimensional topology and
renormalization in four-dimensional quantum field theory becomes more
transparent.

By design, analysis on the Minkowski space has many connections with the
four-dimensional physics.
In our previous paper \cite{FL} we noted several relations between
representation theory underlying quaternionic analysis and Feynman diagrams.
Split quaternionic analysis adds another important connection between the
mathematical and physical structures.
To make the connections with physics more apparent, in both cases we consider
various results of quaternionic analysis over $\BB H$ and $\HR$ in the
Minkowski space realization.
The projectors onto the continuous and discrete components on $\BB M$
are expressed using the kernels:
\begin{equation}  \label{qft-kernels}
\frac1{\det Z + i\epsilon z^0}, \qquad \frac1{\det Z - i\epsilon z^0},
\qquad \frac1{\det Z + i\epsilon}, \qquad \frac1{\det Z - i\epsilon},
\qquad Z=Y-W,
\end{equation}
$z^0$ denotes the time coordinate of $Z$,
in the case of Poisson integrals and their $\nabla_Y$-derivatives in the case
of Cauchy-Fueter integrals.
But these are precisely the massless singular functions that form the
foundation of the theory of Feynman integrals! (See, for example, \cite{BSh}.)
Any other singular function appearing in physics is a linear combination of
the ones listed in (\ref{qft-kernels}).
Since there is a linear relation between these kernels:
$$
\frac1{\det Z + i\epsilon z^0} + \frac1{\det Z - i\epsilon z^0}
= \frac1{\det Z + i\epsilon} + \frac1{\det Z - i\epsilon},
$$
there are only three essentially different kernels in the list
(\ref{qft-kernels}), which is exactly the number of
$\mathfrak{su}(2,2)$-irreducible components of the space of functions on
imaginary Lobachevski space $SL(2,\BB C)/SU(1,1)$!
Thus the problem of separation of the series on $SL(2,\BB R)$ and
$SL(2,\BB C)/SU(1,1)$, which was solved by means of quaternionic analysis,
inevitably leads to the heart of four-dimensional quantum field theory!
Unfortunately, our representation-theoretic approach does not naturally include
the massive versions of the singular functions in (\ref{qft-kernels}).
However, in our forthcoming work \cite{FL3} we will show that,
in spite of the enormous rigidity of quaternionic analysis, it is possible
to make a one-parameter deformation that we expect to include the missing mass.

The paper is organized as follows. In Section 2 we introduce notations and
overview basics of split quaternionic analysis.
In particular, we derive an analogue of the Cauchy-Fueter formula for regular
functions on $\HR$. We also give the bases for the spaces of harmonic,
left- and right-regular on $\HR$.
In Section 3 we study the discrete series component of split quaternionic
analysis (over $\HR$). First we find the polynomial bases of the spaces of
harmonic, left- and right-regular functions ${\cal D}$, ${\cal V}$,
${\cal V}'$ and describe the action of the Lie algebra
$\mathfrak{sl}(4,\BB C)$ in those bases.
After a brief review of Ol'shanskii semigroups, we obtain expansions of the
Poisson and Fueter kernels in terms of the matrix coefficients of the discrete
series of $SU(1,1)$. From these expansions we derive the discrete series
projectors.
In Section 4 we prove that the Cayley transform between $\HR$ and $\BB M$
interchanges the discrete and continuous series components on $\HR$ and $\BB M$.
Using results of Section 3 we obtain the continuous series projectors onto the
same spaces ${\cal D}$, ${\cal V}$ and ${\cal V}'$ in the new
setting of the Minkowski space.
Then we study the minimal representation of $SL(4,\BB R)$ realized
in the space of harmonic functions on $\BB M$, find a polynomial basis for
the $K$-types and use that basis to derive the discrete series projector.
In Section 5 we study the Poisson integrals and their effect on the discrete
and continuous series on $\HR$. We end this section with the proof that
the boundary value of the Poisson integral in the sense of hyperfunctions
yields a diagonal operator with density given by the inverse of the Plancherel
measure of the group $SL(2,\BB R)$.
In the Appendix we give a brief introduction to the special functions that
we use throughout the paper and list their properties.

The first author was supported by the NSF grants DMS-0457444 and DMS-1001633;
the second author was supported by the NSF grant DMS-0904612.

\section{Split Quaternionic Analysis}

\subsection{The Quaternionic Spaces $\HC$, $\HR$ and $\BB M$}

In this article we use notations established in \cite{FL}.
In particular, $e_0$, $e_1$, $e_2$, $e_3$ denote the units of the classical
quaternions $\BB H$ corresponding to the more familiar $1$, $i$, $j$, $k$
(we reserve the symbol $i$ for $\sqrt{-1} \in \BB C$).
Thus $\BB H$ is an algebra over $\BB R$ generated by $e_0$, $e_1$, $e_2$, $e_3$,
and the multiplicative structure is determined by the rules
$$
e_0 e_i = e_i e_0 = e_i, \qquad
(e_i)^2 = - e_0, \qquad
e_ie_j=-e_ie_j, \qquad 1 \le i< j \le 3,
$$
and the fact that $\BB H$ is a division ring.
Next we consider the algebra of complexified quaternions
$\HC = \BB C \otimes \BB H$.
We define a complex conjugation on $\HC$ with respect to $\BB H$:
$$
Z = z^0e_0 + z^1e_1 + z^2e_2 + z^3e_3 \quad \mapsto \quad
Z^c = \B{z^0}e_0 + \B{z^1}e_1 + \B{z^2}e_2 + \B{z^3}e_3,
\qquad z^0,z^1,z^2,z^3 \in \BB C,
$$
so that $\BB H = \{ Z \in \HC ;\: Z^c=Z \}$.
The quaternionic conjugation on $\HC$ is defined by:
$$
Z = z^0e_0 + z^1e_1 + z^2e_2 + z^3e_3 \quad \mapsto \quad
Z^+ = z^0e_0 - z^1e_1 - z^2e_2 - z^3e_3,
\qquad z^0,z^1,z^2,z^3 \in \BB C;
$$
it is an anti-involution:
$$
(ZW)^+ = W^+ Z^+, \qquad \forall Z,W \in \HC.
$$
We will also use an involution
$$
Z \mapsto Z^- = - e_3 Z e_3 \qquad \text{(conjugation by $e_3$).}
$$
Then the complex conjugation, the quaternionic conjugation and the involution
$Z \mapsto Z^-$ commute with each other.

In this article we will be primarily interested in the space of
{\em split quaternions} $\HR$ which is a real form of $\HC$ defined by
$$
\HR = \{ Z \in \HC ;\: Z^{c-} = Z \} = \{\text{$\BB R$ -span of
$e_0$, $\tilde e_1 = i e_1$, $\tilde e_2 = -i e_2$, $e_3$} \}
$$
and the Minkowski space $\BB M$ which is regarded as another real form of $\HC$:
$$
\BB M = \{Z \in \HC ;\: Z^{c+} = - Z \}
= \{\text{$\BB R$ -span of $\tilde e_0 = -i e_0$, $e_1$, $e_2$, $e_3$} \}.
$$

On $\HC$ we have a quadratic form $N$ defined by
$$
N(Z) = ZZ^+ = Z^+Z = (z^0)^2 + (z^1)^2 + (z^2)^2 + (z^3)^2,
$$
hence $Z^{-1} = Z^+/N(Z)$.
The corresponding symmetric bilinear form on $\HC$ is
\begin{equation}  \label{bilinear_form}
\langle Z, W \rangle = \frac12 \tr(Z^+ W) = \frac12 \tr(Z W^+),
\qquad Z,W \in \HC,
\end{equation}
where $\tr Z = 2z^0 = Z + Z^+$.
When this quadratic form is restricted to $\BB H$, $\HR$ and $\BB M$,
it has signature $(4,0)$, $(2,2)$ and $(3,1)$ respectively.

We will use the standard matrix realization of $\BB H$ so that
$$
e_0 = \begin{pmatrix} 1 & 0 \\ 0 & 1 \end{pmatrix}, \qquad
e_1 = \begin{pmatrix} 0 & -i \\ -i & 0 \end{pmatrix}, \qquad
e_2 = \begin{pmatrix} 0 & -1 \\ 1 & 0 \end{pmatrix}, \qquad
e_3 = \begin{pmatrix} -i & 0 \\ 0 & i \end{pmatrix}
$$
and
$$
\BB H = \{ Z \in \HC ;\: Z^c = Z \}
= \biggl\{
Z= \begin{pmatrix} z_{11} & z_{12} \\ z_{21} & z_{22} \end{pmatrix} \in \HC
; \: z_{22} = \overline{z_{11}}, \: z_{21} = - \overline{z_{12}} \biggr\}.
$$
Then $\HC$ can be identified with the algebra of all complex
$2 \times 2$ matrices:
$$
\HC = \biggl\{
Z= \begin{pmatrix} z_{11} & z_{12} \\ z_{21} & z_{22} \end{pmatrix}
; \: z_{ij} \in \BB C \biggr\},
$$
the quadratic form $N(Z)$ becomes $\det Z$, and the involution
$Z \mapsto Z^-$ becomes
$$
Z = \begin{pmatrix} z_{11} & z_{12} \\ z_{21} & z_{22} \end{pmatrix}
\quad \mapsto \quad Z^- =
\begin{pmatrix} 1 & 0 \\ 0 & -1 \end{pmatrix} Z
\begin{pmatrix} 1 & 0 \\ 0 & -1 \end{pmatrix}
= \begin{pmatrix} z_{11} & -z_{12} \\ -z_{21} & z_{22} \end{pmatrix}.
$$
The split quaternions $\HR$ and the Minkowski space $\BB M$ have matrix
realizations
\begin{align*}
\HR &= \biggl\{
Z= \begin{pmatrix} z_{11} & z_{12} \\ z_{21} & z_{22} \end{pmatrix} \in \HC
; \: z_{22} = \overline{z_{11}}, \: z_{21} = \overline{z_{12}} \biggr\},  \\
\BB M &= \biggl\{
Z= \begin{pmatrix} z_{11} & z_{12} \\ z_{21} & z_{22} \end{pmatrix} \in \HC
; \: z_{11}, z_{22} \in i \BB R, \: z_{21} = -\overline{z_{12}} \biggr\}.
\end{align*}
From this realization it is easy to see that the split quaternions form an
algebra over $\BB R$ isomorphic to $\mathfrak{gl} (2, \BB R)$ and
the invertible elements in $\HR$, denoted by $\HR^{\times}$,
are nothing else but $GL(2,\BB R)$.
We regard $SL(2,\BB C)$ as a quadric $\{N(Z) = 1\}$ in $\HC$, and
we also regard $SU(1,1) \simeq SL(2,\BB R)$ as the set of real points of
this quadric:
$$
SU(1,1) = \{ Z \in \HR ; \: N(Z) = 1 \}
= \biggl\{ Z = \begin{pmatrix} z_{11} & z_{12} \\
\overline{z_{12}} & \overline{z_{11}} \end{pmatrix}
\in \HR ;\: \det Z = |z_{11}|^2 - |z_{12}|^2 = 1 \biggr\}.
$$
The group
\begin{equation}  \label{SL(2,C)-embedding}
\biggl\{ \begin{pmatrix} a & 0 \\ 0 & a^c \end{pmatrix};\:
a \in \HC, \: N(a)=1 \biggr\} \quad \subset \quad GL(2,\HC)
\end{equation}
is naturally isomorphic to $SL(2,\BB C)$ and acts on the unit hyperboloid
of one sheet $\tilde H = \{ Y \in \BB M ;\: N(Y)=1 \}$ transitively.
The stabilizer group of $e_3 \in \tilde H$ is
$$
\biggl\{ \begin{pmatrix} a & 0 \\ 0 & a^c \end{pmatrix};\:
a \in \HR, \: N(a)=1 \biggr\} \quad \simeq \quad SU(1,1).
$$
Thus the hyperboloid $\tilde H$ -- also known as the imaginary Lobachevski
space -- is naturally identified with the homogeneous space
$SL(2,\BB C)/SU(1,1)$.

The algebra of split quaternions $\HR$ is spanned over $\BB R$
by the four matrices
$$
e_0 = \begin{pmatrix} 1 & 0 \\ 0 & 1 \end{pmatrix}, \qquad
\tilde e_1 = \begin{pmatrix} 0 & 1 \\ 1 & 0 \end{pmatrix}, \qquad
\tilde e_2 = \begin{pmatrix} 0 & i \\ -i & 0 \end{pmatrix}, \qquad
e_3 = \begin{pmatrix} -i & 0 \\ 0 & i \end{pmatrix},
$$
so
$$
\HR = \biggl\{ x^0 e_0 + x^1 \tilde e_1 + x^2 \tilde e_2 + x^3 e_3 =
\begin{pmatrix} x^0-ix^3 & x^1+ix^2 \\ x^1-ix^2 & x^0+ix^3 \end{pmatrix}
;\: x^0, x^1, x^2, x^3 \in \BB R \biggr\}.
$$
The quaternionic conjugation in this basis becomes
$$
e_0^+ = e_0, \quad \tilde e_1^+ = - \tilde e_1,
\quad \tilde e_2^+ = - \tilde e_2, \quad e_3^+ = -e_3.
$$
The multiplication rules for $\HR$ are:
\begin{center}
$e_0$ commutes with all elements of $\HR$, \\
$\tilde e_1$, $\tilde e_2$, $e_3$ anti-commute, \\
$e_0^2 = \tilde e_1^2 = \tilde e_2^2 = e_0$, \quad $e_3^2 = -e_0$,  \\
$\tilde e_1 \tilde e_2 = e_3$, \quad $\tilde e_2e_3 = -\tilde e_1$,
\quad $e_3 \tilde e_1= - \tilde e_2$.
\end{center}
The elements $e_0$, $\tilde e_1$, $\tilde e_2$, $e_3$
are orthogonal with respect to the bilinear form (\ref{bilinear_form})
and $\langle e_0, e_0 \rangle = \langle e_3, e_3 \rangle = 1$,
$\langle \tilde e_1, \tilde e_1 \rangle =
\langle \tilde e_2, \tilde e_2 \rangle = -1$.

Similarly, the Minkowski space $\BB M$ is a subspace of $\HC$
spanned over $\BB R$ by the four matrices
$$
\tilde e_0 = -ie_0 = \begin{pmatrix} -i & 0 \\ 0 & -i \end{pmatrix}, \qquad
e_1 = \begin{pmatrix} 0 & -i \\ -i & 0 \end{pmatrix}, \qquad
e_2 = \begin{pmatrix} 0 & -1 \\ 1 & 0 \end{pmatrix}, \qquad
e_3 = \begin{pmatrix} -i & 0 \\ 0 & i \end{pmatrix}.
$$

Define a norm on $\HC$ by
$$
\|Z\| = \frac 1{\sqrt 2}
\sqrt{|z_{11}|^2 + |z_{12}|^2 + |z_{21}|^2 + |z_{22}|^2}, \qquad
Z = \begin{pmatrix} z_{11} & z_{12} \\ z_{21} & z_{22} \end{pmatrix} \in \HC,
$$
so that $\|e_i\|=1$, $0 \le i \le 3$.

The (classical) quaternions $\BB H$ are oriented so that
$\{e_0, e_1, e_2, e_3 \}$ is a positive basis.
Let $dV= dz^0 \wedge dz^1 \wedge dz^2 \wedge dz^3$ be the holomorphic 4-form
on $\HC$ determined by $dV(e_0,e_1,e_2,e_3)=1$, then the restriction
$dV \bigl |_{\BB H}$ is the Euclidean volume form corresponding to
$\{e_0, e_1, e_2, e_3 \}$.
On the other hand, the restriction $dV \bigl |_{\HR}$ is also real-valued
and hence determines an orientation on $\HR$ so that
$\{e_0, \tilde e_1, \tilde e_2, e_3 \}$ becomes a positively oriented basis.
The orientation on $\BB M$ was defined in \cite{FL} so that
$\{\tilde e_0, e_1, e_2, e_3 \}$ is a positively oriented basis.

In \cite{FL} we defined a holomorphic 3-form on $\HC$ with values in $\HC$
$$
Dz = e_0 dz^1 \wedge dz^2 \wedge dz^3 - e_1 dz^0 \wedge dz^2 \wedge dz^3
+ e_2 dz^0 \wedge dz^1 \wedge dz^3 - e_3 dz^0 \wedge dz^1 \wedge dz^2
$$
characterized by the property
$$
\langle Z_1, Dz(Z_2,Z_3,Z_4) \rangle =
\frac 12 \tr (Z_1^+, Dz (Z_2,Z_3,Z_4)) = dV(Z_1,Z_2,Z_3,Z_4)
$$
for all $Z_1,Z_2,Z_3,Z_4 \in \HC$.
Let $Dx = Dz \bigl |_{\HR}$ and $D \tilde x = Dz \bigl |_{\BB H}$.

\begin{prop}
The 3-form $Dx$ takes values in $\HR$.
If we write $X = x^0 e_0 + x^1 \tilde e_1 + x^2 \tilde e_2 + x^3 e_3 \in \HR$,
$x^0, x^1, x^2, x^3 \in \BB R$, then $Dx$ is given explicitly by
\begin{equation}  \label{Dx-explicit}
Dx = e_0 dx^1 \wedge dx^2 \wedge dx^3 + \tilde e_1 dx^0 \wedge dx^2 \wedge dx^3
- \tilde e_2 dx^0 \wedge dx^1 \wedge dx^3 - e_3 dx^0 \wedge dx^1 \wedge dx^2.
\end{equation}
\end{prop}

\begin{rem}
Clearly, the form $Dx$ satisfies the property
$$
\langle X_1, Dx(X_2,X_3,X_4) \rangle =
\frac 12 \tr (X_1^+, Dx (X_2,X_3,X_4)) = dV(X_1,X_2,X_3,X_4)
$$
for all $X_1,X_2,X_3,X_4 \in \HR$, which could be used to define it.
\end{rem}

Let $U \subset \HR$ be an open region with piecewise smooth boundary
$\partial U$. We give a canonical orientation to $\partial U$ as follows.
The positive orientation of $U$ is determined by
$\{e_0, \tilde e_1, \tilde e_2, e_3 \}$.
Pick a smooth point $p \in \partial U$ and let $\overrightarrow{n_p}$
be a non-zero vector in $T_p\HR$ perpendicular to $T_p\partial U$ and
pointing outside of $U$.
Then $\{\overrightarrow{\tau_1}, \overrightarrow{\tau_2},
\overrightarrow{\tau_3}\} \subset T_p \partial U$ is positively oriented
in $\partial U$ if and only if
$\{\overrightarrow{n_p}, \overrightarrow{\tau_1}, \overrightarrow{\tau_2},
\overrightarrow{\tau_3}\}$ is positively oriented in $\HR$.
We orient $SU(1,1)$ and, more generally, hyperboloids
$H_R = \{X \in \HR ;\: N(X)=R^2 \}$ as the boundaries of the open sets
$\{X \in \HR ;\: N(X)<R^2 \}$.

\begin{lem}  \label{restrictions}
Let $R \in \BB R$ be a constant, then we have the following restriction
formula:
$$
Dx \bigl|_{H_R} = X \,\frac{dS}{\|X\|},
$$
where $dS$ denotes the restrictions of the Euclidean measure on $\HR$ to $H_R$.
\end{lem}

We recall some notations from \cite{FL}.
Let $\BB S$ be the natural two-dimensional complex representation of
the algebra $\HC$.
We realize it as a column of two complex numbers, then $\HC$ acts on $\BB S$
by matrix multiplication on the left.
Similarly, we denote by $\BB S'$ the dual space of $\BB S$, this time
realized as a row of two complex numbers. We have a right action of $\HC$ on
$\BB S'$ by multiplication on the right.

\subsection{Regular Functions on $\BB H$ and $\HC$}

Recall that regular functions on $\BB H$ are defined using an analogue of the
Cauchy-Riemann equations. Write $\tilde X \in \BB H$ as
$\tilde X = \tilde x^0 e_0 + \tilde x^1 e_1 + \tilde x^2 e_2 + \tilde x^3 e_3$,
$\tilde x^0, \tilde x^1, \tilde x^2, \tilde x^3 \in \BB R$,
and factor the four-dimensional Laplacian operator $\square$ on $\BB H$
as a product of two Dirac operators
$$
\square = \frac {\partial^2}{(\partial \tilde x^0)^2} +
\frac {\partial^2}{(\partial \tilde x^1)^2} +
\frac {\partial^2}{(\partial \tilde x^2)^2} +
\frac {\partial^2}{(\partial \tilde x^3)^2}
= \nabla \nabla^+ = \nabla^+ \nabla,
$$
where
\begin{align*}
\nabla^+ &= e_0 \frac{\partial}{\partial \tilde x^0}
+ e_1 \frac{\partial}{\partial \tilde x^1}
+ e_2 \frac{\partial}{\partial \tilde x^2}
+ e_3 \frac{\partial}{\partial \tilde x^3} \qquad \text{and} \\
\nabla &= e_0 \frac{\partial}{\partial \tilde x^0}
- e_1 \frac{\partial}{\partial \tilde x^1}
- e_2 \frac{\partial}{\partial \tilde x^2}
- e_3 \frac{\partial}{\partial \tilde x^3}.
\end{align*}
The operators $\nabla^+$, $\nabla$ can be applied to functions
on the left and on the right.
For an open subset $U \subset \BB H$ and a differentiable function $f$ on $U$
with values in $\BB H$, $\BB S$ or $\HC$, we say $f$ is {\em left-regular} if
$(\nabla^+ f)(\tilde X)=0$ for all $\tilde X \in U$.
Similarly, a differentiable function $g$ on $U$ with values in $\BB H$,
$\BB S'$ or $\HC$, is {\em right-regular} if $(g \nabla^+)(\tilde X)=0$
for all $\tilde X \in U$.

\begin{prop}
For any ${\cal C}^1$-functions $f: U \to \BB S$ and $g: U \to \BB S'$,
where $U \subset \BB H$ is an open subset,
$$
d(D\tilde x \cdot f) =
- D\tilde x \wedge df = (\nabla^+ f) \,dV \bigr|_{\BB H}, \qquad
d(g \cdot D\tilde x) =
dg \wedge D\tilde x = (g \nabla^+) \,dV \bigr|_{\BB H}.
$$
In particular,
$$
\nabla^+ f = 0 \quad \Longleftrightarrow \quad
D\tilde x \wedge df =0, \qquad
g \nabla^+ = 0 \quad \Longleftrightarrow \quad
dg \wedge D\tilde x =0.
$$
\end{prop}

Let $U^{\BB C} \subset \HC$ be an open set.
Following \cite{FL}, we say that a differential function on $U^{\BB C}$
with values in $\BB C$, $\HC$, $\BB S$ or $\BB S'$ is {\em holomorphic}
if it is holomorphic with respect to the complex variables $z^0, z^1, z^2, z^3$.
Then we define $f^{\BB C}: U^{\BB C} \to \BB S$ to be
{\em holomorphic left-regular} if it is holomorphic and $\nabla^+ f^{\BB C} =0$.
Similarly, $g^{\BB C}: U^{\BB C} \to \BB S'$ is defined to be
{\em holomorphic right-regular} if it is holomorphic and $g^{\BB C} \nabla^+ =0$.
Finally, we call a function $\phi: U^{\BB C} \to \BB C$
{\em holomorphic harmonic} if it is holomorphic and $\square \phi=0$.

If we identify $\HC$ with complex $2 \times 2$ matrices
$\begin{pmatrix} z_{11} & z_{12} \\ z_{21} & z_{22} \end{pmatrix}$,
$z_{ij} \in \BB C$, then a function is holomorphic if and only if it is
holomorphic with respect to the complex variables $z_{ij}$, $1 \le i,j \le 2$.
For holomorphic functions $f^{\BB C}: U^{\BB C} \to \BB S$ and
$g^{\BB C}: U^{\BB C} \to \BB S'$, the following derivatives are equal:
$$
\nabla^+ f^{\BB C} =
e_0 \frac{\partial f^{\BB C}}{\partial z^0} 
+ e_1 \frac{\partial f^{\BB C}}{\partial z^1}
+ e_2 \frac{\partial f^{\BB C}}{\partial z^2}
+ e_3 \frac{\partial f^{\BB C}}{\partial z^3}
=
2 \begin{pmatrix} \frac {\partial}{\partial z_{22}} &
- \frac {\partial}{\partial z_{21}}  \\
- \frac {\partial}{\partial z_{12}} &
\frac {\partial}{\partial z_{11}} \end{pmatrix} f^{\BB C},
$$
$$
g^{\BB C} \nabla^+ =
\frac{\partial g^{\BB C}}{\partial z^0} e_0
+ \frac{\partial g^{\BB C}}{\partial z^1} e_1
+ \frac{\partial g^{\BB C}}{\partial z^2} e_2
+ \frac{\partial g^{\BB C}}{\partial z^3} e_3
=
2 g^{\BB C} \begin{pmatrix} \frac {\partial}{\partial z_{22}} &
- \frac {\partial}{\partial z_{21}}  \\
- \frac {\partial}{\partial z_{12}} &
\frac {\partial}{\partial z_{11}} \end{pmatrix},
$$
$$
\nabla f^{\BB C} =
e_0 \frac{\partial f^{\BB C}}{\partial z^0} 
- e_1 \frac{\partial f^{\BB C}}{\partial z^1}
- e_2 \frac{\partial f^{\BB C}}{\partial z^2}
- e_3 \frac{\partial f^{\BB C}}{\partial z^3}
=
2 \begin{pmatrix} \frac {\partial}{\partial z_{11}} &
\frac {\partial}{\partial z_{21}}  \\
\frac {\partial}{\partial z_{12}} &
\frac {\partial}{\partial z_{22}} \end{pmatrix} f^{\BB C},
$$
$$
g^{\BB C} \nabla =
\frac{\partial g^{\BB C}}{\partial z^0} e_0
- \frac{\partial g^{\BB C}}{\partial z^1} e_1
- \frac{\partial g^{\BB C}}{\partial z^2} e_2
- \frac{\partial g^{\BB C}}{\partial z^3} e_3
=
2 g^{\BB C} \begin{pmatrix} \frac {\partial}{\partial z_{11}} &
\frac {\partial}{\partial z_{21}}  \\
\frac {\partial}{\partial z_{12}} &
\frac {\partial}{\partial z_{22}} \end{pmatrix}.
$$
Since we are interested in holomorphic functions only, we will abuse the
notation and denote by $\nabla$ and $\nabla^+$ the holomorphic differential
operators
$\frac{\partial}{\partial z^0} e_0
- \frac{\partial}{\partial z^1} e_1
- \frac{\partial}{\partial z^2} e_2
- \frac{\partial}{\partial z^3} e_3$
and
$e_0 \frac{\partial}{\partial z^0} 
+ e_1 \frac{\partial}{\partial z^1}
+ e_2 \frac{\partial}{\partial z^2}
+ e_3 \frac{\partial}{\partial z^3}$
respectively.

\begin{prop}
For any holomorphic functions $f^{\BB C}: U^{\BB C} \to \BB S$ and
$g^{\BB C}: U^{\BB C} \to \BB S'$,
$$
\nabla^+ f^{\BB C} = 0 \quad \Longleftrightarrow \quad Dz \wedge df^{\BB C} =0,
\qquad
g^{\BB C} \nabla^+= 0 \quad \Longleftrightarrow \quad dg^{\BB C} \wedge Dz =0.
$$
\end{prop}

\begin{lem}  \label{closed}
We have:
\begin{enumerate}
\item
$\square \frac 1{N(Z)} = 0$;
\item
$\nabla \frac 1{N(Z)}
= \frac 1{N(Z)} \nabla
= -2 \frac {Z^{-1}}{N(Z)} = -2 \frac {Z^+} {N(Z)^2}$;
\item
$\frac {Z^{-1}}{N(Z)} = \frac {Z^+} {N(Z)^2}$
is a holomorphic left- and right-regular function defined wherever
$N(Z) \ne 0$;
\item
The form
$\frac {Z^{-1}}{N(Z)} \cdot Dz = \frac {Z^+} {N(Z)^2} \cdot Dz$ is a closed
holomorphic $\HC$-valued 3-form defined wherever $N(Z) \ne 0$.
\end{enumerate}
\end{lem}

\begin{lem}  \label{inverse-deriv}
For any differentiable function $\phi: U^{\BB C} \to \BB C$, we have:
$$
\nabla \bigl( \phi(Z^+) \bigr) = (\nabla^+ \phi)(Z^+),
\qquad
\nabla^+ \bigl( \phi(Z^+) \bigr) = (\nabla \phi)(Z^+),
$$
$$
\nabla \bigl( \phi(Z^{-1}) \bigr) =
- Z^{-1} \cdot (\nabla \phi)(Z^{-1}) \cdot Z^{-1}.
$$
\end{lem}

We will often use the shifted degree operator
$$
\deg = 1 + z_0 \frac{\partial}{\partial z^0}
+ z_1 \frac{\partial}{\partial z^1}
+ z_2 \frac{\partial}{\partial z^2}
+ z_3 \frac{\partial}{\partial z^3} 
$$
(the degree operator plus identity). By direct computation we obtain

\begin{lem}
\begin{equation}  \label{deg-nabla}
2(\deg + 1) = Z^+\nabla^+ + \nabla Z = \nabla^+Z^+ + Z\nabla.
\end{equation}
\end{lem}

Next we describe actions of the group $GL(2,\HC)$ on the spaces of
left-, right-regular and harmonic functions on $\HC$ with
singularities by conformal (fractional linear) transformations. Let
$$
h = \begin{pmatrix} a' & b' \\ c' & d' \end{pmatrix} \in GL(2,\HC)
\qquad \text{and write} \qquad
h^{-1} = \begin{pmatrix} a & b \\ c & d \end{pmatrix}.
$$
On the space of left-regular $\BB S$-valued functions $GL(2,\HC)$ acts by
$$
\pi_l(h): \: f(Z) \quad \mapsto \quad (\pi_l(h)f)(Z) =
\frac {(cZ+d)^{-1}}{N(cZ+d)} \cdot f \bigl( (aZ+b)(cZ+d)^{-1} \bigr).
$$
On the space of right-regular $\BB S'$-valued functions $GL(2,\HC)$ acts by
$$
\pi_r(h): \: g(Z) \quad \mapsto \quad (\pi_r(h)g)(Z) =
g \bigl( (a'-Zc')^{-1}(-b'+Zd') \bigr) \cdot \frac {(a'-Zc')^{-1}}{N(a'-Zc')}.
$$
On the space of $\BB C$-valued harmonic functions we have two different actions:
\begin{align*}
\pi^0_l(h): \: \phi(Z) \quad &\mapsto \quad \bigl( \pi^0_l(h)\phi \bigr)(Z) =
\frac 1{N(cZ+d)} \cdot \phi \bigl( (aZ+b)(cZ+d)^{-1} \bigr),  \\
\pi^0_r(h): \: \phi(Z) \quad &\mapsto \quad \bigl( \pi^0_r(h)\phi \bigr)(Z) =
\frac 1{N(a'-Zc')} \cdot \phi \bigl( (a'-Zc')^{-1}(-b'+Zd') \bigr).
\end{align*}
These two actions coincide on $SL(2,\HC)$, which is defined
as the connected Lie subgroup of $GL(2,\HC)$ with Lie algebra
$$
\mathfrak{sl}(2,\HC) = \{ x \in \mathfrak{gl}(2,\HC) ;\: \re (\tr x) =0 \}.
$$

\subsection{Regular Functions on $\HR$}

We introduce linear differential operators on $\HR$
\begin{align*}
\nabla^+_{\BB R} &= e_0 \frac{\partial}{\partial x^0}
- \tilde e_1 \frac{\partial}{\partial x^1}
- \tilde e_2 \frac{\partial}{\partial x^2}
+ e_3 \frac{\partial}{\partial x^3} \qquad \text{and} \\
\nabla_{\BB R} &= e_0 \frac{\partial}{\partial x^0}
+ \tilde e_1 \frac{\partial}{\partial x^1}
+ \tilde e_2 \frac{\partial}{\partial x^2}
- e_3 \frac{\partial}{\partial x^3}
\end{align*}
which may be applied to functions on the left and on the right.

\begin{df}
Fix an open subset $U \subset \HR$. A differentiable function
$f: \HR \to \BB S$ is {\em left-regular} if it satisfies
$$
(\nabla^+_{\BB R} f)(X) = e_0 \frac{\partial f}{\partial x^0}(X)
- \tilde e_1 \frac{\partial f}{\partial x^1}(X)
- \tilde e_2 \frac{\partial f}{\partial x^2}(X)
+ e_3 \frac{\partial f}{\partial x^3}(X) =0, \qquad \forall X \in U.
$$
Similarly, a differentiable function $g: \HR \to \BB S'$ is
{\em right-regular} if
$$
(g\nabla^+_{\BB R} )(X) = \frac{\partial g}{\partial x^0}(X)e_0
- \frac{\partial g}{\partial x^1}(X) \tilde e_1
- \frac{\partial g}{\partial x^2}(X) \tilde e_2
+ \frac{\partial g}{\partial x^3}(X)e_3 =0, \qquad \forall X \in U.
$$
\end{df}

We denote by $\square_{2,2}$ the ultrahyperbolic wave operator on $\HR$
which can be factored as follows:
$$
\square_{2,2} = \frac {\partial^2}{(\partial x^0)^2} -
\frac {\partial^2}{(\partial x^1)^2} -
\frac {\partial^2}{(\partial x^2)^2} +
\frac {\partial^2}{(\partial x^3)^2}
=\nabla_{\BB R} \nabla^+_{\BB R} = \nabla^+_{\BB R} \nabla_{\BB R}.
$$ 
Abusing terminology, we call a smooth function $\phi: U \to \BB C$
{\em harmonic} if $\square_{2,2} \phi=0$.

\begin{prop}  \label{d(fDzg)}
For any ${\cal C}^1$-functions $f: U \to \BB S$ and $g: U \to \BB S'$,
$$
d(Dx \cdot f) = - Dx \wedge df = (\nabla^+_{\BB R} f) \,dV, \qquad
d(g \cdot Dx) = dg \wedge Dx = (g \nabla^+_{\BB R}) \,dV,
$$
In particular, 
$$
\nabla^+_{\BB R} f = 0 \quad \Longleftrightarrow \quad Dx \wedge df =0,
\qquad
g \nabla^+_{\BB R} = 0 \quad \Longleftrightarrow \quad dg \wedge Dx =0.
$$
\end{prop}

Let $U^{\BB C} \subset \HC$ be an open set. The restriction relations
$$
Dz \bigl |_{\HR} = Dx,
\qquad
Dz \bigl |_{\BB H} = D\tilde x
$$
imply that the restriction of a holomorphic left- or right-regular function
to $U^{\BB R} = U^{\BB C} \cap \HR$ produces a left- or right-regular function
on $U^{\BB R}$ respectively.
And the restriction of a holomorphic left- or right-regular function
to $U_{\BB H} = U^{\BB C} \cap \BB H$ also produces a left- or right-regular
function on $U_{\BB H}$ respectively.
Conversely, if one starts with, say, a left-regular function on $\HR$,
extends it holomorphically to a left-regular function on $\HC$
and then restricts the extension to $\BB H$, the resulting function
is left-regular on $\BB H$.

\begin{prop}
Let $f^{\BB C}: U^{\BB C} \to \BB S$ and $g^{\BB C}: U^{\BB C} \to \BB S'$ be
holomorphic functions. Then
$$
\nabla^+_{\BB R} f^{\BB C} = \nabla^+ f^{\BB C}, \qquad
\nabla_{\BB R} f^{\BB C} = \nabla f^{\BB C}, \qquad
g^{\BB C} \nabla^+_{\BB R} = g^{\BB C} \nabla^+, \qquad
g^{\BB C} \nabla_{\BB R} = g^{\BB C} \nabla.
$$
\end{prop}

The actions of the group $GL(2,\HR)$ on the spaces of
left-, right-regular and harmonic functions on $\HR$ with
singularities are given by the same formulas as before.

\subsection{Fueter Formula for Holomorphic Regular Functions on $\HR$}

We are interested in extensions of the Cauchy-Fueter formula to functions
on $\HR$.
First we recall the classical version of the integral formula due to Fueter:

\begin{thm} [Cauchy-Fueter Formula \cite{F1, F2}]  \label{Fueter}
Let $U_{\BB H} \subset \BB H$ be an open bounded subset with piecewise
${\cal C}^1$ boundary $\partial U_{\BB H}$.
Suppose that $f(\tilde X)$ is left-regular on a
neighborhood of the closure $\overline{U_{\BB H}}$, then
$$
\frac 1 {2\pi^2} \int_{\partial U_{\BB H}}
\frac {(\tilde X - \tilde X_0)^{-1}}{N(\tilde X - \tilde X_0)}
\cdot D\tilde x \cdot f(\tilde X) =
\begin{cases}
f(\tilde X_0) & \text{if $\tilde X_0 \in U_{\BB H}$;}\\
0 & \text{if $\tilde X_0 \notin \overline{U_{\BB H}}$.}
\end{cases}
$$
If $g(\tilde X)$ is right-regular on a neighborhood of the
closure $\overline{U_{\BB H}}$, then
$$
\frac 1 {2\pi^2} \int_{\partial U_{\BB H}}
g(\tilde X) \cdot D \tilde x \cdot
\frac {(\tilde X - \tilde X_0)^{-1}}{N(\tilde X - \tilde X_0)} =
\begin{cases}
g(\tilde X_0) & \text{if $\tilde X_0 \in U_{\BB H}$;}\\
0 & \text{if $\tilde X_0 \notin \overline{U_{\BB H}}$.}
\end{cases}
$$
\end{thm}

Let $U \subset \HR$ be an open subset, and let $f$ be a ${\cal C}^1$-function
defined on a neighborhood of $\overline{U}$ such that $\nabla_{\BB R}^+ f =0$.
In this subsection we extend the Cauchy-Fueter integral formula to left-regular
functions which can be extended holomorphically to a
neighborhood of $\overline{U}$ in $\HC$.
Observe that the expression in the integral formula
$\frac {(\tilde X - \tilde X_0)^{-1}}{N(\tilde X - \tilde X_0)}
\cdot D\tilde x$ is nothing else but the restriction to $\BB H$ of the
holomorphic 3-form $\frac {(Z-\tilde X_0)^{-1}} {N(Z-\tilde X_0)} \cdot Dz$,
which is the form from Lemma \ref{closed} translated by $\tilde X_0$.
For this reason we expect an integral formula of the kind
$$
f(X_0) = \frac 1 {2\pi^2} \int_{\partial U}
\biggl( \frac {(Z-X_0)^{-1}} {N(Z-X_0)} \cdot Dz \biggr)
\biggl|_{\HR} \cdot f(X),
\qquad \forall X_0 \in U.
$$
However, the integrand is singular wherever $N(Z-X_0)=0$.
We resolve this difficulty by deforming the contour of integration $\partial U$
in such a way that the integral is no longer singular.
Fix an $\epsilon \in \BB R$ and define an $\epsilon$-deformation
$h_{\epsilon}: \HC \to \HC$, $Z \mapsto Z_{\epsilon}$, by:
\begin{center}
\begin{tabular} {lcl}
$z_{11} \quad \mapsto \quad z_{11} + i\epsilon z_{11}$ & \qquad \qquad &
$z_{12} \quad \mapsto \quad z_{12} - i\epsilon z_{12}$  \\
$z_{21} \quad \mapsto \quad z_{21} - i\epsilon z_{21}$ & \qquad \qquad &
$z_{22} \quad \mapsto \quad z_{22} + i\epsilon z_{22}$.
\end{tabular}
\end{center}
In other words, $Z_{\epsilon} = Z + i \epsilon Z^-$.
For $Z_0 \in \HC$ fixed, we use a notation
$$
h_{\epsilon, Z_0} (Z) = Z_0 + h_{\epsilon}(Z-Z_0) = Z + i\epsilon (Z-Z_0)^-.
$$

\begin{lem}
Define a quadratic form $S(Z) = z_{11}z_{22} + z_{12}z_{21}$ on $\HC$. We have:
$$
N(Z_{\epsilon}) = (1-\epsilon^2) N(Z) + i 2\epsilon S(Z),
$$
$$
S(X) = \|X\|^2, \qquad \forall X \in \HR.
$$
\end{lem}

\begin{thm}  \label{holomorphic_Fueter}
Let $U \subset \HR$ be an open bounded subset with piecewise
${\cal C}^1$ boundary $\partial U$, and let $f(X)$ be a function defined on a
neighborhood of the closure $\overline{U}$ such that $\nabla^+_{\BB R} f =0$.
Suppose that $f$ extends to a holomorphic left-regular function on an open
subset $V^{\BB C} \subset \HC$ containing $\overline{U}$, then
$$
- \frac 1 {2\pi^2} \int_{(h_{\epsilon, X_0})_*(\partial U)}
\frac {(Z-X_0)^{-1}} {N(Z-X_0)} \cdot Dz \cdot f^{\BB C}(Z) =
\begin{cases}
f(X_0) & \text{if $X_0 \in U$;}\\
0 & \text{if $X_0 \notin \overline{U}$.}
\end{cases}
$$
for all $\epsilon \ne 0$ sufficiently close to 0.
\end{thm}

\begin{rem}
For all $\epsilon \ne 0$ sufficiently close to 0 the contour of integration
$(h_{\epsilon, X_0})_*(\partial U)$ lies inside $V^{\BB C}$ and the integrand
is non-singular, thus the integrals are well defined.
Moreover, the value of the integral becomes constant when
the parameter $\epsilon$ is sufficiently close to 0.
Of course, there is a similar formula for right-regular functions on $\HR$.
\end{rem}

The proof is similar to the proof of Theorem 51 in \cite{FL};
for this reason we just give an outline.
Since the integrand is a closed form, by Stokes' we are free to deform
the cycle of integration as long as we stay inside the set
\begin{equation}  \label{VminusN}
\{ Z \in V^{\BB C} ;\: N(X_0-Z) \ne 0 \}.
\end{equation}
Let $S_r = \{ \tilde X \in \BB H+X_0 ;\: \|\tilde X-X_0\|^2 = r^2 \}$ be
the sphere of radius $r$ centered at $X_0$ and lying in the subspace of $\HC$
parallel to $\BB H$. We orient $S_r$ as the boundary of the open ball.
One can show that, for $r>0$ sufficiently small, the cycle of integration
$(h_{\epsilon, X_0})_*(\partial U)$ is homologous to $-S_r$ if $X_0 \in U$
and $0$ if $X_0 \notin \overline{U}$ as homology 3-cycles in (\ref{VminusN}).
Then the result follows from the Fueter formula for the regular quaternions
(Theorem \ref{Fueter}). Alternatively, one can let $r \to 0^+$ and show
directly that the integral remains unchanged and at the same time approaches
$-2\pi^2 f(X_0)$ in the same way the Cauchy and Cauchy-Fueter formulas are
proved.

One can drop the assumption that $f(X)$ extends holomorphically to an
open neighborhood of $\B{U}$ in $\HC$ and prove the following
version involving generalized functions:

\begin{thm}[Integral Formula]
Let $U \subset \HR$ be a bounded open region with smooth boundary
$\partial U$.
Let $f: U \to \HC$ be a function which extends to a real-differentiable
function on an open neighborhood $V \subset \HR$ of the closure
$\overline{U}$ such that $\nabla^+_{\BB R} f = 0$.
Then, for any point $X_0 \in \HR$ such that $\partial U$ intersects the cone
$\{ X \in \HR ;\: N(X-X_0) =0 \}$ transversally, we have:
$$
\lim_{\epsilon \to 0} -\frac1{2\pi^2} \int_{X \in \partial U}
\frac {(X-X_0)^+} {\bigl( N(X-X_0) +i\epsilon \|X-X_0\|^2 \bigr)^2}
\cdot Dz \cdot f(X) =
\begin{cases}
f(X_0) & \text{if $X_0 \in U$;}\\
0 & \text{if $X_0 \notin \overline{U}$.}
\end{cases}
$$
\end{thm}

The proof of this theorem is given in \cite{L}.

\subsection{The Matrix Coefficients of $SU(1,1)$}

The matrix coefficients of the (generalized) principal series representations
are functions of $Z \in SU(1,1)$ given by
\begin{equation}  \label{int_t}
t^l_{n\,\underline{m}}(Z) = \frac 1{2\pi i}
\oint (sz_{11}+z_{21})^{l-m} (sz_{12}+z_{22})^{l+m} s^{-l+n} \frac {ds}s,
\qquad Z = \begin{pmatrix} z_{11} & z_{12} \\ z_{21} & z_{22} \end{pmatrix},
\end{equation}
where the integral is taken over the unit circle $\{ s\in \BB C ;\: |s|=1 \}$
traversed once in the counterclockwise direction \cite{V}.
The parameters $m$ and $n$ are either both integers or half-integers:
$$
m,n \in \BB Z \qquad \text{or} \qquad m,n \in \BB Z + \frac 12,
\qquad \text{and} \qquad l \in \BB C.
$$
When the parameters $l,m,n$ range over
$$
l = -\frac 12, -1, -\frac 32, -2, -\frac 52, \dots,
\qquad m,n \in \BB Z + l, \qquad m, n \ge -l,
$$
we get the matrix coefficients of the holomorphic discrete series
representations ($l \le -1$) together with its limit ($l=-1/2$).
When the parameters $l,m,n$ range over
$$
l = -\frac 12, -1, -\frac 32, -2, -\frac 52, \dots,
\qquad m,n \in \BB Z + l, \qquad m, n \le l,
$$
we get the matrix coefficients of the antiholomorphic discrete series
representations ($l \le -1$) together with its limit ($l=-1/2$).
When the parameter $l$ ranges over
$$
l = - \frac 12 +i \lambda, \qquad \lambda \in \BB R,
$$
we get the matrix coefficients of the continuous series representations
and the two limits of the discrete series
($\lambda=0$ and $m,n \in \BB Z +\frac12$).
Formula (\ref{int_t}) involves complex numbers raised
to complex powers, and we need to clarify it so there is no ambiguity.
In the discrete series situation the powers $l \pm m$ are actually
integers, so there is no ambiguity at all.
In the general case $l \in \BB C$ we write
\begin{multline*}
(sz_{11}+z_{21})^{l-m} (sz_{12}+z_{22})^{l+m} s^{-l+n}
= (sz_{11}+z_{21})^{l-m} (z_{12}+ s^{-1} z_{22})^{l+m} s^{m+n}  \\
= |sz_{11}+z_{21}|^{2(l-m)} (z_{12}+ s^{-1} z_{22})^{2m} s^{m+n}
\end{multline*}
(recall that $\overline{z_{22}} = z_{11}$, $\overline{z_{12}} = z_{21}$,
$s^{-1} = \overline{s}$).
The expression $|sz_{11}+z_{21}|^{2(l-m)}$ is well defined because
$|sz_{11}+z_{21}|$ is a positive real number, and the expressions
$(z_{12}+ s^{-1} z_{22})^{2m}$, $s^{m+n}$ are well defined as well since
the powers $2m$, $m+n$ are integers.
For certain values of $l$, $m$, $n$, the matrix coefficients vanish:
\begin{equation}  \label{t=0}
t^l_{n\,\underline{m}}(Z) = 0
\qquad \text{when} \quad l=0, -\frac 12,-1,-\frac32, \dots
\quad \text{and} \quad
\begin{matrix}
m \ge -l, \: n<-l \\ \text{or} \\ m \le l, \: n>l.
\end{matrix}
\end{equation}

We give another expression for the matrix coefficients (\ref{int_t})
also due to Vilenkin \cite{V}.
Fix a parameterization of the group $SU(1,1)$ which is essentially the
$KAK$ decomposition:
\begin{multline}  \label{param}
X(\phi,\tau,\psi) =
\begin{pmatrix} x^0-ix^3 & x^1+ix^2 \\ x^1-ix^2 & x^0+ix^3 \end{pmatrix} =
\begin{pmatrix} \cosh \frac{\tau}2 \cdot e^{i\frac{\phi+\psi}2} &
\sinh \frac{\tau}2 \cdot e^{i\frac{\phi-\psi}2} \\
\sinh \frac{\tau}2 \cdot e^{i\frac{\psi-\phi}2} &
\cosh \frac{\tau}2 \cdot e^{-i\frac{\phi+\psi}2} \end{pmatrix}  \\
=
\begin{pmatrix} e^{i\frac{\phi}2} & 0 \\ 0 & e^{-i\frac{\phi}2} \end{pmatrix}
\begin{pmatrix} \cosh \frac{\tau}2 & \sinh \frac{\tau}2 \\
\sinh \frac{\tau}2 & \cosh \frac{\tau}2 \end{pmatrix}
\begin{pmatrix} e^{i\frac{\psi}2} & 0 \\ 0 & e^{-i\frac{\psi}2} \end{pmatrix}
\end{multline}
with $0 \le \phi < 2\pi$, $0 < \tau < \infty$, $-2\pi \le \psi < 2\pi$.
Then we have the following expressions for the Haar measure on $SU(1,1)$:
$$
\frac {dx^1 \wedge dx^2 \wedge dx^3}{x^0}
= \frac {\sinh \tau}8 \, d\phi \wedge d\tau \wedge d\psi
= Dx \cdot X^{-1} = \frac {dS}{\|X\|}.
$$
In terms of the parameterization (\ref{param}), the coefficients
$t^l_{n\,\underline{m}}(X)$, $X \in SU(1,1)$, can be written as
\begin{equation}  \label{t-german}
t^l_{n\,\underline{m}}(\phi, \tau, \psi) =
e^{-i(n\phi + m\psi)} \cdot \mathfrak{P}^l_{n\,m}(\cosh \tau),
\end{equation}
where, by definition,
\begin{equation*}
\mathfrak{P}^l_{n\,m}(\cosh \tau) =
\frac 1{2\pi i} \oint \Bigl( s\cosh \frac{\tau}2+\sinh \frac{\tau}2 \Bigr)^{l-m}
\Bigl( s\sinh \frac{\tau}2 + \cosh \frac{\tau}2 \Bigr)^{l+m} s^{-l+n} \frac {ds}s.
\end{equation*}
The functions $\mathfrak{P}^l_{n\,m}(\cosh \tau)$ are real-valued when
$l \in \BB R$. Moreover,
$$
\B{\mathfrak{P}^l_{n\,m}(\cosh \tau)} = \mathfrak{P}^{\B{l}}_{n\,m}(\cosh \tau),
$$
\begin{multline*}
\mathfrak{P}^l_{m\,n}(\cosh \tau)
= \mathfrak{P}^l_{-m\,-n}(\cosh \tau)
= (-1)^{m-n} \mathfrak{P}^{-l-1}_{n\,m}(\cosh \tau)  \\
= \frac{\Gamma(l+n+1) \Gamma(l-n+1)}{\Gamma(l+m+1) \Gamma(l-m+1)}
\cdot \mathfrak{P}^l_{n\,m}(\cosh \tau).
\end{multline*}
(When $l \in \frac 12 \BB Z$ two of the $\Gamma$-factors become
infinite and must be transformed by the formula
$\Gamma(x) \cdot \Gamma(1-x) = \pi/\sin(\pi x)$.)
They also satisfy the orthogonality relations:
\begin{equation}  \label{P-orthog}
\int_1^{\infty} \mathfrak{P}^l_{m\,n}(t) \cdot \B{\mathfrak{P}^{l'}_{m\,n}(t)} \,dt =0
\qquad \text{if $l \ne l'$}, \qquad l,l',m,n \in \frac 12 \BB Z,
\end{equation}
\begin{equation}  \label{P-norm}
\int_1^{\infty} \bigl| \mathfrak{P}^l_{n\,m}(t) \bigr|^2 \,dt
= (-1)^{m-n} \frac {-2}{2l+1}
\frac{\Gamma(l+m+1) \Gamma(l-m+1)}{\Gamma(l+n+1) \Gamma(l-n+1)},
\qquad l = -1,-\frac 32,-2, \dots.
\end{equation}
The functions $\mathfrak{P}^l_{m\,n}(\cosh t)$ can be expressed in terms of
the hypergeometric function. We will use
\begin{multline}  \label{P-hypergeometric}
\mathfrak{P}^l_{n\,m}(\cosh \tau)  \\
= \frac{\Gamma(l-m+1) \bigl(\cosh \frac{\tau}2 \bigr)^{m+n}
\bigl(\sinh \frac{\tau}2 \bigr)^{n-m}}{\Gamma(l-n+1)(n-m)!} \cdot
\,_2F_1\Bigl(l+n+1,-l+n;n-m+1; -\sinh^2 \frac{\tau}2 \Bigr)
\end{multline}
(valid for $n \ge m$) found in Subsection 6.5.3 of \cite{VK}.


Let $\HR^+ = \{ X \in \HR ;\: N(X) >0 \}$.
Clearly, the matrix coefficient functions $t^l_{n\,\underline{m}}$
extend from $SU(1,1)$ to $\HR^+$.
Differentiating (\ref{int_t}) under the integral sign we obtain:

\begin{lem}  \label{deriv_calc}
$$
\frac {\partial}{\partial z_{11}} t^l_{n\,\underline{m}}(Z) =
(l-m) t^{l- \frac 12}_{n+ \frac 12 \,\underline{m+ \frac 12}}(Z),\qquad
\frac {\partial}{\partial z_{12}} t^l_{n\,\underline{m}}(Z) =
(l+m) t^{l- \frac 12}_{n+ \frac 12 \,\underline{m- \frac 12}}(Z),
$$
$$
\frac {\partial}{\partial z_{21}} t^l_{n\,\underline{m}}(Z) =
(l-m) t^{l- \frac 12}_{n- \frac 12 \,\underline{m+ \frac 12}}(Z), \qquad
\frac {\partial}{\partial z_{22}} t^l_{n\,\underline{m}}(Z) =
(l+m) t^{l- \frac 12}_{n- \frac 12 \,\underline{m- \frac 12}}(Z).
$$
\end{lem}

We also have the following multiplication identities for matrix coefficients.

\begin{lem}  \label{mult-ids}
$$
\Bigl( t^{l-\frac 12}_{m + \frac 12 \, \underline{n}}(Z),
t^{l-\frac 12}_{m - \frac 12 \, \underline{n}}(Z) \Bigr)
\begin{pmatrix} z_{11} & z_{12} \\ z_{21} & z_{22} \end{pmatrix}
=
\Bigl( t^l_{m \, \underline{n- \frac 12}}(Z),
t^l_{m \, \underline{n+ \frac 12}}(Z) \Bigr)
$$
and
$$
\begin{pmatrix} z_{11} & z_{12} \\ z_{21} & z_{22} \end{pmatrix}
\begin{pmatrix}
(l-m+ \frac 12) t^l_{n \, \underline{m+ \frac 12}}(Z)  \\
(l+m+ \frac 12) t^l_{n \, \underline{m- \frac 12}}(Z)
\end{pmatrix}
=
\begin{pmatrix}
(l-n+1) t^{l + \frac 12}_{n - \frac 12 \, \underline{m}}(Z)  \\
(l+n+1) t^{l + \frac 12}_{n + \frac 12 \, \underline{m}}(Z)
\end{pmatrix}.
$$
\end{lem}

Applying $\square_{2,2} = \frac {\partial^2}{\partial z_{11}\partial z_{22}}
- \frac{\partial^2}{\partial z_{12}\partial z_{21}}$
and using Lemma \ref{deriv_calc} we obtain:

\begin{prop}  \label{t-harmonic}
On $\HR^+$ we have:
$$
\square_{2,2} \, t^l_{n\,\underline{m}}(Z) =0
\qquad \text{and} \qquad
\square_{2,2} \, N(Z)^{-1} \cdot t^l_{n\,\underline{m}}(Z^{-1}) =0.
$$
\end{prop}


Applying $\nabla^+ =
2 \begin{pmatrix} \frac {\partial}{\partial z_{11}} &
\frac {\partial}{\partial z_{21}}  \\
\frac {\partial}{\partial z_{12}} &
\frac {\partial}{\partial z_{22}} \end{pmatrix}$
and using Lemma \ref{deriv_calc} we obtain:

\begin{prop}  \label{left-reg_fns}
The following $\BB S$-valued functions on $\HR^+$
$$
\begin{pmatrix} (l-m+ \frac 12) t^l_{n \, \underline{m+ \frac 12}}(Z)  \\
(l+m+ \frac 12) t^l_{n \, \underline{m- \frac 12}}(Z) \end{pmatrix},
\qquad
\frac 1{N(Z)}
\begin{pmatrix} (l-n+ \frac 12) t^l_{n- \frac 12 \, \underline{m}}(Z^{-1})  \\
(l+n+ \frac 12) t^l_{n+ \frac 12 \, \underline{m}}(Z^{-1}) \end{pmatrix}
$$
are left-regular. Dually, the following $\BB S'$-valued functions on $\HR^+$
$$
\frac 1{N(Z)} \Bigl( t^{l+\frac 12}_{m \, \underline{n- \frac 12}}(Z^{-1}),
t^{l+\frac 12}_{m \, \underline{n+ \frac 12}}(Z^{-1}) \Bigr),
\qquad
\Bigl( t^{l-\frac 12}_{m + \frac 12 \, \underline{n}}(Z),
t^{l-\frac 12}_{m - \frac 12 \, \underline{n}}(Z) \Bigr)
$$
are right-regular.
\end{prop}

Next we would like to extend the matrix coefficients (\ref{int_t})
holomorphically from $\HR^+$ to an open region in $\HC$.
Define an open region in $\HC$:
$$
U = \biggl\{
Z = \begin{pmatrix} z_{11} & z_{12} \\ z_{21} & z_{22} \end{pmatrix} \in \HC
;\: |z_{21}| < |z_{11}|,\: |z_{12}| < |z_{22}|,\: \re (z_{11}z_{22}) >0
\text{ or } \im (z_{11}z_{22}) \ne 0 \biggr\}.
$$

\begin{lem}
The matrix coefficients (\ref{int_t}) can be holomorphically extended
to $U$.
\end{lem}

\pf
We can rewrite
\begin{equation}  \label{t-cont-ext}
t^l_{n\,\underline{m}}(Z) = \frac 1{2\pi i}
\oint \Bigl( 1 + s^{-1}\frac {z_{21}}{z_{11}} \Bigr)^{l-m}
\Bigl( 1 + s \frac {z_{12}}{z_{22}} \Bigr)^{l+m}
(z_{11}z_{22})^{l-m} (z_{22})^{2m} s^{n-m} \frac {ds}s.
\end{equation}
By assumption, $\re (z_{11}z_{22}) >0$ or $\im (z_{11}z_{22}) \ne 0$,
which allows us to choose a single branch of the complex multivalued
function $(z_{11}z_{22})^{l-m}$.
Since $|z_{21}| < |z_{11}|$ and $|z_{12}| < |z_{22}|$,
$$
\re \Bigl( 1 + s^{-1}\frac {z_{21}}{z_{11}} \Bigr) >0
\qquad \text{and} \qquad
\re \Bigl( 1 + s \frac {z_{12}}{z_{22}} \Bigr)>0.
$$
Thus we can choose single branches of the multivalued functions
$\bigl( 1 + s^{-1}\frac {z_{21}}{z_{11}} \bigr)^{l-m}$ and
$\bigl( 1 + s \frac {z_{12}}{z_{22}} \bigr)^{l+m}$ as well.
\qed

Let $U^+$ be another open region in $\HC$:
\begin{multline*}
U^+ = \{ Z \in \HC ;\: Z^+ \in U \} \\
= \biggl\{
Z = \begin{pmatrix} z_{11} & z_{12} \\ z_{21} & z_{22} \end{pmatrix} \in \HC
;\: |z_{12}| < |z_{11}|,\: |z_{21}| < |z_{22}|,\: \re (z_{11}z_{22}) >0
\text{ or } \im (z_{11}z_{22}) \ne 0 \biggr\}.
\end{multline*}
Observe that, for $Z \in U \cup U^+$, $N(Z)$ is never a negative real number
or zero, hence arbitrary complex powers of $N(Z)$ are well defined.
Letting $\tilde Z = Z^{-1}$ and writing
$$
t^l_{n\,\underline{m}}(Z) = t^l_{n\,\underline{m}}(\tilde Z^{-1}) =
N(\tilde Z)^{-2l} \cdot t^l_{n\,\underline{m}}(\tilde Z^+),
$$
we see that the matrix coefficients $t^l_{n\,\underline{m}}(Z)$'s
can be holomorphically extended to $U^+$ as well.
Now we define $\HC^+ = U \cup U^+$.

\begin{lem}
The matrix coefficients (\ref{int_t}) can be holomorphically extended
to $\HC^+$.
\end{lem}

\pf
Since the matrix coefficients extend to both $U$ and $U^+$, all we need to
show is that we do not run into problems with multivaluedness on $U \cap U^+$.
The matrix coefficients of the discrete series and their limits
are rational functions in $z_{ij}$'s and certainly extend to $\HC^+$.
For the matrix coefficients of the continuous series the result follows
from the following general observation:

Let $U^{\BB R}$ be an open region in $\BB R^n$, and let $U^{\BB C}_1$ and
$U^{\BB C}_2$ be two connected open regions in $\BB C^n$ such that
$U^{\BB R} \subset U^{\BB C}_i \cap \BB R^n$, $i=1,2$, and every loop in
$U^{\BB C}_1 \cup U^{\BB C}_2$ is homotopic to a loop in $U^{\BB R}$.
Suppose a real-analytic function $f$ on $U^{\BB R}$ has holomorphic
extensions $f_1$ to $U^{\BB C}_1$ and $f_2$ to $U^{\BB C}_2$,
then $f_1=f_2$ on $U^{\BB C}_1 \cap U^{\BB C}_2$ and $f$ has a unique
holomorphic extension to $U^{\BB C}_1 \cup U^{\BB C}_2$.
\qed

We view $\HC^+$ as an open neighborhood of $\HR^+$ in $\HC$, and
in light of the lemma we can regard the matrix coefficients (\ref{int_t})
as functions on $\HC^+$.
Lemma \ref{deriv_calc} and Proposition \ref{left-reg_fns} formally extend to
complex variables and $\HR^+$ replaced with $\HC^+$.

We conclude this subsection with a list of relations between the matrix
coefficients which follow from the relations between
$\mathfrak{P}^l_{n\,m}(t)$'s.
When $\re l = -1/2$, the functions $t^{\B{l}}_{n\,\underline{m}}(Z)$ and
$N(Z)^{-1} \cdot t^l_{n\,\underline{m}} ( Z/N(Z))$ on $\HC^+$ are proportional:
\begin{equation}  \label{m-coeff-symmetry-2}
t^{\B{l}}_{n\,\underline{m}}(Z) =
(-1)^{m-n} \frac{\Gamma(l-n+1)\Gamma(l+n+1)}{\Gamma(l-m+1)\Gamma(l+m+1)} \cdot
\frac 1{N(Z)} \cdot t^l_{n\,\underline{m}} \bigl( Z/N(Z) \bigr),
\qquad \re l = -\frac 12.
\end{equation}
We also have
\begin{equation}  \label{t(Z^+)}
t^l_{n\,\underline{m}}(Z^+) = 
(-1)^{m-n} \frac{\Gamma(l-m+1)\Gamma(l+m+1)}{\Gamma(l-n+1)\Gamma(l+n+1)} \cdot
t^l_{-m\,\underline{-n}}(Z).
\end{equation}
Combining this with (\ref{m-coeff-symmetry-2}) we obtain
\begin{equation}  \label{t(Z^{-1})}
N(Z)^{-1} \cdot t^l_{n\,\underline{m}} \bigl( Z^{-1} \bigr)
= t^{\B{l}}_{-m\,\underline{-n}}(Z),
\qquad \re l = -1/2.
\end{equation}

\subsection{Schwartz Space ${\cal S}(\HR^+)$ and Invariant Pairings}

\begin{df}
Let $f$ be a function on $\HR^+$ (with values in $\BB C$, $\HC$, $\BB S$ or
$\BB S'$). We say $f$ is {\em quasi-regular at the origin} if
$$
\lim_{t \to 0^+} t^{1+\delta}f(tX) =0
\qquad
\text{for all $X \in \HR^+$ and all $\delta>0$.}
$$

If $\phi$ is a solution of the wave equation $\square_{2,2}\phi=0$ on $\HR^+$,
we say $\phi$ is {\em quasi-regular at infinity} if
$$
\pi^0_l \begin{pmatrix} 0 & 1 \\ 1 & 0 \end{pmatrix} \phi =
\pi^0_r \begin{pmatrix} 0 & 1 \\ 1 & 0 \end{pmatrix} \phi =
\frac 1{N(X)} \cdot \phi(X^{-1})
$$
is quasi-regular at the origin.

Similarly, we say that a left-regular function $f$ or a right-regular
function $g$ on $\HR^+$ is {\em quasi-regular at infinity} if
$$
\pi_l \begin{pmatrix} 0 & 1 \\ 1 & 0 \end{pmatrix} f =
\frac {X^{-1}}{N(X)} \cdot f(X^{-1})
\qquad \text{or} \qquad
\pi_r \begin{pmatrix} 0 & 1 \\ 1 & 0 \end{pmatrix} g =
g(X^{-1}) \cdot \frac {X^{-1}}{N(X)}
$$
is quasi-regular at the origin.
\end{df}

For example, let $f$ be a homogeneous function on $\HR^+$ of homogeneity
degree $\mu \in \BB C$.
Then $f$ is quasi-regular at the origin if and only if $\re \mu \ge -1$.
If $f$ is a solution of the wave equation, then $f$ is quasi-regular
at infinity if and only if $\re \mu \le -1$.
In particular, all matrix coefficients $t^l_{n\,\underline{m}}(Z)$ are 
quasi-regular at infinity; moreover, the matrix coefficients of the
continuous series and the limits of the discrete series ($\re l = -1/2$)
are also quasi-regular at the origin.
Finally, if $f$ is a left- or right-regular function, then $f$ is
quasi-regular at infinity if and only if $\re \mu \le -2$.

For the purposes of this article, a suitable reference for the Schwartz
functions on Lie groups in the sense of Harish-Chandra is \cite{Va}.
Our definition of the Schwartz space on $\HR^+$, denoted by ${\cal S}(\HR^+)$,
is motivated by the following needs:
\begin{itemize}
\item
If $f \in {\cal S}(\HR^+)$, then the restrictions of $f$ to all hyperboloids
$H_R$, $R>0$, should be Schwartz functions in the sense of Harish-Chandra;
\item
Since we often apply the $\deg$ operator, if $f \in {\cal S}(\HR^+)$, then the
restrictions $(\operatorname{deg})^d f \bigr|_{H_R}$ should be Schwartz
functions as well, for all $R>0$ and integers $d \ge 0$;
\item
Since we often perform calculations with matrix coefficients
$t^l_{n\,\underline{m}}(Z)$, the matrix coefficients of the discrete series and
wave packets formed out of the matrix coefficients of the continuous series
should be in ${\cal S}(\HR^+)$.
\end{itemize}

The group $\HR^+ \times \HR^+$ acts on $\HR^+$ by multiplication:
$$
(a,b): X \mapsto aXb^{-1}, \qquad (a,b) \in \HR^+ \times \HR^+, \: X \in \HR^+.
$$
Hence we get an action of $\HR^+ \times \HR^+$ on the space of smooth
functions on $\HR^+$:
$$
(a,b): f(X) \mapsto f(a^{-1}Xb), \qquad f \in {\cal C}^{\infty} (\HR^+).
$$
Differentiating, we obtain an action of the Lie algebra
$\mathfrak{gl}(2,\BB R) \times \mathfrak{gl}(2,\BB R)$ which extends to
an action of its universal enveloping algebra
${\cal U}(\mathfrak{gl}(2,\BB R)) \otimes {\cal U}(\mathfrak{gl}(2,\BB R))$
on ${\cal C}^{\infty} (\HR^+)$.

\begin{df}
We define ${\cal S}(\HR^+)$ -- the {\em Schwartz space on $\HR^+$} --
to be the space of smooth functions $f$ on $\HR^+$ such that, for each
$a \in {\cal U}(\mathfrak{gl}(2,\BB R))\otimes{\cal U}(\mathfrak{gl}(2,\BB R))$,
each $n \in \BB N$ and all $R>1$, there exists a constant $C(a,n,R)>0$
such that
$$
|a \cdot f(X)| \le C(a,n,R) \cdot \frac  {(1+\log \|X\|)^{-n}}{\|X\|}
$$
for all $X \in \HR^+$ such that $R^{-2} \le N(X) \le R^2$.
\end{df}

We shall regard ${\cal S}(\HR^+)$ as a Fr\'echet space with respect to the
seminorms
$$
\mu_{a,n,R}(f) = \sup_{\{X \in \HR^+ ;\: R^{-2} \le N(X) \le R^2\}}
|a \cdot f(X)| \cdot \|X\| \cdot (1+\log \|X\|)^n.
$$
Denoting by ${\cal C}^{\infty}_c(\HR^+)$ the space of compactly supported
smooth functions on $\HR^+$ we get maps
$$
{\cal C}^{\infty}_c(\HR^+) \hookrightarrow {\cal S}(\HR^+)
\to L^2\bigl( SU(1,1) \bigr ),
$$
where the first map is the inclusion and the second map is the
restriction map $f \mapsto f \bigr|_{SU(1,1)}$.
These maps are continuous with dense images.

Let
$$
{\cal H}(\HR^+) = \{ \phi \in {\cal S}(\HR^+) ;\: \square_{2,2} \phi =0\}
$$
be the space of ``harmonic'' functions in ${\cal S}(\HR^+)$.
Similarly, we denote by $\BB S(\HR^+)$ the space of $\BB S$-valued
left-regular functions and, respectively, $\BB S'(\HR)$ the space of
$\BB S'$-valued right-regular functions on $\HR^+$ with both components in
${\cal S}(\HR^+)$.
Essentially by definition, the matrix coefficient functions listed in
Proposition \ref{t-harmonic} with $l \le -1$ (i.e. belonging to the
discrete series) together with the wave packets formed out the functions
with $\re l = -1/2$ (i.e. belonging to the continuous series) form a
dense subset of ${\cal H}(\HR^+)$.
Similarly, the left- and right-regular functions listed in Proposition
\ref{left-reg_fns} generate the spaces $\BB S(\HR^+)$ and $\BB S'(\HR)$
respectively. That is, the left- and right-regular functions of
Proposition \ref{left-reg_fns} with entries belonging to the discrete series
together with the wave packets formed out of the functions with entries in
the continuous series form dense subsets of $\BB S(\HR^+)$ and $\BB S'(\HR)$
respectively.

Fix an $R>0$ and define a bilinear form on ${\cal H} (\HR^+)$ by
\begin{equation}  \label{b-form}
\langle \phi_1 , \phi_2 \rangle_R = - \frac 1{2\pi^2}
\int_{X \in H_R} (\deg \phi_1)(X) \cdot \phi_2(X) \, \frac{dS}{\|X\|},
\qquad \phi_1, \phi_2 \in {\cal H}(\HR^+).
\end{equation}
This form is not symmetric, not $\mathfrak{sl}(2,\HC)$-invariant and
depends on the choice of $R$. In Subsection \ref{extension} we will extend
the space ${\cal H}(\HR^+)$ to $\widehat{\cal H}(\HR^+)$ and define a symmetric
$\mathfrak{sl}(2,\HC)$-invariant nondegenerate bilinear pairing on it.

\begin{lem}
Let $\tilde G(H_R) \subset GL(2,\HC)$ be the subgroup consisting
of all elements of $GL(2,\HC)$ with entries in $\HR$ and
preserving the hyperboloid $H_R = \{X \in \HR ;\: N(X)=R^2\}$. Then
$$
Lie(\tilde G(H_R)) =
\biggl\{ \begin{pmatrix} A & R^2C^+ \\ C & D \end{pmatrix} ;\:
A,C,D \in \HR, \: \re A = \re D \biggr\}.
$$
\end{lem}

\pf
The Lie algebra of $\tilde G(H_R)$ consists of all matrices
$\begin{pmatrix} A & B \\ C & D \end{pmatrix}$, $A,B,C,D \in \HR$,
which generate vector fields tangent to $H_R$.
Such a matrix $\begin{pmatrix} A & B \\ C & D \end{pmatrix}$
generates a vector field
$$
\frac d{dt} \bigl( (1+tA)X+tB \bigr) (tCX+1+tD)^{-1} \Bigr|_{t=0}
= AX+B-XCX-XD.
$$
A vector field is tangent to $H_R$ if and only if it is orthogonal with
respect to (\ref{bilinear_form}) to the vector field $X$ for $N(X)=R^2$:
$$
0=\re \bigl( (AX+B-XCX-XD)X^+ \bigr) = \re (R^2A-R^2D+BX^+-R^2XC),
\qquad \forall X \in H_R.
$$
It follows that $\re A = \re D$ and $B=R^2C^+$.
\qed

\begin{cor}  \label{G(H_R)}
Let $G(H_R) \subset GL(2, \HC)$ be the connected subgroup with Lie algebra
\begin{multline*}
\mathfrak g(H_R) = \bigl\{ x \in Lie(\tilde G(H_R));\: \re(\tr x)=0 \bigr\}  \\
= \biggl\{ \begin{pmatrix} A & R^2C^+ \\ C & D \end{pmatrix} ;\:
A,C,D \in \HR, \: \re A = \re D =0 \biggr\}.
\end{multline*}
Then $G(H_R)$ preserves the hyperboloid $H_R = \{ X \in \HR ;\: N(X)=R^2 \}$
and the open sets $\{ X \in \HR ;\: N(X)>R^2 \}$,
$\{ X \in \HR ;\: N(X)<R^2 \}$.
\end{cor}

The Lie algebra $\mathfrak g(H_R)$ and the Lie group $G(H_R)$ are isomorphic
to $\mathfrak{so}(3,2) = \mathfrak{sp}(2, \BB R)$ and $SO^+(3,2)$ respectively
(see, for example, \cite{H}).

\begin{prop}  \label{invar-pairing-prop}
The bilinear pairing $\langle \phi_1 , \phi_2 \rangle_R$
is invariant under the $\pi^0_l$ action of $\mathfrak{g}(H_R)$.
\end{prop}

\pf
First we find a convenient pair of subgroups generating $G(H_R)$:

\begin{lem}  \label{G(H_R)'}
The group $G(H_R)$ is generated by $SU(1,1) \times SU(1,1)$ realized as
the subgroup of diagonal matrices 
$\begin{pmatrix} a & 0 \\ 0 & d \end{pmatrix}$, $a,d \in SU(1,1) \subset \HR$,
and the one-parameter group
$$
G(H_R)' = \biggl\{ \begin{pmatrix}
\cosh t & R\sinh t \\ R^{-1}\sinh t & \cosh t \end{pmatrix}
;\: t \in \BB R \biggr\}.
$$
\end{lem}

Clearly, the bilinear pairing is invariant under the $\pi^0_l$ action of
$SU(1,1) \times SU(1,1)$.
Thus it is sufficient to show it is invariant under a one-dimensional
Lie algebra $\mathfrak{g}(H_R)'=Lie(G(H_R)')$.

\begin{lem}  \label{Jacobian-lemma}
Fix an element
$g = \begin{pmatrix} \cosh t & R\sinh t \\
R^{-1} \sinh t & \cosh t \end{pmatrix} \in G(H_R)'$
and consider its conformal action on $H_R$:
$$
\pi_l(g):\: X \mapsto \tilde X =
(\cosh t X - R\sinh t) (-R^{-1}\sinh t X + \cosh t)^{-1}.
$$
Then the Jacobian $J$ of this map is
$$
\pi_l(g)^* \biggl( \frac{dS}{\|X\|} \biggr) = J\,\frac{dS}{\|X\|}
= \frac 1{N(-R^{-1}\sinh t X + \cosh t)} \cdot
\frac{R^2 - (\re \tilde X)^2}{R^2 - (\re X)^2} \,\frac{dS}{\|X\|}.
$$
\end{lem}

Let $g=\exp\begin{pmatrix} 0 & Rt \\ R^{-1}t & 0 \end{pmatrix} \in G(H_R)'$.
For $t \to 0$ and modulo terms of order $t^2$, we have:
\begin{align}
\tilde X &= X+t(X^2-R^2)/R,  \label{X-tilde-1}  \\
N(-R^{-1}\sinh tX + \cosh t) &= 1 - 2t \re X /R, \label{X-tilde-2}\\
R^2 - (\re \tilde X)^2 &=
\bigl( R^2 -(\re X)^2 \bigr) (1+4t \re X/R). \label{X-tilde-3}
\end{align}

\begin{lem}  \label{deg-x-tilde}
Modulo terms of order $t^2$ we have
$$
\deg_X \bigl( \pi^0_l(g) \phi \bigr) = \deg_X \biggl(
\frac {\phi(\tilde X)}{N(-R^{-1}\sinh tX + \cosh t)} \biggr)
= (1 + 4t \re X /R) \cdot (\deg_X \phi)(\tilde X).
$$
\end{lem}

Letting $\tilde\phi_i = \pi^0_l(g) \phi_i$, $i=1,2$, and
continuing to work modulo terms of order $t^2$, we get
\begin{multline*}
-2\pi^2 \cdot \langle \tilde\phi_1 , \tilde\phi_2 \rangle_R
= \int_{X \in H_R} \deg_X (\tilde \phi_1)(X)
\cdot \tilde \phi_2(X) \,\frac{dS}{\|X\|}  \\
= \int_{X \in H_R} (\deg_X \phi_1)(\tilde X) \cdot \phi_2(\tilde X) \cdot
\frac{1 + 4t \re X /R}{N(-R^{-1}\sinh tX + \cosh t)} \,\frac{dS}{\|X\|}  \\
= \int_{\tilde X \in H_R}
(\deg_X \phi_1)(\tilde X) \cdot \phi_2(\tilde X) \,\frac{dS}{\|\tilde X\|}
= -2\pi^2 \cdot \langle \phi_1 , \phi_2 \rangle_R.
\end{multline*}
This proves that the bilinear pairing $\langle \phi_1 , \phi_2 \rangle_R$
is $\mathfrak{g}(H_R)$-invariant.
\qed

Next we calculate the pairings between the matrix coefficient functions.

\begin{prop}  \label{t-orthogonal}
The matrix coefficients (\ref{int_t}) satisfy the following orthogonality
relationships:
\begin{multline*}
\Bigl\langle t^{l'}_{n'\,\underline{m'}}(X),
\frac 1{N(X)} \cdot t^l_{m\,\underline{n}}(X^{-1}) \Bigr\rangle_R
= - \Bigl\langle \frac 1{N(X)} \cdot t^l_{m\,\underline{n}}(X^{-1}),
t^{l'}_{n'\,\underline{m'}}(X) \Bigr\rangle_R  \\
= \begin{cases}
\delta_{ll'} \delta_{mm'} \delta_{nn'} & \text {if $l,l' \ne - 1/2$;}\\
0 & \text{if $l$ or $l'=-1/2$,}
\end{cases}
\end{multline*}
\begin{multline*}
R^{-2(2l+1)} \cdot \bigl\langle t^{l}_{n\,\underline{m}}(X),
t^{l'}_{-n'\,\underline{-m'}}(X) \bigr\rangle_R  \\
=
-R^{2(2l+1)} \cdot \Bigl\langle \frac 1{N(X)} \cdot t^l_{n\,\underline{m}}(X^{-1}),
\frac 1{N(X)} \cdot t^{l'}_{-n'\,\underline{-m'}}(X^{-1}) \Bigr\rangle_R  \\
=
(-1)^{m-n}
\frac{\Gamma(l+m+1) \Gamma(l-m+1)}{\Gamma(l+n+1) \Gamma(l-n+1)} \cdot
\begin{cases}
\delta_{ll'} \delta_{mm'} \delta_{nn'} & \text {if $l,l' \ne - 1/2$;} \\
0 & \text{if $l$ or $l'=-1/2$,}
\end{cases}
\end{multline*}
for all $R>0$, provided that $l$ or $l'$ lies in
$\{-\frac12, -1,-\frac32,\dots\}$.
In particular, we get a nondegenerate pairing between the spaces of
harmonic functions spanned by $t^l_{n\,\underline{m}}(X)$'s and by
$N(X)^{-1} \cdot t^l_{n\,\underline{m}}(X^{-1})$'s with $l=-1,-\frac32,-2,\dots$,
which is independent of the choice of $R>0$.
\end{prop}

\pf
Since the functions $N(X)^{-1} \cdot t^l_{m\,\underline{n}}(X^{-1})$ and
$t^{l'}_{n'\,\underline{m'}}(X)$ are homogeneous of degree $-2(l+1)$ and $2l'$
respectively, it is sufficient to prove the result for $R=1$ only.
This is done by integrating in coordinates (\ref{param}) using identities
(\ref{t-german}) and (\ref{t(Z^+)}) and the orthogonality relations for
$\mathfrak{P}^l_{n\,m}$'s (\ref{P-orthog})-(\ref{P-norm}).
\qed

The $\mathfrak{gl}(2,\HC)$-invariant pairing between left-regular and
right-regular functions on $\HR^+$ is given by the formula
$$
\langle g, f \rangle = - \frac 1{2\pi^2}
\int_{X \in SU(1,1)} g(X) \cdot Dx \cdot f(X)
$$
(provided that the integral converges).
The fact that this bilinear pairing is $\mathfrak{gl}(2,\HC)$-invariant can
be seen as follows. First, one can  realize that if $SU(1,1)$ is replaced
with an arbitrary hyperboloid $H_R$, the pairing remains unchanged
(this follows, for example, from Proposition \ref{t-orthogonal2} below).
Then, exactly as in the proof of Proposition \ref{invar-pairing-prop},
 one can show that the pairing is invariant with respect to each algebra
$\mathfrak{g}(H_R)$.
But for $R_1 \ne R_2$ the algebras $\mathfrak{g}(H_{R_1})$ and
$\mathfrak{g}(H_{R_2})$ generate all of $\mathfrak{sl}(2,\HR)$.
Finally, one can check the invariance under diagonal matrices and get
the $\mathfrak{gl}(2,\HC)$-invariance.

\begin{prop} \label{t-orthogonal2}
We have the following orthogonality relationships:
$$
\left\langle
\frac 1{N(X)} \Bigl( t^{l}_{m \, \underline{n- \frac 12}}(X^{-1}),
t^{l}_{m \, \underline{n+ \frac 12}}(X^{-1}) \Bigr),
\begin{pmatrix}
(l'-m') t^{l'-\frac 12}_{n' \, \underline{m'+ \frac 12}}(X)  \\
(l'+m') t^{l'-\frac 12}_{n' \, \underline{m'- \frac 12}}(X)
\end{pmatrix} \right\rangle
= \delta_{ll'} \delta_{mm'} \delta_{nn'},
$$
$$
\left\langle
\Bigl( t^{l'- \frac 12}_{m'+\frac12 \, \underline{n'}}(X),
t^{l'- \frac 12}_{m'-\frac12 \, \underline{n'}}(X) \Bigr),
\frac 1{N(X)}
\begin{pmatrix}
(l-n+\frac12) t^{l}_{n-\frac12 \, \underline{m}} (X^{-1})  \\
(l+n+\frac12)t^{l}_{n+\frac12 \, \underline{m}} (X^{-1})
\end{pmatrix}
\right\rangle
= \delta_{ll'} \delta_{mm'} \delta_{nn'},
$$
$$
\left\langle \Bigl( t^{l}_{n+\frac12 \, \underline{m}}(X),
t^{l}_{n-\frac12 \, \underline{m}}(X) \Bigr),
\begin{pmatrix}
(l'-m') t^{l'-\frac12}_{n' \, \underline{m'+ \frac 12}}(X)  \\
(l'+m') t^{l'-\frac12}_{n' \, \underline{m'- \frac 12}}(X)
\end{pmatrix} \right\rangle = 0,
$$
$$
\left\langle \frac 1{N(X)} \Bigl( t^{l}_{n \, \underline{m- \frac 12}}(X^{-1}),
t^{l}_{n \, \underline{m+ \frac 12}}(X^{-1}) \Bigr),
\frac 1{N(X)}
\begin{pmatrix}
(l'-n') t^{l'-\frac12}_{n'-\frac12 \, \underline{m'}} (X^{-1})  \\
(l'+n')t^{l'-\frac12}_{n'+\frac12 \, \underline{m'}} (X^{-1})
\end{pmatrix}
\right\rangle = 0,
$$
provided that $l$ or $l'-\frac12$ lies in the discrete series set
$\{-1,-\frac32,-2,\dots\}$.
\end{prop}

\pf
Using Lemmas \ref{restrictions} and \ref{mult-ids}, we obtain:
\begin{multline*}
\Bigl( t^{l}_{m \, \underline{n- \frac 12}}(X^{-1}),
t^{l}_{m \, \underline{n+ \frac 12}}(X^{-1}) \Bigr)
\cdot X \cdot
\begin{pmatrix}
(l'-m') t^{l'-\frac 12}_{n' \, \underline{m'+ \frac 12}}(X)  \\
(l'+m') t^{l'-\frac 12}_{n' \, \underline{m'- \frac 12}}(X)
\end{pmatrix} \\
=
(l'-n'+1/2) \cdot t^{l'}_{n' - \frac 12 \, \underline{m'}}(X)
\cdot t^{l}_{m \, \underline{n- \frac 12}}(X^{-1})
+
(l'+n'+1/2) \cdot t^{l'}_{n' + \frac 12 \, \underline{m'}}(X)
\cdot t^{l}_{m \, \underline{n+ \frac 12}}(X^{-1}),
\end{multline*}
\begin{multline*}
\Bigl( t^{l'- \frac 12}_{m'+\frac12 \, \underline{n'}}(X),
t^{l'- \frac 12}_{m'-\frac12 \, \underline{n'}}(X) \Bigr)
\cdot X \cdot
\begin{pmatrix}
(l-n+\frac12) t^{l}_{n-\frac12 \, \underline{m}} (X^{-1})  \\
(l+n+\frac12)t^{l}_{n+\frac12 \, \underline{m}} (X^{-1})
\end{pmatrix}  \\
=
(l-n+1/2) \cdot t^{l'}_{m'\, \underline{n'-\frac12}}(X)
\cdot t^{l}_{n-\frac12 \, \underline{m}} (X^{-1})
+ (l+n+1/2) \cdot t^{l'}_{m'\, \underline{n'+\frac12}}(X)
\cdot t^{l}_{n+\frac12 \, \underline{m}} (X^{-1}),
\end{multline*}
then the result follows from Proposition \ref{t-orthogonal}.
\qed

\subsection{Fueter Formula for Hyperboloids}

We fix $0< R' < R$ and let $U = \{ X \in \HR ;\: R'^2 < N(X) < R^2 \}$
be the open region in $\HR^+$ bounded by two hyperboloids.
In this section we essentially substitute this unbounded set $U$
and bounded functions $f$ into the integral formula from
Theorem \ref{holomorphic_Fueter} and prove that the resulting identity
still holds.

Recall the deformation $h_{\epsilon, Z_0}: \HC \to \HC$,
$Z \mapsto Z + i \epsilon (Z-Z_0)^-$.
First we deal with a technical issue arising from the fact that
the image of $U$ under $h_{\epsilon,Z_0}$ does not lie inside $\HC^+$.
Thus we modify the deformation $h_{\epsilon, Z_0}$ as
$$
\tilde h_{\epsilon, Z_0}: \HC \to \HC, \qquad
Z \mapsto Z + \epsilon (1+i) (Z-Z_0)^-.
$$
Note that
$$
\tilde h_{\epsilon, Z_0}^* \bigl(N(Z-Z_0) \bigr) =
\bigl( N(Z-Z_0) + 2\epsilon S(Z-Z_0) \bigr)
+ 2i\epsilon \bigl( S(Z-Z_0) + \epsilon N(Z-Z_0) \bigr),
$$
so, for $|\epsilon|$ sufficiently small, the denominator of
$\frac{(Z-X_0)^{-1}}{N(Z-X_0)}$ is bounded away
from zero on $(\tilde h_{\epsilon, X_0})_*(\partial U)$.

\begin{lem}
Fix any $Z_0 \in \HC$, then, for $\epsilon >0$ sufficiently small,
$\tilde h_{\epsilon, Z_0} (\overline{U}) \subset \HC^+$.
\end{lem}

\pf
Let $Z = \begin{pmatrix} z_{11} & z_{12} \\ z_{21} & z_{22} \end{pmatrix}$,
$Z_0 = \begin{pmatrix} y_{11} & y_{12} \\ y_{21} & y_{22} \end{pmatrix}$
and $\tilde h_{\epsilon, Z_0}(Z) =
\begin{pmatrix} z_{11}' & z_{12}' \\ z_{21}' & z_{22}' \end{pmatrix}$.
We can write
$$
\begin{pmatrix} z_{11}' & z_{12}' \\ z_{21}' & z_{22}' \end{pmatrix}
=
\begin{pmatrix} (1+\epsilon(1+i)) z_{11} & (1-\epsilon(1+i))z_{12} \\
(1-\epsilon(1+i)) z_{21} & (1+\epsilon(1+i))z_{22} \end{pmatrix}
- \epsilon(1+i) \begin{pmatrix} y_{11} & -y_{12} \\ -y_{21} & y_{22} \end{pmatrix}.
$$
Since $|z_{11}|>|z_{21}|$ and $|z_{22}|>|z_{12}|$ for all $Z \in \B{U}$,
by taking $\epsilon >0$ sufficiently small we can arrange that
$$
\bigl| (1+\epsilon(1+i)) z_{11} \bigr| - \bigl| (1-\epsilon(1+i)) z_{21} \bigr|
> 2\sqrt{2} \epsilon |y_{ij}|,
$$
$$
\bigl| (1+\epsilon(1+i)) z_{22} \bigr| - \bigl| (1-\epsilon(1+i)) z_{12} \bigr|
> 2\sqrt{2} \epsilon |y_{ij}|.
$$
Then $\tilde h_{\epsilon, Z_0}(Z)$ satisfies the inequalities
$|z_{21}'| < |z_{11}'|$ and $|z_{12}'| < |z_{22}'|$.
Finally,
\begin{multline*}
z_{11}'z_{22}' =
\bigl( (1+\epsilon(1+i)) z_{11} - \epsilon(1+i) y_{11} \bigr)
\bigl( (1+\epsilon(1+i)) z_{22} - \epsilon(1+i) y_{22} \bigr)  \\
= (1+2\epsilon)z_{11}z_{22} + 2i\epsilon(1+\epsilon)z_{11}z_{22}
+ 2i\epsilon^2y_{11}y_{22}
- \epsilon(1+i)(1+\epsilon(1+i)) (z_{11}y_{22}+z_{22}y_{11}).
\end{multline*}
The first term $(1+2\epsilon)z_{11}z_{22}$ is a positive real number,
the second and the third terms $2i\epsilon(1+\epsilon)z_{11}z_{22}$ and
$2i\epsilon^2y_{11}y_{22}$ are purely imaginary,
finally the fourth term
$\epsilon(1+i)(1+\epsilon(1+i)) (z_{11}y_{22}+z_{22}y_{11})$
becomes smaller in magnitude than $(1+2\epsilon)z_{11}z_{22}$
when $\epsilon>0$ is sufficiently small.
This proves $\re (z_{11}'z_{22}') >0$.
\qed

\begin{thm}  \label{split-Fueter}
Suppose that $f(Z)$ is a left-regular function on $\HC^+$ such that its
components are bounded on closed sets
\begin{equation}  \label{closed_sets}
\biggl\{
Z = \begin{pmatrix} z_{11} & z_{12} \\ z_{21} & z_{22} \end{pmatrix} \in \HC;\quad
\begin{matrix} c \le |z_{11}/z_{21}| \le 2c, \quad |\B{z_{11}}/z_{22}-1|\le d, \\
c \le |z_{22}/z_{12}| \le 2c, \quad |\B{z_{12}}/z_{21}-1| \le d, \end{matrix} \quad
|z_{11}|,|z_{22}| \ge L_0 \biggr\}
\end{equation}
for all $c>1$ sufficiently close to 1 and some fixed values $d, L_0>0$.
Let $X_0 \in \HR$, then for $\epsilon >0$ sufficiently small
$$
- \frac 1 {2\pi^2} \int_{(\tilde h_{\epsilon, X_0})_*(\partial U)}
\frac {(Z-X_0)^{-1}} {N(Z-X_0)} \cdot Dz \cdot f(Z) =
\begin{cases}
f(X_0) & \text{if $X_0 \in U$;}\\
0 & \text{if $X_0 \notin \overline{U}$.}
\end{cases}
$$
\end{thm}

\pf
For $L \in \BB R$, let $B_L$ denote the open ball
$\{ X \in \HR ;\: \|X\|<L \}$ of radius $L$, and set $U_L = U \cap B_L$.
Clearly, the closure of $U_L$ is compact, and the proof of
Theorem \ref{holomorphic_Fueter} shows that for $L$ sufficiently large
$$
\int_{(\tilde h_{\epsilon, X_0})_*(\partial U_L)}
\frac {(Z-X_0)^{-1}} {N(Z-X_0)} \cdot Dz \cdot f(Z)
=
\begin{cases}
- 2\pi^2 f(X_0) & \text{if $X_0 \in U$;} \\
0 & \text{if $X_0 \notin \overline{U}$.}
\end{cases}
$$
The support of $\partial U_L$ consists of a portion which lies in
$H_R \cup H_{R'}$ and its complement.
We define 3-chains $C_L^1$ and $C_L^2$ by
$$
\partial U_L = C_L^1 + C_L^2
$$
with supports $|C_L^1| \subset H_R \cup H_{R'}$ and $|C_L^2| \subset B_L$.
Then
$$
\int_{(\tilde h_{\epsilon, X_0})_*(\partial U)}
\frac {(Z-X_0)^{-1}} {N(Z-X_0)} \cdot Dz \cdot f(Z)
=
\lim_{L \to \infty} \int_{(\tilde h_{\epsilon, X_0})_*(\partial U_L - C_L^2)}
\frac {(Z-X_0)^{-1}} {N(Z-X_0)} \cdot Dz \cdot f(Z).
$$
Thus the theorem will follow if we can show
$$
\lim_{L \to \infty} \int_{(\tilde h_{\epsilon, X_0})_*(C_L^2)}
\frac {(Z-X_0)^{-1}} {N(Z-X_0)} \cdot Dz \cdot f(Z) =0.
$$
First we choose $\epsilon$ small enough so that
$$
\biggl| \frac{1+\epsilon(1+i)}{1+\epsilon(1-i)} -1 \biggr| <\frac d2.
$$
Then we choose a number $c$ such that
$$
1 < c < \biggl| \frac{1+\epsilon(1+i)}{1-\epsilon(1+i)} \biggr| <2c.
$$
As $L$ tends to infinity, the supports of $C_L^2$ will lie inside the set
(\ref{closed_sets}), hence the components of $f$ are bounded on $|C_L^2|$.
The volume of $|C_L^2|$ grows as $O(L)$, the denominator of
$\frac {(Z-X_0)^+} {N(Z-X_0)^2}$ grows as $O(L^4)$,
and the numerator grows as $O(L)$.
This proves that the integral over $C_L^2$ tends to zero.
\qed

\begin{rem}
From (\ref{t-cont-ext}) one can see that the matrix coefficients
(\ref{int_t}) are bounded on the closed sets (\ref{closed_sets})
as long as $\re l \le 0$.
\end{rem}

\section{The Discrete Series Component on $\HR$}

\subsection{The Polynomial Algebra on $\HR^+$}  \label{polynomial_algebra}

Let us consider the matrix coefficients $t^l_{n\,\underline{m}}(Z)$'s of
the discrete series together with the limits of the discrete series.
It is easy to see that they are rational functions in
$z_{11}, z_{12}, z_{21}, z_{22}$, hence uniquely extend to $\HC$.
In this subsection we identify certain spaces of rational functions spanned by
these matrix coefficients.

We think of the matrix coefficients $t^l_{n\,\underline{m}}(Z)$ and
$N(Z)^{-2l-1} \cdot t^l_{n\,\underline{m}}(Z)$ as meromorphic functions on $\HC$
and introduce the following spaces:
$$
{\cal D}^-_{discr} =  \text{$\BB C$-span of $t^l_{n\,\underline{m}}(Z)$ and
$N(Z)^{-2l-1} \cdot t^l_{n\,\underline{m}}(Z)$},
\qquad \begin{matrix} l = -1, -\frac 32, -2, \dots, \\
m,n \in \BB Z+l, \qquad m, n \ge -l, \end{matrix}
$$
(span of the matrix coefficients of the holomorphic discrete series),
$$
{\cal D}^-_{lim} =  \text{$\BB C$-span of $t^{-\frac12}_{n\,\underline{m}}(Z)$},
\qquad m,n \in \BB Z+1/2, \qquad m, n \ge 1/2,
$$
(span of the matrix coefficients of the limit of the holomorphic discrete
series),
$$
{\cal D}^+_{discr} =  \text{$\BB C$-span of $t^l_{n\,\underline{m}}(Z)$ and
$N(Z)^{-2l-1} \cdot t^l_{n\,\underline{m}}(Z)$},
\qquad \begin{matrix} l = -1, -\frac 32, -2, \dots, \\
m,n \in \BB Z+l, \qquad m, n \le l, \end{matrix}
$$
(span of the matrix coefficients of the antiholomorphic discrete series),
$$
{\cal D}^+_{lim} =  \text{$\BB C$-span of $t^{-\frac12}_{n\,\underline{m}}(Z)$},
\qquad m,n \in \BB Z+1/2, \qquad m, n \le -1/2,
$$
(span of the matrix coefficients of the limit of the antiholomorphic discrete
series),
$$
{\cal D}^- = {\cal D}^-_{discr} \oplus {\cal D}^-_{lim},
\qquad
{\cal D}^+ = {\cal D}^+_{discr} \oplus {\cal D}^+_{lim}.
$$
Note that the matrix entries $z_{11}$ and $z_{22}$ are invertible on $\HR^+$.
We denote by
$$
\BB C [z_{11}^{-1}, z_{12}, z_{21}, z_{22}]^{\le 0} \quad \subset \quad
\BB C [z_{11}^{-1}, z_{12}, z_{21}, z_{22}]
$$
the subalgebra of polynomials in $\BB C [z_{11}^{-1}, z_{12}, z_{21}, z_{22}]$
spanned by monomials $z_{11}^{d_{11}} \cdot z_{12}^{d_{12}} \cdot
z_{21}^{d_{21}} \cdot z_{22}^{d_{22}}$ with $d_{11}+d_{12}+d_{21}+d_{22} \le 0$.
Similarly we can define a subalgebra
$$
\BB C [z_{11}, z_{12}, z_{21}, z_{22}^{-1}]^{\le 0} \quad \subset \quad
\BB C [z_{11}, z_{12}, z_{21}, z_{22}^{-1}].
$$
Let ${\cal D}^-_{\le 0}$ (respectively ${\cal D}^-_{\ge 0}$) be the $\BB C$-span
of $t^l_{n\,\underline{m}}(Z)$ (respectively
$t^l_{n\,\underline{m}}(Z) \cdot N(Z)^{-2l-1}$) with $l=-\frac12,-1,-\frac32,\dots$,
$m,n \in \BB Z+l$, $m, n \ge -l$.
Then  ${\cal D}^- = {\cal D}^-_{\le 0} + {\cal D}^-_{\ge 0}$ and
${\cal D}^-_{\ge 0}$ is the image of ${\cal D}^-_{\le 0}$ under the
inversion map $\phi(Z) \mapsto N(Z)^{-1} \cdot \phi \bigl(Z/N(Z)\bigr)$.
Similarly, we can define ${\cal D}^+_{\le 0}$, ${\cal D}^+_{\ge 0}$ so that
${\cal D}^+ = {\cal D}^+_{\le 0} + {\cal D}^+_{\ge 0}$ and ${\cal D}^+_{\ge 0}$ is
the image of ${\cal D}^+_{\le 0}$ under the inversion map.

\begin{prop}
We have:
\begin{align*}
{\cal D}^-_{\le 0} &=
\bigl\{ \phi \in z_{11}^{-1} \cdot \BB C [z_{11}^{-1}, z_{12}, z_{21}, z_{22}]^{\le 0}
;\: \square_{2,2} \phi =0 \bigr\}, \\
{\cal D}^+_{\le 0} &=
\bigl\{ \phi \in z_{22}^{-1} \cdot \BB C [z_{11}, z_{12}, z_{21}, z_{22}^{-1}]^{\le 0}
;\: \square_{2,2} \phi =0 \bigr\}.
\end{align*}
\end{prop}

\begin{rem}
Note that ${\cal D}^-$ is a proper subspace of
$\bigl\{ \phi \in z_{11}^{-1} \cdot
\BB C [z_{11}^{-1}, z_{12}, z_{21}, z_{22}] ;\: \square_{2,2} \phi =0 \bigr\}$,
since the functions $z_{11}^{d_{11}} \cdot z_{12}^{d_{12}}$ and 
$z_{11}^{d_{11}} \cdot z_{21}^{d_{21}}$ with $d_{11}<0$, $d_{12},d_{21} \ge -d_{11}$
are not in ${\cal D}^-$. Similarly, ${\cal D}^+$ is a proper subspace of
$\bigl\{ \phi \in z_{22}^{-1} \cdot \BB C [z_{11}, z_{12}, z_{21}, z_{22}^{-1}] ;\:
\square_{2,2} \phi =0 \bigr\}$.
\end{rem}


\pf
Recall the integral formula for matrix coefficients (\ref{int_t}).
In the holomorphic case $m, n \ge -l$, we have $l+m \ge 0$ and
$$
(sz_{12}+z_{22})^{l+m} =
(z_{22})^{l+m} \cdot \Bigl( 1 + s \frac {z_{12}}{z_{22}} \Bigr)^{l+m} =
(z_{22})^{l+m} \cdot \sum_{q=0}^{l+m}
\begin{pmatrix} l+m \\ q \end{pmatrix}
\cdot \Bigl( s \frac {z_{12}}{z_{22}} \Bigr)^q;
$$
on the other hand $l-m <0$ and
$$
(sz_{11}+z_{21})^{l-m} = (s z_{11})^{l-m} \cdot
\Bigl( 1 + \frac {z_{21}}{sz_{11}} \Bigr)^{l-m} =
(s z_{11})^{l-m} \cdot \sum_{r=0}^{\infty}
\begin{pmatrix} l-m \\ r \end{pmatrix} \cdot
\Bigl( \frac {z_{21}}{s z_{11}} \Bigr)^r.
$$
Substituting into (\ref{int_t}) we see that $t^l_{n\,\underline{m}}(Z)$ is
the coefficient of $s^0$ in the expression
$$
\sum_{\genfrac{}{}{0pt}{}{0 \le q \le l+m}{r \ge 0}}
\begin{pmatrix} l+m \\ q \end{pmatrix} \cdot
\begin{pmatrix} l-m \\ r \end{pmatrix} \cdot
z_{11}^{l-m-r} \cdot z_{12}^q \cdot z_{21}^r \cdot z_{22}^{l+m-q} \cdot s^{n-m+q-r}.
$$
Therefore, $q = r+m-n$ and
\begin{equation}  \label{t-expression}
t^l_{n\,\underline{m}}(Z) = \sum_{r=\max \{0, n-m\}}^{l+n}
\begin{pmatrix} l+m \\ r+m-n \end{pmatrix} \cdot
\begin{pmatrix} l-m \\ r \end{pmatrix} \cdot
z_{11}^{l-m-r} \cdot z_{12}^{r+m-n} \cdot z_{21}^r \cdot z_{22}^{l+n-r}.
\end{equation}
Let $d_{ij}$ denote the degree of $z_{ij}$.
Then all the monomials in $t^l_{n\,\underline{m}}(Z)$ have
$d_{11}+d_{12}+d_{21}+d_{22}=2l \le-1$.
This proves
$$
{\cal D}^-_{\le 0} \subset
\bigl\{ \phi \in z_{11}^{-1} \cdot \BB C [z_{11}^{-1}, z_{12}, z_{21}, z_{22}]
^{\le 0} ;\: \square_{2,2} \phi =0 \bigr\}.
$$

It remains to prove the other inclusion.
Observe that each $t^l_{n\,\underline{m}}(Z)$ is a finite linear combination
of monomials $z_{11}^{d_{11}} \cdot z_{12}^{d_{12}} \cdot
z_{21}^{d_{21}} \cdot z_{22}^{d_{22}}$ with
$$
d_{11}+d_{12}=l-n, \qquad
d_{21}+d_{22}=l+n, \qquad
d_{11}+d_{21}=l-m, \qquad
d_{12}+d_{22}=l+m.
$$
These sums of degrees can be conveniently arranged in the following form:
$$
\begin{matrix}
 & d_{11}+d_{12} &  \\
d_{11}+d_{21} & & d_{12}+d_{22}  \\
 & d_{21}+d_{22} &
\end{matrix}
\quad = \quad
\begin{matrix}
 & l-n & \\
l-m & & l+m  \\
 & l+n &
\end{matrix}.
$$
Each degree invariant
\begin{equation}  \label{degree_invar}
\begin{matrix}
 & d_{11}+d_{12} &  \\
d_{11}+d_{21} & & d_{12}+d_{22}  \\
 & d_{21}+d_{22} &
\end{matrix}
\quad = \quad
\begin{matrix}
 & D_1 & \\
D_2 & & D_3  \\
 & D_4 &
\end{matrix}
\end{equation}
has an obvious restriction $D_1+D_4=D_2+D_3$.
These observations lead to the following lemma:

\begin{lem}
The matrix coefficients $t^l_{n\,\underline{m}}(Z)$ with $m,n \ge -l$
produce all possible degree invariants (\ref{degree_invar}) satisfying
$D_3,D_4 \ge 0$ and $D_1+D_4 = D_2+D_3 \le -1$.
Moreover, the functions of the type $t^l_{n\,\underline{m}}(Z)$ can be
uniquely recovered from their degree invariant.
\end{lem}

The degree invariants of $z_{11}^{-1}$, $z_{12}$, $z_{21}$, $z_{22}$ are
respectively
$$
\begin{matrix}
 & -1 & \\
-1 & & 0  \\
 & 0 &
\end{matrix}, \qquad
\begin{matrix}
 & 1 & \\
0 & & 1  \\
 & 0 &
\end{matrix}, \qquad
\begin{matrix}
 & 0 & \\
1 & & 0  \\
 & 1 &
\end{matrix}, \qquad
\begin{matrix}
 & 0 & \\
0 & & 1  \\
 & 1 &
\end{matrix}.
$$
Thus the monomials in
$z_{11}^{-1} \cdot \BB C [z_{11}^{-1}, z_{12}, z_{21}, z_{22}]$
automatically satisfy $D_3,D_4 \ge 0$.
The operator $\square_{2,2}$ is homogeneous with respect to these invariants
in the sense that it sends
$$
\begin{matrix}
 & D_1 & \\
D_2 & & D_3  \\
 & D_4 &
\end{matrix}
\quad \mapsto \quad
\begin{matrix}
 & D_1-1 & \\
D_2-1 & & D_3-1  \\
 & D_4-1 &
\end{matrix}.
$$
Thus it is enough to show that the space of harmonic polynomials in
$z_{11}^{-1} \cdot \BB C [z_{11}^{-1}, z_{12}, z_{21}, z_{22}]$
with a given degree invariant (\ref{degree_invar})
is at most one-dimensional.
Fix a degree invariant (\ref{degree_invar}),
then any function
$\phi \in z_{11}^{-1} \cdot \BB C [z_{11}^{-1}, z_{12}, z_{21}, z_{22}]$
with that invariant must be a finite linear combination of
$$
z_{12}^{D_1} \cdot z_{21}^{D_2} \cdot z_{22}^{D_3-D_1} \cdot
\Bigl( \frac{z_{12}}{z_{22}} \cdot \frac{z_{21}}{z_{11}} \Bigr)^r
=
z_{11}^{-r} \cdot z_{12}^{D_1+r} \cdot z_{21}^{D_2+r} \cdot z_{22}^{D_3-D_1-r},
\qquad r=1,\dots, D_3-D_1.
$$
Thus
$$
\phi = \sum_{r=1}^{D_3-D_1} a_r \cdot
z_{11}^{-r} \cdot z_{12}^{D_1+r} \cdot z_{21}^{D_2+r} \cdot z_{22}^{D_3-D_1-r}.
$$
We spell out the equation $\square_{2,2} \phi =0$:
\begin{multline*}
0 = \sum_{r=1}^{D_3-D_1} a_r \cdot r(D_3-D_1-r) \cdot
z_{11}^{-r-1} \cdot z_{12}^{D_1+r} \cdot z_{21}^{D_2+r} \cdot z_{22}^{D_3-D_1-r-1}  \\
+ \sum_{r=1}^{D_3-D_1} a_r \cdot (D_1+r)(D_2+r) \cdot
z_{11}^{-r} \cdot z_{12}^{D_1+r-1} \cdot z_{21}^{D_2+r-1} \cdot z_{22}^{D_3-D_1-r}
\end{multline*}
which gives us recursive equations
\begin{equation}  \label{recursive}
(r-1)(D_3-D_1-r+1) \cdot a_{r-1} = - (D_1+r)(D_2+r) \cdot a_r,
\qquad r=1,\dots, D_3-D_1.
\end{equation}
Hence the space of solutions in
$z_{11}^{-1} \cdot \BB C [z_{11}^{-1}, z_{12}, z_{21}, z_{22}]^{\le 0}$
is at most one-dimensional.

The case $m,n \le l$ is completely analogous.
\qed

\subsection{The Action of $\mathfrak{sl}(4,\BB C)$}

For convenience we restate Lemma 17 of \cite{FL} describing the $\pi^0_l$ and
$\pi^0_r$ actions of $\mathfrak{gl}(4,\BB C) = \mathfrak{gl}(2,\HC)$
on the spaces of harmonic functions
\begin{align*}
\pi_l^0 \begin{pmatrix} A & 0 \\ 0 & 0 \end{pmatrix} &:
\phi \mapsto - \tr \bigl( A \cdot X \cdot \partial \phi \bigr)  \\
\pi_r^0 \begin{pmatrix} A & 0 \\ 0 & 0 \end{pmatrix} &:
\phi \mapsto - \tr \bigl( A \cdot (X \cdot \partial \phi + \phi) \bigr)  \\
\pi_l^0 \begin{pmatrix} 0 & B \\ 0 & 0 \end{pmatrix} =
\pi_r^0 \begin{pmatrix} 0 & B \\ 0 & 0 \end{pmatrix} &:
\phi \mapsto - \tr \bigl( B \cdot \partial \phi \bigr)  \\
\pi_l^0 \begin{pmatrix} 0 & 0 \\ C & 0 \end{pmatrix} =
\pi_r^0 \begin{pmatrix} 0 & 0 \\ C & 0 \end{pmatrix} &:
\phi \mapsto \tr \Bigl( C \cdot \bigl(
X \cdot (\partial \phi) \cdot X + X\phi \bigr) \Bigr)
= \tr \Bigl( C \cdot \bigl(
X \cdot \partial (X\phi) \bigr) - X\phi \Bigr)  \\
\pi_l^0 \begin{pmatrix} 0 & 0 \\ 0 & D \end{pmatrix} &: \phi \mapsto
\tr \Bigl( D \cdot \bigl( (\partial \phi) \cdot X + \phi \bigr) \Bigr)
= \tr \Bigl( D \cdot \bigl( \partial (X\phi) - \phi \bigr) \Bigr) \\
\pi_r^0 \begin{pmatrix} 0 & 0 \\ 0 & D \end{pmatrix} &: \phi \mapsto
\tr \bigl( D \cdot (\partial \phi) \cdot X \bigr)
= \tr \Bigl( D \cdot \bigl( \partial (X\phi) - 2\phi \bigr) \Bigr),
\end{align*}
where $\partial = \begin{pmatrix} \partial_{11} & \partial_{21} \\
\partial_{12} & \partial_{22} \end{pmatrix} = \frac 12 \nabla$.
Recall the representations ${\cal H}^+$ and ${\cal H}^-$ of
$\mathfrak{sl}(2,\HC)$ realized in the space of harmonic functions on $\BB H$
(see Subsection 2.5 in \cite{FL}).

\begin{thm}  \label{D-component}
We have ${\cal D}^- \simeq {\cal H}^-$ and ${\cal D}^+ \simeq {\cal H}^+$
as representations of 
$\mathfrak{sl}(4,\BB C) = \mathfrak{sl}(2,\HC) \subset \mathfrak{gl}(2,\HC)$.
Moreover, ${\cal D}^-$ and ${\cal D}^+$ are irreducible representations of
$\mathfrak{sl}(4,\BB C)$ with highest (or lowest) weight vectors
$t^{-\frac12}_{\frac12\,\underline{\frac12}}(Z) = \frac1{z_{11}}$ and
$t^{-\frac12}_{-\frac12\,\underline{-\frac12}}(Z) = \frac1{z_{22}}$ respectively.
\end{thm}

\pf
Using Lemmas \ref{deriv_calc} and \ref{mult-ids} we can compute the
Lie algebra actions $\pi^0_l$ and $\pi^0_r$ of $\mathfrak{gl}(2,\BB H)$ on
${\cal D}^-$ and ${\cal D}^+$:
\begin{align*}
\pi_l^0 \begin{pmatrix} A & 0 \\ 0 & 0 \end{pmatrix} &: t^{l}_{n\,\underline{m}}
\mapsto -\tr \biggl( A \cdot \begin{pmatrix}
(l-n) t^{l}_{n\,\underline{m}} & (l-n+1) t^{l}_{n-1\,\underline{m}} \\
(l+n+1) t^{l}_{n+1\,\underline{m}} & (l+n) t^{l}_{n\,\underline{m}}
\end{pmatrix} \biggr)  \\
\pi_r^0 \begin{pmatrix} A & 0 \\ 0 & 0 \end{pmatrix} &: t^{l}_{n\,\underline{m}}
\mapsto - \tr \biggl( A \cdot \begin{pmatrix}
(l-n+1) t^{l}_{n\,\underline{m}} & (l-n+1) t^{l}_{n-1\,\underline{m}} \\
(l+n+1) t^{l}_{n+1\,\underline{m}} & (l+n+1) t^{l}_{n\,\underline{m}}
\end{pmatrix} \biggr)  \\
\pi_l^0 \begin{pmatrix} 0 & B \\ 0 & 0 \end{pmatrix} =
\pi_r^0 \begin{pmatrix} 0 & B \\ 0 & 0 \end{pmatrix} &: t^{l}_{n\,\underline{m}}
\mapsto - \tr \Biggl( B \cdot \begin{pmatrix}
(l-m) t^{l-\frac12}_{n+\frac12\,\underline{m+\frac12}} &
(l-m) t^{l-\frac12}_{n-\frac12\,\underline{m+\frac12}} \\
(l+m) t^{l-\frac12}_{n+\frac12\,\underline{m-\frac12}} &
(l+m) t^{l-\frac12}_{n-\frac12\,\underline{m-\frac12}}
\end{pmatrix} \Biggr)  \\
\pi_l^0 \begin{pmatrix} 0 & 0 \\ C & 0 \end{pmatrix} =
\pi_r^0 \begin{pmatrix} 0 & 0 \\ C & 0 \end{pmatrix} &: t^{l}_{n\,\underline{m}}
\mapsto \tr \Biggl( C \cdot \begin{pmatrix}
(l-n+1) t^{l+\frac12}_{n-\frac12\,\underline{m-\frac12}} &
(l-n+1) t^{l+\frac12}_{n-\frac12\,\underline{m+\frac12}} \\
(l+n+1) t^{l+\frac12}_{n+\frac12\,\underline{m-\frac12}} &
(l+n+1) t^{l+\frac12}_{n+\frac12\,\underline{m+\frac12}}
\end{pmatrix} \Biggr)  \\
\pi_l^0 \begin{pmatrix} 0 & 0 \\ 0 & D \end{pmatrix} &: t^{l}_{n\,\underline{m}}
\mapsto \tr \Biggl( D \cdot \begin{pmatrix}
(l-m+1) t^{l}_{n\,\underline{m}} & (l-m) t^{l}_{n\,\underline{m+1}} \\
(l+m) t^{l}_{n\,\underline{m-1}} & (l+m+1) t^{l}_{n\,\underline{m}}
\end{pmatrix} \Biggr)  \\
\pi_r^0 \begin{pmatrix} 0 & 0 \\ 0 & D \end{pmatrix} &: t^{l}_{n\,\underline{m}}
\mapsto \tr \Biggl( D \cdot \begin{pmatrix}
(l-m) t^{l}_{n\,\underline{m}} & (l-m) t^{l}_{n\,\underline{m+1}} \\
(l+m) t^{l}_{n\,\underline{m-1}} & (l+m) t^{l}_{n\,\underline{m}}
\end{pmatrix} \Biggr)
\end{align*}
and similarly for $t^l_{n\,\underline{m}}(Z) \cdot N(Z)^{-2l-1}$'s.
This together with equation (\ref{t=0}) shows that the $\pi_l^0$-action of
$\mathfrak{sl}(2,\HC)$ preserves ${\cal D}^{\mp}$ and that
$t^{-\frac12}_{\frac12\,\underline{\frac12}}(Z)$ and
$t^{-\frac12}_{-\frac12\,\underline{-\frac12}}(Z)$
generate ${\cal D}^-$ and ${\cal D}^+$ respectively.
Let $\tilde e_3=ie_3$. Consider an element
\begin{equation}  \label{Cayley-composition}
\gamma_0 = \frac12 \begin{pmatrix} -i\tilde e_3 -i & \tilde e_3-1 \\
\tilde e_3 -1 & i\tilde e_3 +i \end{pmatrix} \in GL(2,\HC)
\qquad \text{with} \qquad
\gamma_0^{-1} = \frac12 \begin{pmatrix} i\tilde e_3 +i & \tilde e_3-1 \\
\tilde e_3 -1 & -i\tilde e_3 -i \end{pmatrix}.
\end{equation}
Then
\begin{align*}
\pi_l^0(\gamma_0) : \qquad {\cal H}^+ \ni 1 \quad &\mapsto \quad
\frac1{iz_{22}}=-i t^{-\frac12}_{-\frac12\,\underline{-\frac12}}(Z) \in {\cal D}^+, \\
{\cal H}^- \ni \frac1{N(Z)} \quad &\mapsto \quad
\frac{i}{z_{11}}=i t^{-\frac12}_{\frac12\,\underline{\frac12}}(Z) \in {\cal D}^-.
\end{align*}
(This is essentially the composition of the Cayley transform from \cite{FL}
with another Cayley-type transform which will be introduced in
Proposition \ref{Cayley}.)
This proves ${\cal D}^{\mp} \simeq {\cal H}^{\mp}$, irreducibility and the
statement about the highest (or lowest) weight vectors.
\qed

Representations ${\cal D}^-$ and ${\cal D}^+$ are dual to each other.
Define a bilinear pairing between ${\cal D}^-$ and ${\cal D}^+$ by declaring
$$
\Bigl\langle t^{l'}_{n'\,\underline{m'}}(Z), \frac 1{N(Z)} \cdot
t^l_{m\,\underline{n}}(Z^{-1}) \Bigr\rangle_{{\cal D}^- \times {\cal D}^+}
= \Bigl\langle \frac 1{N(Z)} \cdot t^l_{m\,\underline{n}}(Z^{-1}),
t^{l'}_{n'\,\underline{m'}}(Z) \Bigr\rangle_{{\cal D}^- \times {\cal D}^+} =
\delta_{ll'} \delta_{mm'} \delta_{nn'},
$$
$$
\bigl\langle t^{l}_{n\,\underline{m}}(Z),
t^{l'}_{n'\,\underline{m'}}(Z) \bigr\rangle_{{\cal D}^- \times {\cal D}^+}
=
\Bigl\langle \frac 1{N(Z)} \cdot t^l_{m\,\underline{n}}(Z^{-1}), \frac 1{N(Z)}
\cdot t^{l'}_{m'\,\underline{n'}}(Z^{-1}) \Bigr\rangle_{{\cal D}^- \times {\cal D}^+} =0,
$$
$l = -1, -\frac32, -2, \dots$.
In the second line we exclude $l=-1/2$ because by (\ref{t(Z^{-1})}) we have
$$
N(Z)^{-1} \cdot t^{-\frac12}_{n\,\underline{m}} \bigl( Z^{-1} \bigr)
= t^{-\frac12}_{-m\,\underline{-n}}(Z).
$$
By Proposition \ref{t-orthogonal}, this pairing partially agrees with the
bilinear form (\ref{b-form}) up to a sign.

\begin{prop}  \label{invar-pairing}
This bilinear pairing on ${\cal D}^- \times {\cal D}^+$
is $\mathfrak{gl}(2,\HC)$-invariant:
$$
\bigl\langle \pi^0_l(Z) \phi_1, \phi_2 \bigr\rangle_{{\cal D}^- \times {\cal D}^+}
+ \bigl\langle \phi_1, \pi^0_r(Z) \phi_2 \bigr\rangle_{{\cal D}^- \times {\cal D}^+} =0,
\qquad \forall Z \in \mathfrak{gl}(2,\HC).
$$
\end{prop}

\pf
Since the elements
$\begin{pmatrix} 0 & B \\ 0 & 0 \end{pmatrix} \in \mathfrak{sl}(2,\HC)$,
$B \in \HC$, together with their conjugates by
$\begin{pmatrix} 0 & 1 \\ 1 & 0 \end{pmatrix} \in GL(2,\HC)$
generate $\mathfrak{sl}(2,\HC)$, to prove $\mathfrak{sl}(2,\HC)$-invariance
it is enough to check the invariance under the action of
$\begin{pmatrix} 0 & B \\ 0 & 0 \end{pmatrix}$. By Lemma 17 from \cite{FL},
$$
\pi_l^0 \begin{pmatrix} 0 & B \\ 0 & 0 \end{pmatrix} =
\pi_r^0 \begin{pmatrix} 0 & B \\ 0 & 0 \end{pmatrix} :
\phi \mapsto \tr \bigl( B \cdot (-\partial \phi ) \bigr).
$$
Applying Lemmas \ref{deriv_calc} and \ref{mult-ids} repeatedly we find that
$$
\begin{pmatrix} \partial_{11} & \partial_{21} \\ \partial_{12} & \partial_{22}
\end{pmatrix} t^{l}_{n \,\underline{m}}(Z)
=
\begin{pmatrix}
(l-m) t^{l- \frac 12}_{n+ \frac 12 \,\underline{m+ \frac 12}}(Z) &
(l-m) t^{l- \frac 12}_{n- \frac 12 \,\underline{m+ \frac 12}}(Z) \\
(l+m) t^{l- \frac 12}_{n+ \frac 12 \,\underline{m- \frac 12}}(Z) &
(l+m) t^{l- \frac 12}_{n- \frac 12 \,\underline{m- \frac 12}}(Z)
\end{pmatrix},
$$
\begin{multline*}
\begin{pmatrix} \partial_{11} & \partial_{21} \\ \partial_{12} & \partial_{22}
\end{pmatrix} \Bigl( \frac 1{N(Z)} \cdot t^l_{m\,\underline{n}}(Z^{-1}) \Bigr)  \\
= - \frac 1{N(Z)} \begin{pmatrix}
(l-m+1) t^{l+\frac 12}_{m-\frac12\,\underline{n-\frac12}}(Z^{-1}) &
(l-m+1) t^{l+\frac12}_{m-\frac12\,\underline{n+\frac12}}(Z^{-1}) \\
(l+m+1) t^{l+\frac12}_{m+\frac12 \,\underline{n-\frac12}}(Z^{-1}) &
(l+m+1) t^{l+\frac12}_{m+\frac12\,\underline{n+\frac12}}(Z^{-1})
\end{pmatrix}.
\end{multline*}
We conclude that the bilinear pairing is invariant under the action of
$\begin{pmatrix} 0 & B \\ 0 & 0 \end{pmatrix}$.
Using Lemma 17 from \cite{FL} it is easy to see that the pairing is invariant
under the scalar matrices as well, and the proposition follows.
\qed

\subsection{Ol'shanskii Semigroups}

We introduce complex Ol'shanskii semigroups
$\Gamma^{\pm} \subset GL(2,\BB C) \subset \HC$.
Following \cite{KouO} and \cite{HN}, consider a Hermitian form $H$
on $\BB C^2$ defined by
$$
H(\zeta, \eta) = - \zeta_1 \overline{\eta_1} + \zeta_2 \overline{\eta_2}.
$$
It is easy to check that if $X \in i \mathfrak{u}(1,1)$, then
$H(X\zeta, \zeta) \in \BB R$ for all $\zeta \in \BB C^2$.
Consider a cone in $i \mathfrak{u}(1,1)$ defined by
$$
C = \{ X \in i \mathfrak{u}(1,1) ;\: H(X\zeta, \zeta) \le 0, \:
\forall \zeta \in \BB C^2 \}.
$$

\begin{lem}
The cone $C$ is closed, convex, pointed (i.e. $C \cap -C = \{0\}$),
generating (i.e. $C -C = i \mathfrak{u}(1,1)$ or, equivalently,
has non-empty interior),
hyperbolic (i.e. for every $X \in C$, the operator $ad X$ has real eigenvalues
and, for every $X$ in the interior of $C$, $ad X$ is diagonalizable)
and $Ad(U(1,1))$-invariant.
\end{lem}

The set
$$
\bar\Gamma^- = U(1,1) \cdot \exp (C)
$$
is called a closed complex Ol'shanskii semigroup contained in $GL(2,\BB C)$.
The interior of $C$ is
$$
C^0 = \{ X \in i \mathfrak{u}(1,1) ;\: H(X\zeta, \zeta) < 0, \:
\forall \zeta \in \BB C^2 \setminus \{0\} \},
$$
and the corresponding open Ol'shanskii semigroup is
$$
\Gamma^- = U(1,1) \cdot \exp (C^0).
$$
The semigroup $\Gamma^-$ is an open subset of $GL(2,\BB C)$ and hence of $\HC$.
In fact, $\Gamma^-$ is the interior of $\bar\Gamma^-$.
If we re-define for a moment the complex conjugation on $\HC$ to be relative
to $\HR = \mathfrak{u}(1,1)$ (i.e. by identifying $\HC$ with
$\mathfrak{u}(1,1) \oplus i \mathfrak{u}(1,1)$),
then the conjugate Ol'shanskii semigroup $\overline{\bar\Gamma^-}$ coincides
with $\bar\Gamma^+ = U(1,1) \cdot \exp (-C)$.
Similarly, $\overline{\Gamma^-} = \Gamma^+ = U(1,1) \cdot \exp (-C^0)$.

\begin{lem}
We have
$$
C = \biggl\{ X =
\begin{pmatrix} a & \gamma \\ -\overline{\gamma} & -b \end{pmatrix};\:
a, b \in \BB R,\: a,b \ge 0,\: \gamma \in \BB C, \: |\gamma|^2 \le ab \biggr\}
\quad \subset \quad i \mathfrak{u}(1,1)
$$
and
$$
C^0 = \biggl\{ X =
\begin{pmatrix} a & \gamma \\ -\overline{\gamma} & -b \end{pmatrix};\:
a, b \in \BB R,\: a,b>0,\: \gamma \in \BB C, \: |\gamma|^2 < ab \biggr\}
\quad \subset \quad i \mathfrak{u}(1,1).
$$
In particular, each $X \in C^0$ has two distinct real eigenvalues --
one positive and one negative.
\end{lem}

Fix a maximal compact subgroup $U(1) \times U(1)$ of $U(1,1)$;
its complexification $K_{\BB C}$ consists of diagonal matrices, and
every element in $\Gamma^-$ is $SU(1,1)$-conjugate to an element in
$$
\Gamma^- \cap K_{\BB C} =
\biggl\{ \begin{pmatrix} \lambda_1 & 0 \\ 0 & \lambda_2 \end{pmatrix} ;\:
\lambda_1, \lambda_2 \in \BB C,\: |\lambda_1| > 1 > |\lambda_2| > 0 \biggr\}.
$$
Similarly, every element in $\Gamma^+$ is $SU(1,1)$-conjugate to an element in
$$
\Gamma^+ \cap K_{\BB C} =
\biggl\{ \begin{pmatrix} \lambda_1 & 0 \\ 0 & \lambda_2 \end{pmatrix} ;\:
\lambda_1, \lambda_2 \in \BB C,\: |\lambda_2| > 1 > |\lambda_1| > 0 \biggr\}.
$$


\begin{prop}
The Ol'shanskii semigroups $\bar\Gamma^-$ and $\Gamma^-$ have the following
crucial property:
$$
\text{$\gamma_0 \gamma$ and $\gamma \gamma_0 \in \Gamma^-$
\qquad whenever $\gamma_0 \in \Gamma^-$ and $\gamma \in \bar\Gamma^-$}.
$$
\end{prop}

In particular, if we let
$$
W \in \HR^+
\qquad \text{and} \qquad Z \in \sqrt{N(W)} \cdot
\begin{pmatrix} \sigma & 0 \\ 0 & \sigma^{-1} \end{pmatrix}
\cdot SU(1,1)\quad \subset \HC,
$$
then $WZ^{-1} \in \Gamma^-$ for all $1> \sigma >0$ and
$WZ^{-1} \in \Gamma^+$ for all $\sigma >1$.

\begin{lem} \label{1/N-bound}
For $W \in \Gamma^- \cup \Gamma^+$ and $Z \in U(1,1)$, the function $N(Z-W)$
is never zero.
Moreover, if we fix $W \in \Gamma^- \cup \Gamma^+$, there exist $c>0$ and
$\epsilon >0$ depending on $W$ such that
$$
\frac 1{N(Z-W)} \le \frac c{\|Z\|}
$$
for all $Z \in \HR$ with $1 - \epsilon \le N(Z) \le 1+ \epsilon$.
\end{lem}

\pf
For concreteness, let us suppose $W \in \Gamma^-$; the other case is similar.
We rewrite $N(Z-W)$ as $N(Z^{-1}W -1) \cdot N(Z)$. Then $Z^{-1}W \in \Gamma^-$
and hence has eigenvalues different from $1$. Therefore, $N(Z^{-1}W -1) \ne 0$.

If $Z \in \HR$, then $N(Z-W) = N(Z) \cdot N(Z^{-1}W -1)$.
Every element $W \in \Gamma^-$ is $SU(1,1)$-conjugate to a diagonal matrix
$\begin{pmatrix} \lambda_1 & 0 \\ 0 & \lambda_2 \end{pmatrix}$ with
$\lambda_1, \lambda_2 \in \BB C$ and $|\lambda_1| > 1 > |\lambda_2| > 0$.
Since the determinant function is $Ad(SU(1,1))$-invariant,
without loss of generality we may assume that
$W=\begin{pmatrix} \lambda_1 & 0 \\ 0 & \lambda_2 \end{pmatrix}$.
Then
\begin{multline*}
\det \biggl[
\begin{pmatrix} \overline{z_{11}} & -z_{12} \\ -\overline{z_{12}} & z_{11}
\end{pmatrix} \begin{pmatrix} \lambda_1 & 0 \\ 0 & \lambda_2 \end{pmatrix}
- \begin{pmatrix} 1 & 0 \\ 0 & 1 \end{pmatrix} \biggr]
=
\det \begin{pmatrix} \lambda_1\overline{z_{11}} -1 & -\lambda_2 z_{12} \\
-\lambda_1 \overline{z_{12}} & \lambda_2 z_{11} -1 \end{pmatrix}  \\
= (\lambda_1\overline{z_{11}} -1)(\lambda_2 z_{11} -1)
- \lambda_1\lambda_2|z_{12}|^2
= \lambda_1\lambda_2N(Z) + 1 - (\lambda_1\overline{z_{11}} + \lambda_2 z_{11}).
\end{multline*}
Now the term $\lambda_1\lambda_2N(Z)+1$ stays bounded while
$\| \lambda_1\overline{z_{11}} + \lambda_2 z_{11} \|$ grows proportionally to
$\|Z\|$ as $Z \to \infty$.
\qed

As we saw in Subsection \ref{polynomial_algebra}, the matrix coefficients
of the discrete series and their limits are rational functions in $z_{ij}$'s.
More precisely, the matrix coefficients of the holomorphic discrete series
and its limit lie in $\BB C [ z_{11}, z_{12}, z_{21}, z_{22}, z_{11}^{-1}]$.
Hence they extend uniquely as holomorphic functions to the open set
$\{ Z \in \HC ;\: z_{11} \ne 0 \}$.
On the other hand, Lemma 2.6 from \cite{KouO} implies that
$\Gamma^- \subset \{ Z \in \HC ;\: z_{11} \ne 0 \}$.
In particular, the matrix coefficients of the holomorphic discrete series
and its limit extend holomorphically to $\Gamma^-$.
Similarly, the matrix coefficients of the antiholomorphic discrete series
and its limit lie in $\BB C [ z_{11}, z_{12}, z_{21}, z_{22}, z_{22}^{-1}]$,
extend holomorphically to the open set $\{ Z \in \HC ;\: z_{22} \ne 0 \}$.
Since $\Gamma^+ \subset \{ Z \in \HC ;\: z_{22} \ne 0 \}$, these coefficients
extend to $\Gamma^+$.

\subsection{Matrix coefficient expansions for $\frac 1{N(Z-W)}$
and $\frac {(Z-W)^{-1}}{N(Z-W)}$}

In this subsection we derive matrix coefficient expansions
for $\frac 1{N(Z-W)}$ and $\frac {(Z-W)^{-1}}{N(Z-W)}$.
Using these expansions we obtain projectors onto the discrete series
components for the spaces of solutions of $\square_{2,2}\phi=0$ and
left- and right-regular functions on $\HR^+$.

\begin{prop}  \label{expansion2}
We have the following matrix coefficient expansions
$$
- \frac 1{N(Z-W)} =
\sum_{\genfrac{}{}{0pt}{}{l,m,n}{m,n \ge -l \ge 1}}
t^l_{n\underline{m}}(W) \cdot \frac 1{N(Z)} \cdot t^l_{m\underline{n}}(Z^{-1})
+ \sum_{\genfrac{}{}{0pt}{}{l,m,n}{m,n \le l \le -1/2}}
\frac 1{N(W)} \cdot t^l_{n\underline{m}}(W^{-1}) \cdot t^l_{m\underline{n}}(Z)
$$
which converges pointwise absolutely whenever $WZ^{-1} \in \Gamma^-$ and
$$
- \frac 1{N(Z-W)} =
\sum_{\genfrac{}{}{0pt}{}{l,m,n}{m,n \le l \le -1}}
t^l_{n\underline{m}}(W) \cdot \frac 1{N(Z)} \cdot t^l_{m\underline{n}}(Z^{-1})  \\
+ \sum_{\genfrac{}{}{0pt}{}{l,m,n}{m,n \ge -l \ge 1/2}}
\frac 1{N(W)} \cdot t^l_{n\underline{m}}(W^{-1}) \cdot t^l_{m\underline{n}}(Z)
$$
which converges pointwise absolutely whenever $WZ^{-1} \in \Gamma^+$.
\end{prop}

\begin{rem}
From (\ref{t(Z^{-1})}) we have
$$
N(W)^{-1} \cdot t^{-\frac12}_{n\underline{m}}(W^{-1}) \cdot t^{-\frac12}_{m\underline{n}}(Z)
= t^{-\frac12}_{-m\underline{-n}}(W) \cdot t^{-\frac12}_{n\underline{m}}(Z) =
t^{-\frac12}_{-m\underline{-n}}(W) \cdot N(Z)^{-1} \cdot
t^{-\frac12}_{-n\underline{-m}}(Z^{-1}).
$$
Hence, in each expansion of $\frac 1{N(Z-W)}$, the terms corresponding to the
limits of the discrete series (terms with $l=-1/2$) can be included in either
the first or the second sum.
\end{rem}

\pf
Since every element in $\Gamma^-$ is $SU(1,1)$-conjugate to an element in
$\Gamma^- \cap K_{\BB C}$, to prove the first expansion we can assume that
$WZ^{-1} \in \Gamma^- \cap K_{\BB C}$, and so
$WZ^{-1} = \begin{pmatrix} \lambda_1 & 0 \\ 0 & \lambda_2 \end{pmatrix}$
for some $\lambda_1, \lambda_2 \in \BB C$ with
$|\lambda_1| > 1 > |\lambda_2| > 0$.
The discrete series characters are:
$$
\tilde\Theta_l^+(WZ^{-1}) = - \frac{\lambda_2^{2l+1}}{\lambda_1-\lambda_2},
\qquad
\tilde\Theta_l^+(ZW^{-1}) =
\frac{\lambda_1 \cdot \lambda_2^{-2l}}{\lambda_1-\lambda_2},
$$
\begin{equation}  \label{Theta}
\tilde\Theta_l^-(WZ^{-1}) = \frac{\lambda_1^{2l+1}}{\lambda_1-\lambda_2},
\qquad
\tilde\Theta_l^-(ZW^{-1}) =
-\frac{\lambda_2 \cdot \lambda_1^{-2l}}{\lambda_1-\lambda_2}.
\end{equation}
Using the multiplicativity property of matrix coefficients
\begin{equation}  \label{t-mult}
t^l_{n\underline{m}}(Z_1Z_2) =
\sum_k t^l_{n\underline{k}}(Z_1) \cdot t^l_{k\underline{m}}(Z_2),
\end{equation}
the expansion reduces to a geometric series computation:
\begin{multline*}
\sum_{\genfrac{}{}{0pt}{}{l,m,n}{m,n \ge -l \ge 1}}
t^l_{n\underline{m}}(W) \cdot \frac 1{N(Z)} \cdot t^l_{m\underline{n}}(Z^{-1})
+ \sum_{\genfrac{}{}{0pt}{}{l,m,n}{m,n \le l \le -1/2}}
\frac 1{N(W)} \cdot t^l_{n\underline{m}}(W^{-1}) \cdot t^l_{m\underline{n}}(Z)  \\
= \frac 1{N(Z)} \cdot \biggl( \sum_{\genfrac{}{}{0pt}{}{l,n}{n \ge -l \ge 1}}
t^l_{n\underline{n}}(WZ^{-1})
+ N(ZW^{-1}) \cdot \sum_{\genfrac{}{}{0pt}{}{l,n}{n \le l \le -1/2}}
t^l_{n\underline{n}}(ZW^{-1}) \biggr)  \\
= \frac 1{N(Z)} \cdot \biggl( \sum_{l \le -1} \tilde\Theta_l^- (WZ^{-1})
+ N(ZW^{-1}) \cdot \sum_{l \le -1/2} \tilde\Theta_l^+(ZW^{-1}) \biggr)  \\
= \frac 1{N(Z)} \cdot \biggl( \sum_{l \le -1}
\frac{\lambda_1^{2l+1}}{\lambda_1 - \lambda_2}
+ \sum_{l \le -1/2} \frac{\lambda_2^{-(2l+1)}}{\lambda_1 - \lambda_2} \biggr)
=
\frac 1{N(Z)} \cdot \frac 1{\lambda_1 - \lambda_2} \cdot
\biggl(\frac 1{\lambda_1-1} + \frac 1{1-\lambda_2} \biggr)  \\
=
- \frac 1{N(Z)} \cdot \frac 1{(\lambda_1-1)(\lambda_2-1)}
=
- \frac 1{N(Z)} \cdot \frac 1{N(WZ^{-1} -1)}
= - \frac 1{N(Z-W)}.
\end{multline*}
The other matrix coefficient expansion is proved in the same way.
\qed

For $R>0$, define operators on ${\cal H}(\HR^+)$ by
\begin{align*}
\bigl( \S_R^- \phi \bigr)(W) &=
-\frac1{2\pi^2} \int_{X \in H_R}
\frac {(\deg \phi)(X)}{N(X-W)} \cdot \frac {dS}{\|X\|}
= \Bigl\langle \phi(X), \frac 1{N(X-W)} \Bigr\rangle_R,
\qquad W \in R \cdot \Gamma^-,  \\
\bigl( \S_R^+ \phi \bigr)(W) &=
-\frac1{2\pi^2} \int_{X \in H_R}
\frac {(\deg \phi)(X)}{N(X-W)} \cdot \frac {dS}{\|X\|}
= \Bigl\langle \phi(X), \frac 1{N(X-W)} \Bigr\rangle_R,
\qquad W \in R \cdot \Gamma^+.
\end{align*}
By Lemma \ref{1/N-bound}, these integrals are well defined.
From the orthogonality relations for matrix coefficients
(Proposition \ref{t-orthogonal}) we obtain:

\begin{thm}  \label{discrete_ser_proj}
The operators $\S_R^-$ and $\S_R^+$ are continuous linear operators
${\cal H}(\HR^+) \to {\cal H}(\HR^+)$.
The operator $\S_R^-$ annihilates the continuous series,
the antiholomorphic discrete series ${\cal D}_{discr}^+$ and sends
$$
\begin{matrix}
t^l_{n\,\underline{m}}(X) \quad \mapsto \quad - t^l_{n\,\underline{m}}(W)
- R^{2(2l+1)} \cdot N(W)^{-2l-1} \cdot t^l_{n\,\underline{m}}(W), \\
\quad \\ N(X)^{-2l-1} \cdot t^l_{n\,\underline{m}}(X) \quad \mapsto \quad
R^{-2(2l+1)} \cdot t^l_{n\,\underline{m}}(W) + N(W)^{-2l-1} \cdot t^l_{n\,\underline{m}}(W),
\end{matrix}
\quad
\begin{matrix}
l = -1, -\frac 32, -2, \dots, \\ m,n \in \BB Z+l,\\ m, n \ge -l.
\end{matrix}
$$
The operator $\S_R^+$ annihilates the continuous series,
the holomorphic discrete series ${\cal D}_{discr}^-$ and sends
$$
\begin{matrix}
t^l_{n\,\underline{m}}(X) \quad \mapsto \quad - t^l_{n\,\underline{m}}(W)
- R^{2(2l+1)} \cdot N(W)^{-2l-1} \cdot t^l_{n\,\underline{m}}(W),\\
\quad \\ N(X)^{-2l-1} \cdot t^l_{n\,\underline{m}}(X) \quad \mapsto \quad
R^{-2(2l+1)} \cdot t^l_{n\,\underline{m}}(W) + N(W)^{-2l-1} \cdot t^l_{n\,\underline{m}}(W),
\end{matrix}
\quad
\begin{matrix}
l = -1, -\frac 32, -2, \dots, \\ m,n \in \BB Z+l,\\ m, n \le l.
\end{matrix}
$$
\end{thm}

Note that the closure of $R \cdot \Gamma^-$ is $R \cdot \Gamma$
which contains $H_R$. Also, the functions $t^l_{n\,\underline{m}}(W)$ and
$R^{2(2l+1)} \cdot N(W)^{-2l-1} \cdot t^l_{n\,\underline{m}}(W)$ agree on $H_R$.
Thus, the values of the holomorphic discrete series component of a function
$\phi \in {\cal H}(\HR^+)$ on $H_R$ can be determined by continuity.
Similarly, one can recover the antiholomorphic discrete series component
of $\phi$ on $H_R$.

Differentiating the matrix coefficient expansions for $\frac 1{N(Z-W)}$
we get two expansions for $\frac {(Z-W)^{-1}}{N(Z-W)}$:

\begin{prop}
We have the following matrix coefficient expansions:
\begin{multline*}
\frac {(Z-W)^{-1}}{N(Z-W)} =
\sum_{\genfrac{}{}{0pt}{}{l,m,n}{m,n \le l \le -1/2}}
\frac 1{N(W)}
\begin{pmatrix}
(l-n+1) t^l_{n-1\, \underline{m}} (W^{-1})  \\
(l+n) t^l_{n \, \underline{m}} (W^{-1})
\end{pmatrix}
\Bigl( t^{l- \frac 12}_{m+\frac12 \, \underline{n-\frac12}}(Z),
t^{l- \frac 12}_{m-\frac12 \, \underline{n-\frac12}}(Z) \Bigr)  \\
- \sum_{\genfrac{}{}{0pt}{}{l,m,n}{m,n \ge -l \ge 1}}
\begin{pmatrix}
(l-m+1) t^{l}_{n \, \underline{m}}(W)  \\
(l+m) t^{l}_{n \, \underline{m-1}}(W)
\end{pmatrix}
\Bigl( N(Z)^{-1} \cdot t^{l+\frac12}_{m-\frac12 \, \underline{n- \frac 12}}(Z^{-1}),
N(Z)^{-1} \cdot t^{l+\frac12}_{m-\frac12 \, \underline{n+ \frac 12}}(Z^{-1}) \Bigr)
\end{multline*}
which converges pointwise absolutely whenever $WZ^{-1} \in \Gamma^-$ and
\begin{multline*}
\frac {(Z-W)^{-1}}{N(Z-W)} =
\sum_{\genfrac{}{}{0pt}{}{l,m,n}{m,n \ge -l \ge 1/2}}
\frac 1{N(W)}
\begin{pmatrix}
(l-n) t^l_{n\, \underline{m}} (W^{-1})  \\
(l+n+1) t^l_{n+1 \, \underline{m}} (W^{-1})
\end{pmatrix}
\Bigl( t^{l- \frac 12}_{m+\frac12 \, \underline{n+\frac12}}(Z),
t^{l- \frac 12}_{m-\frac12 \, \underline{n+\frac12}}(Z) \Bigr)  \\
- \sum_{\genfrac{}{}{0pt}{}{l,m,n}{m,n \le l \le -1}}
\begin{pmatrix}
(l-m) t^{l}_{n \, \underline{m+1}}(W)  \\
(l+m+1) t^{l}_{n \, \underline{m}}(W)
\end{pmatrix}
\Bigl( N(Z)^{-1} \cdot t^{l+\frac12}_{m+\frac12 \, \underline{n- \frac 12}}(Z^{-1}),
N(Z)^{-1} \cdot t^{l+\frac12}_{m+\frac12 \, \underline{n+ \frac 12}}(Z^{-1}) \Bigr)
\end{multline*}
which converges pointwise absolutely whenever $WZ^{-1} \in \Gamma^+$.
\end{prop}

\pf
The two expansions are proved by applying $\nabla_Z$ to the expansions
in Proposition \ref{expansion2}, where the subscript $Z$ in $\nabla_Z$
indicates that the differentiation is done with respect to this variable.
We give a proof of the second formula
\begin{multline}  \label{exp-sum1+2}
\frac {(Z-W)^{-1}}{N(Z-W)} = - \frac 12 \nabla_Z \frac 1{N(Z-W)}
=
\frac 1{N(W)} \sum_{\genfrac{}{}{0pt}{}{l,m,n'}{m,n' \ge -l \ge 1/2}}
t^{l}_{n'\underline{m}}(W^{-1}) \cdot
\begin{pmatrix} \partial_{11} & \partial_{21} \\
\partial_{12} & \partial_{22} \end{pmatrix}_Z t^{l}_{m\underline{n'}}(Z)  \\
-
\sum_{\genfrac{}{}{0pt}{}{l,m',n}{m',n \le l \le -1}}
t^{l}_{n\underline{m'}}(W) \cdot
\begin{pmatrix} \partial_{11} & \partial_{21} \\
\partial_{12} & \partial_{22} \end{pmatrix}_Z
\biggl( \frac {-1}{N(Z)} \cdot t^{l}_{m'\underline{n}}(Z^{-1}) \biggr).
\end{multline}
Using the derivative formulas from Lemma \ref{deriv_calc}
we can expand the first sum in (\ref{exp-sum1+2}) as:
\begin{multline*}
\sum_{\genfrac{}{}{0pt}{}{l,m,n'}{m,n' \ge -l \ge 1/2}}
t^{l}_{n'\underline{m}} (W^{-1}) \cdot
\begin{pmatrix}
(l-n') t^{l- \frac 12}_{m+ \frac 12 \, \underline{n'+ \frac 12}}(Z) &
(l-n') t^{l- \frac 12}_{m- \frac 12 \, \underline{n'+ \frac 12}}(Z)  \\
(l+n') t^{l- \frac 12}_{m+ \frac 12 \, \underline{n'- \frac 12}}(Z) &
(l+n') t^{l- \frac 12}_{m- \frac 12 \, \underline{n'- \frac 12}}(Z)
\end{pmatrix}
=  \\
\sum_{\genfrac{}{}{0pt}{}{l,m,n'}{m,n' \ge -l \ge 1/2}}
\begin{pmatrix}
(l-n') t^{l}_{n'\underline{m}} (W^{-1}) \cdot
t^{l- \frac 12}_{m+ \frac 12 \, \underline{n'+ \frac 12}}(Z) &
(l-n') t^{l}_{n'\underline{m}} (W^{-1}) \cdot
t^{l- \frac 12}_{m- \frac 12 \, \underline{n'+ \frac 12}}(Z)  \\
(l+n') t^{l}_{n'\underline{m}} (W^{-1}) \cdot
t^{l- \frac 12}_{m+ \frac 12 \, \underline{n'- \frac 12}}(Z) &
(l+n') t^{l}_{n'\underline{m}} (W^{-1}) \cdot
t^{l- \frac 12}_{m- \frac 12 \, \underline{n'- \frac 12}}(Z)
\end{pmatrix}.
\end{multline*}
Replacing $n'$ with $n+1$ in the second row we get:
\begin{multline*}
\sum_{\genfrac{}{}{0pt}{}{l,m,n}{m,n \ge -l \ge 1/2}}
\begin{pmatrix}
(l-n) t^l_{n\, \underline{m}} (W^{-1} \bigr) \cdot
t^{l- \frac 12}_{m+\frac12 \, \underline{n+\frac12}}(Z) &
(l-n) t^l_{n\, \underline{m}} (W^{-1}) \cdot
t^{l- \frac 12}_{m-\frac12 \, \underline{n+\frac12}}(Z)  \\
(l+n+1) t^l_{n+1 \, \underline{m}} (W^{-1}) \cdot
t^{l- \frac 12}_{m+\frac12 \, \underline{n+\frac12}}(Z) &
(l+n+1) t^l_{n+1 \, \underline{m}} (W^{-1}) \cdot
t^{l- \frac 12}_{m-\frac12 \, \underline{n+\frac12}}(Z)
\end{pmatrix}  \\
=
\sum_{\genfrac{}{}{0pt}{}{l,m,n}{m,n \ge -l \ge 1/2}}
\begin{pmatrix}
(l-n) t^l_{n\, \underline{m}} (W^{-1})  \\
(l+n+1) t^l_{n+1 \, \underline{m}} (W^{-1})
\end{pmatrix}
\Bigl( t^{l- \frac 12}_{m+\frac12 \, \underline{n+\frac12}}(Z),
t^{l- \frac 12}_{m-\frac12 \, \underline{n+\frac12}}(Z) \Bigr).
\end{multline*}

Similarly, we can expand the second sum in (\ref{exp-sum1+2}) using the
derivative formulas from Lemmas \ref{inverse-deriv}, \ref{deriv_calc}
and the multiplication identities from Lemma \ref{mult-ids}.
\qed

From Lemma \ref{1/N-bound} and Proposition \ref{t-orthogonal2} we obtain:

\begin{thm}  \label{Hilbert-transform}
For each left-regular function $f \in \BB S(\HR^+)$ and $R>0$, the integral
$$
-\frac1{2\pi^2}
\int_{X \in H_R} \frac {(X-W)^{-1}}{N(X-W)} \cdot Dx \cdot f(X),
\qquad W \in R \cdot (\Gamma^- \cup \Gamma^+),
$$
converges.

Let $\P_0^-$ and $\P_{\infty}^-$ be the
projections onto the holomorphic discrete components quasi-regular at the
origin and infinity respectively. Similarly, let $\P_0^+$ and
$\P_{\infty}^+$ be the projections onto the antiholomorphic
discrete components quasi-regular at the origin and infinity respectively.
This integral gives
\begin{align*}
&\bigl( \P_0^-(f) - \P_{\infty}^-(f) \bigr)(W)
\qquad \text{if $W \in R \cdot \Gamma^-$}\\
\text{and} \qquad
&\bigl( \P_0^+(f) - \P_{\infty}^+(f) \bigr)(W)
\qquad \text{if $W \in R \cdot \Gamma^+$.} 
\end{align*}
Similar statement holds for right-regular functions $g \in \BB S'(\HR^+)$.
\end{thm}

\subsection{The Discrete Series Projector}

Recall that in Theorems \ref{discrete_ser_proj} and \ref{Hilbert-transform}
integration takes place over $X \in H_R$ and the variable $W$ lies in
$R \cdot \Gamma^-$ or $R \cdot \Gamma^+$.
But $\HR \cap R \cdot (\Gamma^- \cup \Gamma^+) = \varnothing$.
In this subsection we deform the contour of integration so that the resulting
integrals still provide projections onto the holomorphic and antiholomorphic
discrete series components, but $W$ will lie in open regions containing $H_R$.
In particular, these regions will have non-empty intersection with $\HR$.

Observe that, for any $\sigma>0$, the sets
$$
\begin{pmatrix} \sigma & 0 \\ 0 & \sigma^{-1} \end{pmatrix} \cdot \HR^+
\qquad \text{and} \qquad
\HR^+ \cdot \begin{pmatrix} \sigma & 0 \\ 0 & \sigma^{-1} \end{pmatrix}
\qquad \text{lie inside $\HC^+$}.
$$
We introduce cycles in $\HC^+$
$$
C_{R,\sigma} =
\begin{pmatrix} \sigma & 0 \\ 0 & \sigma^{-1} \end{pmatrix} \cdot H_R
\qquad \text{and} \qquad
C'_{R,\sigma} =
H_R \cdot \begin{pmatrix} \sigma & 0 \\ 0 & \sigma^{-1} \end{pmatrix}.
$$
If $\sigma>1$ these cycles lie in $R \cdot \Gamma^-$, and if $0<\sigma<1$
these cycles lie in $R \cdot \Gamma^+$. In other words, the sets
$$
R \cdot
\begin{pmatrix} \sigma & 0 \\ 0 & \sigma^{-1} \end{pmatrix} \cdot \Gamma^-
\quad \text{and} \quad
R \cdot \Gamma^- \cdot
\begin{pmatrix} \sigma & 0 \\ 0 & \sigma^{-1} \end{pmatrix}
\quad \text{contain $H_R$ when $0<\sigma<1$},
$$
$$
R \cdot \begin{pmatrix}
\sigma & 0 \\ 0 & \sigma^{-1} \end{pmatrix} \cdot \Gamma^+
\quad \text{and} \quad
R \cdot \Gamma^+ \cdot
\begin{pmatrix} \sigma & 0 \\ 0 & \sigma^{-1} \end{pmatrix}
\quad \text{contain $H_R$ when $\sigma>1$}.
$$

\begin{thm}  \label{discr_proj-harm}
For each $R,\sigma>0$, we have well defined continuous linear operators
$\S^{\mp}_{R,\sigma}$ and $\S'^{\mp}_{R,\sigma}: {\cal H}(\HR^+) \to {\cal H}(\HR^+)$
\begin{align*}
(\S^{\mp}_{R,\sigma} \phi)(W) &= -\frac1{2\pi^2} \int_{Z \in C_{R,\sigma}}
\frac{Z^{-1} \cdot Dz}{N(Z-W)} \cdot (\deg \phi)(Z),
\qquad W \in R \cdot
\begin{pmatrix} \sigma & 0 \\ 0 & \sigma^{-1} \end{pmatrix} \cdot \Gamma^{\mp}, \\
(\S'^{\mp}_{R,\sigma} \phi)(W) &= -\frac1{2\pi^2} \int_{Z \in C'_{R,\sigma}}
\frac{Dz \cdot Z^{-1}}{N(Z-W)} \cdot (\deg \phi)(Z),
\qquad W \in R \cdot \Gamma^{\mp} \cdot
\begin{pmatrix} \sigma & 0 \\ 0 & \sigma^{-1} \end{pmatrix}.
\end{align*}
Operators $\S^-_{R,\sigma}$ and $\S'^-_{R,\sigma}$
annihilate the continuous series, the antiholomorphic discrete series
${\cal D}_{discr}^+$ and send
$$
\begin{matrix}
t^l_{n\,\underline{m}}(Z) \quad \mapsto \quad - t^l_{n\,\underline{m}}(W)
- R^{2(2l+1)} \cdot N(W)^{-2l-1} \cdot t^l_{n\,\underline{m}}(W),\\
\quad \\ N(Z)^{-2l-1} \cdot t^l_{n\,\underline{m}}(Z) \quad \mapsto \quad
R^{-2(2l+1)} \cdot t^l_{n\,\underline{m}}(W) + N(W)^{-2l-1} \cdot t^l_{n\,\underline{m}}(W),
\end{matrix} \quad
\begin{matrix}
l = -1, -\frac 32, -2, \dots, \\ m,n \in \BB Z+l,\\ m, n \ge -l.
\end{matrix}
$$
Operators $\S^+_{R,\sigma}$ and $\S'^+_{R,\sigma}$
annihilate the continuous series, the holomorphic discrete series
${\cal D}_{discr}^-$ and send
$$
\begin{matrix}
t^l_{n\,\underline{m}}(Z) \quad \mapsto \quad - t^l_{n\,\underline{m}}(W)
- R^{2(2l+1)} \cdot N(W)^{-2l-1} \cdot t^l_{n\,\underline{m}}(W),\\
\quad \\ N(Z)^{-2l-1} \cdot t^l_{n\,\underline{m}}(Z) \quad \mapsto \quad
R^{-2(2l+1)} \cdot t^l_{n\,\underline{m}}(W) + N(W)^{-2l-1} \cdot t^l_{n\,\underline{m}}(W),
\end{matrix} \quad
\begin{matrix}
l = -1, -\frac 32, -2, \dots, \\ m,n \in \BB Z+l,\\ m, n \le l.
\end{matrix}
$$
\end{thm}

\pf
From Proposition 11 of \cite{FL} we see that the differential form
$Z^{-1} \cdot Dz$ is invariant under the map
$Z \mapsto \begin{pmatrix} \sigma & 0 \\ 0 & \sigma^{-1} \end{pmatrix} Z$.
Similarly, $Dz \cdot Z^{-1}$ is invariant under the map
$Z \mapsto Z \begin{pmatrix} \sigma & 0 \\ 0 & \sigma^{-1} \end{pmatrix}$.
Let
$\tilde Z = \begin{pmatrix} \sigma^{-1} & 0 \\ 0 & \sigma \end{pmatrix} Z$,
$\tilde W = \begin{pmatrix} \sigma^{-1} & 0 \\ 0 & \sigma \end{pmatrix} W$,
and define $\tilde\phi$ by $\tilde \phi(\tilde Z) = \phi(Z)$.
Using (\ref{int_t}) and (\ref{t(Z^+)}) one can show
$$
t^l_{n\,\underline{m}}(\tilde Z) = \sigma^{-2n} \cdot t^{l}_{n\,\underline{m}}(Z),
$$
which implies that the map $\phi \mapsto \tilde\phi$ preserves the
holomorphic discrete series, antiholomorphic discrete series and continuous
components of ${\cal H}(\HR^+)$. We have:
$$
(\S^{\mp}_{R,\sigma} \phi)(W)
= \frac{-1}{2\pi^2} \int_{Z \in C_{R,\sigma}}
\frac{Z^{-1} \cdot Dz  \cdot (\deg \phi)(Z)}{N(Z-W)}
= \frac{-1}{2\pi^2} \int_{\tilde Z \in H_R} \frac{(\deg \tilde\phi)(\tilde Z)}
{N(\tilde Z-\tilde W)} \,\frac {dS}{\|\tilde Z\|}
=
(\S_R^{\mp} \tilde\phi)(\tilde W),
$$
and similarly for $\S'^{\mp}_{R,\sigma}$. Then the result follows from
Theorem \ref{discrete_ser_proj}.
\qed

Now we redo this for left- and right-regular functions. The same argument
that was used to prove Theorem \ref{discr_proj-harm} also proves:

\begin{thm}
For each left-regular function $f \in \BB S(\HR^+)$ we have:
$$
-\frac1{2\pi^2} \int_{Z \in C_{R, \sigma}}
\frac {(Z-W)^{-1}}{N(Z-W)} \cdot Dz \cdot f(Z)
= \begin{cases}
\bigl( \P_0^-(f) - \P_{\infty}^-(f) \bigr)(W) & \text{if 
$W \in R \cdot
\begin{pmatrix} \sigma & 0 \\ 0 & \sigma^{-1} \end{pmatrix} \cdot \Gamma^-$}; \\
\bigl( \P_0^+(f) - \P_{\infty}^+(f) \bigr)(W) & \text{if 
$W \in R \cdot \begin{pmatrix}
\sigma & 0 \\ 0 & \sigma^{-1} \end{pmatrix} \cdot \Gamma^+$.}
\end{cases}
$$
Similarly, for each right-regular function $g \in \BB S'(\HR^+)$ we have:
$$
-\frac1{2\pi^2} \int_{Z \in C'_{R, \sigma}}
g(Z) \cdot Dz \cdot \frac {(Z-W)^{-1}}{N(Z-W)}
= \begin{cases}
\bigl( \P_0^-(g) - \P_{\infty}^-(g) \bigr)(W) & \text{if
$W \in R \cdot \Gamma^- \cdot
\begin{pmatrix} \sigma & 0 \\ 0 & \sigma^{-1} \end{pmatrix}$;} \\
\bigl( \P_0^+(g) - \P_{\infty}^+(g) \bigr)(W) & \text{if
$W \in R \cdot \Gamma^+ \cdot
\begin{pmatrix} \sigma & 0 \\ 0 & \sigma^{-1} \end{pmatrix}$.}
\end{cases}
$$
\end{thm}

\subsection{Second Order Pole in $\HR$}

We consider the linear span of the functions
$t^l_{n\,\underline{m}}(Z) \cdot N(Z)^k$ over
$$
l = -1, -3/2, -2, \dots, \qquad
m,n \in \BB Z+l , \qquad m, n \ge -l, \qquad 0 \le k \le -2l-2
$$
and over
$$
l = -1, -3/2, -2, \dots, \qquad
m,n \in \BB Z+l, \qquad m, n \le l, \qquad 0 \le k \le -2l-2.
$$
We denote the first span by ${\cal D}^{--}$ and the second one by ${\cal D}^{++}$.
It is easy to check that
$$
{\cal D}^{--} \subsetneq
z_{11}^{-2} \cdot \BB C [ z_{11}^{-1}, z_{12}, z_{21}, z_{22} ]^{\le 0}
\quad \subset \quad \BB C [ z_{11}^{-1}, z_{12}, z_{21}, z_{22}]
$$
and
$$
{\cal D}^{++} \subsetneq
z_{22}^{-2} \cdot \BB C [ z_{11}, z_{12}, z_{21}, z_{22}^{-1} ]^{\le 0}
\quad \subset \quad \BB C [ z_{11}, z_{12}, z_{21}, z_{22}^{-1}].
$$
We have the following analogue of Proposition 19 in \cite{FL}:

\begin{prop}
Let ${\cal D}^{--}_{\le 0}$ (respectively ${\cal D}^{--}_{\ge 0}$) denote
the span of $t^l_{n\,\underline{m}}(Z) \cdot N(Z)^k$ with
$l=-1,-\frac32,\dots$, $m,n \in \BB Z+l$, $m, n \ge -l$, $0 \le k <-l$
(respectively $-l \le k \le -2l-2$).
Then ${\cal D}^{--} = {\cal D}^{--}_{\le 0} + {\cal D}^{--}_{\ge 0}$, where
${\cal D}^{--}_{\ge 0}$ is the image of ${\cal D}^{--}_{\le 0}$ under the
inversion map $F(Z) \mapsto N(Z)^{-2} \cdot F\bigl(Z/N(Z)\bigr)$ and
$$
{\cal D}^{--}_{\le 0} =
z_{11}^{-2} \cdot \BB C [z_{11}^{-1}, z_{12}, z_{21}, z_{22}]^{\le 0}.
$$

Similar statement holds for ${\cal D}^{++}$.
\end{prop}



Recall an action $\rho_1$ of $GL(2,\HC)$ on the space of functions on $\HC$
with singularities defined in \cite{FL} by
\begin{multline*}
\rho_1(h): \: F(Z) \mapsto \bigl( \rho_1(h)F \bigr)(Z) =
\frac {F \bigl( (aZ+b)(cZ+d)^{-1} \bigr)}{N(cZ+d) \cdot N(a'-Zc')},  \\
h = \begin{pmatrix} a' & b' \\ c' & d' \end{pmatrix},\:
h^{-1} = \begin{pmatrix} a & b \\ c & d \end{pmatrix} \in GL(2,\HC).
\end{multline*}
Differentiating, we obtain an action of
$\mathfrak{gl}(2,\HC) \simeq \mathfrak{gl}(4,\BB C)$ also denoted by $\rho_1$;
this action was described in Lemma 68 in \cite{FL}.
Recall the representations $\Zh^+$ and $\Zh^-$ of
$\mathfrak{sl}(2,\HC)$ introduced in \cite{FL}.

\begin{thm}
We have ${\cal D}^{--} \simeq \Zh^-$ and ${\cal D}^{++} \simeq \Zh^+$
as representations of $\mathfrak{gl}(4,\BB C)$.
Moreover, ${\cal D}^{--}$ and ${\cal D}^{++}$ are irreducible representations of
$\mathfrak{sl}(4,\BB C)$ with highest (or lowest) weight vectors
$t^{-1}_{1\,\underline{1}}(Z)=z_{11}^{-2}$ and
$t^{-1}_{-1\,\underline{-1}}(Z)=z_{22}^{-2}$ respectively.
\end{thm}

\pf
The proof proceeds in the same way as that of Theorem \ref{D-component}.
First, we check using Lemma 68 in \cite{FL} and Lemmas \ref{deriv_calc},
\ref{mult-ids} that the spaces ${\cal D}^{--}$ and ${\cal D}^{++}$ are
invariant under the $\rho_1$ action and that
$t^{-1}_{1\,\underline{1}}(Z)$ and $t^{-1}_{-1\,\underline{-1}}(Z)$
generate ${\cal D}^{--}$ and ${\cal D}^{++}$ respectively.
Recall the element $\gamma_0$ from (\ref{Cayley-composition}).
Then
\begin{align*}
\rho_l(\gamma_0) : \qquad \Zh^+ \ni 1 \quad &\mapsto \quad
z_{22}^{-2}= t^{-1}_{-1\,\underline{-1}}(Z) \in {\cal D}^{++}, \\
\Zh^- \ni \frac1{N(Z)} \quad &\mapsto \quad
z_{11}^{-2}= t^{-1}_{1\,\underline{1}}(Z) \in {\cal D}^{--}.
\end{align*}
This proves ${\cal D}^{\mp\mp} \simeq \Zh^{\mp}$, irreducibility and the
statement about the highest (or lowest) weight vectors.
\qed

Next we establish two expansions for $\frac 1{N(Z-W)^2}$.

\begin{prop}
We have the following matrix coefficient expansions
$$
\frac1{N(Z-W)^2} =
\sum_{\text{\tiny $\begin{matrix} k,l,m,n \\ m,n \ge -l \ge 1 \\ 0 \le k \le -2l-2 \end{matrix}$}}
-(2l+1) t^l_{n\underline{m}}(W) \cdot N(W)^k \cdot
t^l_{m\underline{n}}(Z^{-1}) \cdot N(Z)^{-k-2}
$$
which converges pointwise absolutely whenever $WZ^{-1} \in \Gamma^-$ and
$$
\frac1{N(Z-W)^2} =
\sum_{\text{\tiny $\begin{matrix} k,l,m,n \\ m,n \le l \le -1 \\ 0 \le k \le -2l-2 \end{matrix}$}}
-(2l+1) t^l_{n\underline{m}}(W) \cdot N(W)^k \cdot
t^l_{m\underline{n}}(Z^{-1}) \cdot N(Z)^{-k-2}
$$
which converges pointwise absolutely whenever $WZ^{-1} \in \Gamma^+$.
\end{prop}

\begin{rem}
Note that, unlike the expansions of $\frac 1{N(Z-W)}$ given in
Proposition \ref{expansion2}, the matrix coefficients of the limits of
the discrete series do not enter the expansions of $\frac 1{N(Z-W)^2}$.
\end{rem}

\pf
The proof is similar to that of Proposition \ref{expansion2}.
Thus we can assume $WZ^{-1} \in \Gamma^- \cap K_{\BB C}$, and so
$WZ^{-1} = \begin{pmatrix} \lambda_1 & 0 \\ 0 & \lambda_2 \end{pmatrix}$
for some $\lambda_1, \lambda_2 \in \BB C$ with
$|\lambda_1| > 1 > |\lambda_2| > 0$.
Using (\ref{Theta}) and (\ref{t-mult}) and letting indices $k$, $l$, $m$, $n$
run over
$$
l \le -1, \qquad m,n \ge -l, \qquad 0 \le k \le -2l-2,
$$
we obtain:
\begin{multline*}
\sum_{k,l,m,n} (2l+1) t^l_{n\underline{m}}(W) \cdot N(W)^k \cdot
t^l_{m\underline{n}}(Z^{-1}) \cdot N(Z)^{-k-2}  \\
= \frac 1{N(Z)^2}
\sum_{k,l,n} (2l+1) t^l_{n\underline{n}}(WZ^{-1}) \cdot N(WZ^{-1})^k
= \frac 1{N(Z)^2}
\sum_{k,l} (2l+1) \tilde\Theta_l^-(WZ^{-1}) \cdot N(WZ^{-1})^k  \\
= \frac 1{N(Z)^2}
\sum_{k,l} (2l+1) \frac{\lambda_1^{2l+1}}{\lambda_1-\lambda_2}
\cdot (\lambda_1 \lambda_2)^k
= \frac 1{N(Z)^2} \sum_{l \le -1} (2l+1) \frac{\lambda_1^{2l+1} - \lambda_2^{-2l-1}}
{(\lambda_1-\lambda_2)(1-\lambda_1 \lambda_2)}  \\
= - \frac 1{N(Z)^2} \frac 1{(1-\lambda_1)^2 (1-\lambda_2)^2}
= - \frac 1{N(Z)^2 \cdot N(1-WZ^{-1})^2}
= - \frac 1{N(Z-W)^2}.
\end{multline*}


The other matrix coefficient expansion is proved in the same way.
\qed

Let $\tilde e_3=ie_3= \begin{pmatrix} 1 & 0 \\ 0 & -1 \end{pmatrix} \in \HC$.
We realize the group $U(2,2)$ as the subgroup of elements of $GL(2,\HC)$
preserving the Hermitian form on $\BB C^4$ given by the $4 \times 4$ matrix
$\begin{pmatrix} \tilde e_3 & 0 \\ 0 & -\tilde e_3 \end{pmatrix}$. Explicitly,
\begin{align*}
U(2,2) &= \Biggl\{ \begin{pmatrix} a & b \\ c & d \end{pmatrix};\:
a,b,c,d \in \HC,\:
\begin{matrix} a^* \tilde e_3 a = \tilde e_3 + c^* \tilde e_3 c \\
d^* \tilde e_3 d = \tilde e_3 + b^* \tilde e_3 b \\
a^* \tilde e_3 b = c^* \tilde e_3 d \end{matrix} \Biggr\}  \\
&= \Biggl\{ \begin{pmatrix} a & b \\ c & d \end{pmatrix};\:
a,b,c,d \in \HC,\:
\begin{matrix} a^* \tilde e_3 a = \tilde e_3 + b^* \tilde e_3 b \\
d^* \tilde e_3 d = \tilde e_3 + c^* \tilde e_3 c \\
a \tilde e_3 c^* = b \tilde e_3 d^* \end{matrix} \Biggr\}.
\end{align*}
The Lie algebra of $U(2,2)$ is
$$
\mathfrak{u}(2,2) = \biggl\{
\begin{pmatrix} A & B \\ \tilde e_3B^*\tilde e_3 & D \end{pmatrix} ;\:
A,B,D \in \HC ,\: \tilde e_3A=-(\tilde e_3A)^*, \tilde e_3D=-(\tilde e_3D)^*
\biggr\}.
$$
If $\begin{pmatrix} a & b \\ c & d \end{pmatrix} \in U(2,2)$,
then $\begin{pmatrix} a & b \\ c & d \end{pmatrix}^{-1} =
\tilde e_3 \begin{pmatrix} a^* & -c^* \\ -b^* & d^* \end{pmatrix} \tilde e_3$.
The group $U(2,2)$ acts on $\HC$ by conformal transformations
preserving $U(1,1)$, where we identify $U(1,1)$ with a subset of $\HC$:
\begin{equation}  \label{u(1,1)}
U(1,1) = \{ Z \in \HC ;\: Z^* \tilde e_3 Z = \tilde e_3 \}.
\end{equation}
We orient $U(1,1)$ so that $\{\tilde e_0, \tilde e_1, \tilde e_2, e_3\}$
is a positive basis of the tangent space at $1 \in U(1,1)$.

We have a symmetric bilinear pairing on ${\cal D}^{--} \oplus {\cal D}^{++}$:
$$
\langle F_1,F_2 \rangle_1 = \frac i{2\pi^3} \int_{U(1,1)} F_1(Z) \cdot F_2(Z) \,dV
$$
(recall that $dV= dz^0 \wedge dz^1 \wedge dz^2 \wedge dz^3$
is a holomorphic 4-form on $\HC$).
This pairing is $\mathfrak{gl}(2,\HC)$-invariant, which follows immediately
from Lemma 61 in \cite{FL} and
$dZ^4 = dz_{11} \wedge dz_{12} \wedge dz_{21} \wedge dz_{22} =4dV$.
Related to this bilinear pairing we have a $\mathfrak{u}(2,2)$-invariant
inner product that is nondegenerate on ${\cal D}^{--}$ and ${\cal D}^{++}$:
$$
(F_1,F_2)_1 = \frac i{2\pi^3}
\int_{U(1,1)} F_1(Z) \cdot \B{F_2(Z)} \,\frac{dV}{N(Z)^2}.
$$
(Checking the $\mathfrak{u}(2,2)$-invariance requires some verification.)

Every element in $U(1,1)$ can be uniquely written as
$e^{i\theta} \cdot Z$ with $Z \in SU(1,1)$ and $\theta \in [0,\pi)$.
Thus we can identify $U(1,1)$ with $[0,\pi) \times SU(1,1)$ and
$dV/N(Z)^2$ restricted to $U(1,1)$ becomes $i d\theta \wedge Dz \cdot Z^{-1}$
(since both sides express a $U(1,1)$-invariant volume form and coincide at $1$).
From Proposition \ref{t-orthogonal} we immediately obtain:

\begin{prop}
We have the following orthogonality relations:
$$
\bigl\langle t^{l'}_{n'\,\underline{m'}}(Z) \cdot N(Z)^{k'},
t^l_{m\underline{n}}(Z^{-1}) \cdot N(Z)^{-k-2} \bigr\rangle_1
= - \frac1{2l+1} \delta_{kk'}\delta_{ll'} \delta_{mm'} \delta_{nn'},
$$
where the indices $k,l,m,n$ are
$l = -1, -\frac 32, -2, \dots$, $m,n \in \BB Z +l$, $m,n \ge -l$ or $m,n \le l$,
$k \in \BB Z$ and similarly for $k',l',m',n'$.
\end{prop}

From these orthogonality relations and expansions for $\frac 1{N(Z-W)^2}$
we obtain a formula similar to Proposition 73 in \cite{FL}:

\begin{thm}
Let $\P^{--}$ and $\P^{++}$ denote the projections of
${\cal D}^{--} \oplus {\cal D}^{++}$ onto ${\cal D}^{--}$ and ${\cal D}^{++}$
respectively.
For each function $F \in {\cal D}^{--} \oplus {\cal D}^{++}$ and $R>0$,
$$
\frac{i}{2\pi^3} \int_{Z \in R \cdot U(1,1)} \frac{F(Z)}{N(Z-W)^2} \,dV =
\begin{cases}
(\P^{--} F)(W) & \text{if $W \in R \cdot \Gamma^-$;}\\
(\P^{++} F)(W) & \text{if $W \in R \cdot \Gamma^+$.}
\end{cases}
$$
In particular, the integral converges absolutely for
$W \in R \cdot (\Gamma^- \cup \Gamma^+)$.
\end{thm}

\section{Separation of the Series for $SL(2,\BB C)/SU(1,1)$}

\subsection{The Cayley Transform between $\HR$ and $\BB M$}

Recall that
$\tilde e_3=ie_3= \begin{pmatrix} 1 & 0 \\ 0 & -1 \end{pmatrix} \in \HC$
and that $U(1,1)$ is identified with a subset of $\HC$ via (\ref{u(1,1)}).

\begin{prop}  \label{Cayley}
Consider an element
$$\gamma = \frac 1{\sqrt{2}}
\begin{pmatrix} \tilde e_3 & i \tilde e_3 \\ 1 & -i \end{pmatrix}
\in GL(2,\HC) \qquad \text{with} \qquad
\gamma^{-1} = \frac 1{\sqrt{2}}
\begin{pmatrix} \tilde e_3 & 1 \\ -i \tilde e_3 & i \end{pmatrix}.
$$
The fractional linear map on $\HC$
$$
\pi_l(\gamma): \: Z \mapsto i(\tilde e_3 Z +1)(\tilde e_3 Z -1)^{-1}
$$
maps $\BB M \to U(1,1) \subset \HC$ (with singularities) and sends
the unit one-sheeted hyperboloid $\tilde H = \{ Y \in \BB M ;\: N(Y) = 1 \}$
into $SU(1,1) = \{ X \in \HR ;\: N(X)=1 \}$.
The singularities of $\pi_l(\gamma)$ on $\BB M$ lie along the
one-sheeted hyperboloid
$\{Y \in \BB M ;\: N(Y)=1,\: \re(e_3Y)=0 \}$.

Conversely, the fractional linear map on $\HC$
$$
\pi_l(\gamma^{-1}): \: Z \mapsto \tilde e_3 (Z+i)(Z-i)^{-1}
$$
maps $U(1,1) \to \BB M$ (with singularities) and sends $SU(1,1)$ into the
hyperboloid $\tilde H$ The singularities of $\pi_l(\gamma^{-1})$ on $\HR$
lie along the two-sheeted hyperboloid $\{ X \in SU(1,1) ;\: \re X =0 \}$.
\end{prop}

The map $\pi_l(\gamma)$ and its inverse were studied in \cite{KouO}, we
think of these maps as quaternionic analogues of Cayley transform.
For future reference we spell out that if
$Y = \begin{pmatrix} y_{11} & y_{12} \\ y_{21} & y_{22} \end{pmatrix} \in \BB M$
and $X = i(\tilde e_3 Y+1)(\tilde e_3 Y -1)^{-1}$, then
$N(\tilde e_3Y -1) = 1-N(Y)+y_{22}-y_{11}$ and
\begin{multline}  \label{Y-tilde}
X = \frac{-i}{N(\tilde e_3Y -1)}
\begin{pmatrix} 1+N(Y)+y_{11}+y_{22} & 2y_{12} \\ -2y_{21} & 1+N(Y)-y_{11}-y_{22}
\end{pmatrix}  \\
=
i - \frac{2i}{N(\tilde e_3Y -1)}
\begin{pmatrix} 1+y_{22} & y_{12} \\ -y_{21} & 1-y_{11} \end{pmatrix}.
\end{multline}
We orient the hyperboloid $\tilde H \subset \BB M$ as the boundary of the open
set $\{Y \in \BB M ;\: N(Y)<1 \}$ and denote by $\sgnx(Y)$ the sign of the
$y^3$-coordinate of $Y= y^0 \tilde e_0 + y^1 e_1 + y^2 e_2 + y^3 e_3 \in \BB M$.

\begin{lem}  \label{orientations}
The restriction of $\pi_l(\gamma)$ to $\tilde H \to SU(1,1)$
preserves the orientations for $\{ Y \in \tilde H ;\: \sgnx(Y) >0 \}$
and reverses the orientations for $\{ Y \in \tilde H ;\: \sgnx(Y) <0 \}$.

Conversely, the restriction of $\pi_l(\gamma^{-1})$ to $SU(1,1) \to \tilde H$
preserves the orientations for $\{ X \in SU(1,1) ;\: \re X <0 \}$
and reverses the orientations for $\{ X \in SU(1,1) ;\: \re X >0 \}$.
\end{lem}


\begin{thm}  \label{inner_product_pullback}
For any $\phi_1(X), \phi_2(X) \in {\cal H}(\HR^+)$ that extend holomorphically
to an open neighborhood of $SU(1,1)$ in $\HC$ we have:
\begin{multline*}
\langle \phi_1 , \phi_2 \rangle = - \frac 1{2\pi^2}
\int_{X \in SU(1,1)} (\deg_X \phi_1)(X) \cdot \phi_2(X) \,\frac{dS}{\|X\|}  \\
=
\frac i{2\pi^2} \int_{Y \in \tilde H} \sgnx(Y) \cdot
\bigl(\deg_Y (\pi_l^0(\gamma)\phi_1) \bigr)(Y) \cdot
(\pi_l^0(\gamma)\phi_2)(Y) \,\frac{dS}{\|Y\|},
\end{multline*}
where
$$
\bigl( \pi_l^0(\gamma)\phi_j \bigr) (Y) =
\frac {-2}{N(\tilde e_3 Y -1)} \cdot
\phi_j \bigl( i(\tilde e_3 Y +1)(\tilde e_3 Y -1)^{-1} \bigr),
\qquad j=1,2.
$$
\end{thm}

\pf
The main ingredient of the proof is the following lemma which is obtained
by direct computation.

\begin{lem}  \label{degrees}
For $Y \in \BB \HC$ and
$X = i(\tilde e_3 Y +1)(\tilde e_3 Y -1)^{-1} \in SU(1,1)$ we have:
\begin{multline*}
\frac {N(Y)+1}{N(\tilde e_3 Y -1)^2} \cdot (\deg_X \phi)(X)
+i\frac{1-N(Y)}{N(\tilde e_3 Y -1)^2} (\partial_{11}\phi+\partial_{22}\phi)(X) \\
= \deg_Y \biggl( \frac 1{N(\tilde e_3 Y -1)} \cdot
\phi \bigl( i(\tilde e_3 Y +1)(\tilde e_3 Y -1)^{-1} \bigr) \biggr).
\end{multline*}
\end{lem}

Writing $X = i(\tilde e_3 Y+1)(\tilde e_3 Y -1)^{-1}$,
$dS/\|Y\| = i Dy/Y$ on $\tilde H$ and using Lemma 10 together with
Proposition 11 from \cite{FL}, we can rewrite
\begin{multline*}
- \int_{X \in SU(1,1)} (\deg_X \phi_1)(X) \cdot \phi_2(X) \,\frac{dS}{\|X\|}
=
- \int_{X \in SU(1,1)}
(\deg_X \phi_1)(X) \cdot \frac {Dx}X \cdot \phi_2(X)  \\
=
8 \int_{Y \in \tilde H} \sgnx(Y) \cdot (\deg_X \phi_1)(X) \cdot
\frac {(\tilde e_3 -Y)^{-1}}{N(\tilde e_3 - Y)} \cdot Dy \cdot
\frac {(\tilde e_3 Y -1)^{-1}}{N(\tilde e_3 Y -1)} \cdot
\frac{\tilde e_3 Y -1}{\tilde e_3 Y +1} \cdot \phi_2(X)  \\
=
8i \int_{Y \in \tilde H} \sgnx(Y) \cdot (\deg_X \phi_1)(X) \cdot
\frac {\phi_2(X)}{N(\tilde e_3 Y -1)^3} \,\frac {dS}{\|Y\|}  \\
=
4i \int_{Y \in \tilde H} \sgnx(Y) \cdot \deg_Y \biggl( \frac 1{N(\tilde e_3 Y -1)}
\cdot \phi \bigl( i(\tilde e_3 Y +1)(\tilde e_3 Y -1)^{-1} \bigr) \biggr)
\cdot \frac {\phi_2(X)}{N(\tilde e_3 Y -1)} \,\frac {dS}{\|Y\|}  \\
=
i\int_{Y \in \tilde H} \sgnx(Y) \cdot \bigl(\deg_Y (\pi_l^0(\gamma)\phi_1) \bigr)(Y)
\cdot (\pi_l^0(\gamma)\phi_2)(Y) \,\frac{dS}{\|Y\|}.
\end{multline*}
\qed


\subsection{The Continuous Series Projector on $\BB M$}
\label{M-cont}

Recall that ${\cal D}^-$ and ${\cal D}^+$ are irreducible representations of
$\mathfrak{sl}(4,\BB C)$ with highest (or lowest) weight vectors
$t^{-\frac12}_{\frac12\,\underline{\frac12}}(Z) = \frac1{z_{11}}$ and
$t^{-\frac12}_{-\frac12\,\underline{-\frac12}}(Z) = \frac1{z_{22}}$ respectively
(Theorem \ref{D-component}). From (\ref{Y-tilde}) we can see that
$$
\biggl( \pi_l^0(\gamma) \Bigl( \frac1{z_{11}} \Bigr)\biggr)(Y)
= \frac {-2i}{N(Y+1)} \qquad \text{and} \qquad
\biggl( \pi_l^0(\gamma) \Bigl( \frac1{z_{22}} \Bigr)\biggr)(Y)
= \frac {-2i}{N(Y-1)}.
$$
These functions have singularities along the set
$\{Y \in \BB M ;\: N(Y)=-1,\: \re Y =0 \}$ which has codimension 2.
It is easy to see that they define a tempered distribution on $\BB M$,
and so their Fourier transforms are $L^2$-functions on the light cone
in the space dual to $\BB M$.
The space of functions on $\BB M$ that arise as Fourier transforms of
the $L^2$-functions on the light cone is known as the continuous series
component of $\BB M$.
This proves that the Cayley transform switches the discrete series
component on $\HR$ and the continuous series component on $\BB M$.

To describe the images of $\Gamma^-$ and $\Gamma^+$ under the Cayley
transform $\pi_l(\gamma^{-1})$ we recall the generalized upper and lower
half-planes introduced in Section 3.5 in \cite{FL}:
\begin{align*}
\BB T^- &= \{ Z = W_1 + i W_2 \in \HC ;\: W_1, W_2 \in \BB M ,\:
\text{$iW_2$ is positive definite} \},  \\
\BB T^+ &= \{ Z = W_1 + i W_2 \in \HC ;\: W_1, W_2 \in \BB M ,\:
\text{$iW_2$ is negative definite} \}.
\end{align*}
In the context of our paper,  Lemma 1.1 in \cite{KouO} can be restated as
follows.

\begin{lem}  \label{Gamma-T}
The Cayley transform $\pi_l(\gamma^{-1})$ sends $\Gamma^-$ and $\Gamma^+$
biholomorphically into respectively $\BB T^-_0$ and $\BB T^+_0$, where
$$
\BB T^-_0 = \{ Z \in \BB T^- ;\: N(Z- \tilde e_3) \ne 0 \}
\qquad \text{and} \qquad
\BB T^+_0 = \{ Z \in \BB T^+ ;\: N(Z- \tilde e_3) \ne 0 \}.
$$
\end{lem}


Recall that $\tilde H$ is the unit hyperboloid of one sheet in $\BB M$, and
let ${\cal H}(\tilde H)$ denote the space of holomorphic functions $\phi$
defined on some connected open neighborhood $U_{\phi}$ of $\tilde H$ in $\HC$
which can be written as $\phi(Z) = \bigl(\pi_l^0(\gamma)\tilde\phi\bigr)(Z)$
for some  $\tilde\phi \in {\cal H}(\HR^+)$ and $Z \in U_{\phi} \cap U(1,1)$.
Clearly, functions in ${\cal H}(\tilde H)$ are harmonic.
We define operators $\widetilde{\S}^{\mp}$ on ${\cal H}(\tilde H)$ by
\begin{align*}
(\widetilde{\S}^-\phi)(Z) &=
\frac i{2\pi^2} \int_{Y \in \tilde H} \sgnx(Y) \cdot
\frac {(\deg \phi)(Y)}{N(Y-Z)}\,\frac {dS}{\|Y\|}, \qquad Z \in \BB T^-_0,\\
(\widetilde{\S}^+\phi)(Z) &=
\frac i{2\pi^2} \int_{Y \in \tilde H} \sgnx(Y) \cdot
\frac {(\deg \phi)(Y)}{N(Y-Z)}\,\frac {dS}{\|Y\|}, \qquad Z \in \BB T^+_0.
\end{align*}

\begin{thm}
Let $\phi \in {\cal H}(\tilde H)$, then $\widetilde{\S}^{\mp}\phi$
are well defined functions on $\BB T^-_0$ and $\BB T^+_0$ respectively.
The operator $\widetilde{\S}^-$ annihilates the discrete series on $\BB M$,
the image of the antiholomorphic discrete series ${\cal D}_{discr}^+$ and sends
\begin{align*}
\bigl(\pi_l^0(\gamma)t^l_{n\,\underline{m}}\bigr)(Y) \quad &\mapsto \quad
- \bigl(\pi_l^0(\gamma) \bigl(t^l_{n\,\underline{m}}
+ N(Z)^{-2l-1} \cdot t^l_{n\,\underline{m}}\bigr)\bigr)(Z), \\
\bigl(\pi_l^0(\gamma) \bigl(N(Y)^{-2l-1} \cdot t^l_{n\,\underline{m}}\bigr)\bigr)(Y)
\quad &\mapsto \quad \bigl(\pi_l^0(\gamma) \bigl(t^l_{n\,\underline{m}}
+ N(Z)^{-2l-1} \cdot t^l_{n\,\underline{m}}\bigr)\bigr)(Z),
\end{align*}
$$
l = -1, -3/2, -2, \dots, \qquad m,n \in \BB Z+l, \qquad m, n \ge -l.
$$
The operator $\widetilde{\S}^+$ annihilates the discrete series on $\BB M$,
the image of the holomorphic discrete series ${\cal D}_{discr}^-$ and sends
\begin{align*}
\bigl(\pi_l^0(\gamma)t^l_{n\,\underline{m}}\bigr)(Y) \quad &\mapsto \quad
- \bigl(\pi_l^0(\gamma) \bigl(t^l_{n\,\underline{m}}
+ N(Z)^{-2l-1} \cdot t^l_{n\,\underline{m}}\bigr)\bigr)(Z), \\
\bigl(\pi_l^0(\gamma) \bigl(N(Y)^{-2l-1} \cdot t^l_{n\,\underline{m}}\bigr)\bigr)(Y)
\quad &\mapsto \quad \bigl(\pi_l^0(\gamma) \bigl(t^l_{n\,\underline{m}}
+ N(Z)^{-2l-1} \cdot t^l_{n\,\underline{m}}\bigr)\bigr)(Z),
\end{align*}
$$
l = -1, -3/2, -2, \dots, \qquad m,n \in \BB Z+l, \qquad m, n \le l.
$$
\end{thm}

\pf
Let $X=i(\tilde e_3 Y +1)(\tilde e_3 Y -1)^{-1}$,
$Z'=i(\tilde e_3 Z +1)(\tilde e_3 Z -1)^{-1}$, and write
$\phi(Y) = \bigl( \pi_l^0(\gamma)\tilde\phi\bigr)(Y)$ for some analytic
$\tilde\phi \in {\cal H}(\HR^+)$.
By Lemma 10 from \cite{FL},
$$
\Bigl( \pi_l^0(\gamma)_X \frac1{N(X-Z')} \Bigr) (Y) =
-\frac12 N(\tilde e_3 Z -1) \cdot \frac1{N(Y-Z)}.
$$
By Theorem \ref{inner_product_pullback} we have
\begin{multline*}
(\widetilde{\S}^{\mp}\phi)(Z)
=
\frac i{2\pi^2} \int_{Y \in \tilde H} \sgnx(Y) \cdot
\frac{(\deg \phi)(Y)}{N(Y-Z)}\,\frac {dS}{\|Y\|} \\
=
\frac1{\pi^2 N(\tilde e_3Z-1)} \int_{X \in SU(1,1)}
\frac{(\deg\tilde\phi)(X)}{N(X-Z')}\,\frac {dS}{\|X\|}
=
\frac{-2}{N(\tilde e_3Z-1)} \cdot (\S_1^{\mp} \tilde \phi)(Z').
\end{multline*}
Then the result follows from Theorem \ref{discrete_ser_proj}
and Lemma \ref{Gamma-T}.
\qed

Note that, for $Z \in \BB M$ and $s \in \BB R$, we have
$Z+se_0 \in \BB T^-$ if $s>0$ and $Z+se_0 \in \BB T^+$ if $s<0$.
Moreover, $Z+se_0 \in \BB T^-_0 \sqcup \BB T^+_0$ if $|s|$ is sufficiently small.
Taking limits $s \to 0^{\pm}$ we obtain:

\begin{cor}  \label{M-cont-cor}
Let $Z \in \BB M$ and let $Y$ range over $\tilde H$, write
$Y = t \tilde e_0 + y^1 e_1 + y^2 e_2 + y^3 e_3$,
$Z= \tilde t \tilde e_0 + z^1 e_1 + z^2 e_2 + z^3 e_3$,
then the operator on ${\cal H}(\tilde H)$
$$
\phi(Y) \quad\mapsto\quad \lim_{\epsilon \to 0^+} \frac i{2\pi^2} \int_{Y \in \tilde H}
\sgnx(Y) \cdot \frac {(\deg \phi)(Y)}{N(Y-Z) + i\epsilon\sgn(t-\tilde t)}
\,\frac {dS}{\|Y\|}
$$
annihilates the discrete series on $\BB M$, the image of the antiholomorphic
discrete series ${\cal D}_{discr}^+$ and sends
\begin{align*}
\bigl(\pi_l^0(\gamma)t^l_{n\,\underline{m}}\bigr)(Y) \quad &\mapsto \quad
- \bigl(\pi_l^0(\gamma) \bigl(t^l_{n\,\underline{m}}
+ N(Z)^{-2l-1} \cdot t^l_{n\,\underline{m}}\bigr)\bigr)(Z), \\
\bigl(\pi_l^0(\gamma) \bigl(N(Y)^{-2l-1} \cdot t^l_{n\,\underline{m}}\bigr)\bigr)(Y)
\quad &\mapsto \quad \bigl(\pi_l^0(\gamma) \bigl(t^l_{n\,\underline{m}}
+ N(Z)^{-2l-1} \cdot t^l_{n\,\underline{m}}\bigr)\bigr)(Z),
\end{align*}
$$
l = -1, -3/2, -2, \dots, \qquad m,n \in \BB Z+l, \qquad m, n \ge -l.
$$
Similarly, the operator on ${\cal H}(\tilde H)$
$$
\phi(Y) \quad\mapsto\quad \lim_{\epsilon \to 0^+} \frac i{2\pi^2} \int_{Y \in \tilde H}
\sgnx(Y) \cdot \frac {(\deg \phi)(Y)}{N(Y-Z) - i\epsilon\sgn(t-\tilde t)}
\,\frac {dS}{\|Y\|}
$$
annihilates the discrete series on $\BB M$, the image of the holomorphic
discrete series ${\cal D}_{discr}^-$ and sends
\begin{align*}
\bigl(\pi_l^0(\gamma)t^l_{n\,\underline{m}}\bigr)(Y) \quad &\mapsto \quad
- \bigl(\pi_l^0(\gamma) \bigl(t^l_{n\,\underline{m}}
+ N(Z)^{-2l-1} \cdot t^l_{n\,\underline{m}}\bigr)\bigr)(Z), \\
\bigl(\pi_l^0(\gamma) \bigl(N(Y)^{-2l-1} \cdot t^l_{n\,\underline{m}}\bigr)\bigr)(Y)
\quad &\mapsto \quad \bigl(\pi_l^0(\gamma) \bigl(t^l_{n\,\underline{m}}
+ N(Z)^{-2l-1} \cdot t^l_{n\,\underline{m}}\bigr)\bigr)(Z),
\end{align*}
$$
l = -1, -3/2, -2, \dots, \qquad m,n \in \BB Z+l, \qquad m, n \le l.
$$
\end{cor}




Next we redo this for left-regular functions.
We denote by $\BB S(\tilde H)$ the space of holomorphic left-regular
$\BB S$-valued functions $f$ that are defined on some connected open
neighborhood $U_f$ of $\tilde H$ in $\HC$ and can be written as
$f(Z) = \bigl(\pi_l(\gamma)\tilde f\bigr)(Z)$ for some 
$\tilde f \in \BB S(\HR^+)$ and $Z \in U_f \cap U(1,1)$.
Similarly we can define the space of right-regular functions $\BB S'(\tilde H)$.

\begin{thm}  \label{M-cont-proj-reg}
For each left-regular function $f \in \BB S(\tilde H)$ we have
$$
\frac1{2\pi^2} \int_{Y \in \tilde H} \sgnx(Y) \cdot
\frac {(Y-Z)^{-1}}{N(Y-Z)} \cdot Dy \cdot f(Y)
= \begin{cases}
\bigl( \widetilde{\P}_{\infty}^-(f) - \widetilde{\P}_0^-(f) \bigr)(Z) &
\text{if $Z \in \BB T^-_0$;} \\
\bigl( \widetilde{\P}_{\infty}^+(f) - \widetilde{\P}_0^+(f) \bigr)(Z) &
\text{if $Z \in \BB T^+_0$,}
\end{cases}
$$
where $\widetilde{\P}_0^-$, $\widetilde{\P}_{\infty}^-$ denote the projections
onto the $\pi_l(\gamma)$-images of the holomorphic discrete components of
$\BB S(\HR^+)$ quasi-regular at the origin and infinity respectively, and
$\widetilde{\P}_0^+$, $\widetilde{\P}_{\infty}^+$ are the projections onto the
$\pi_l(\gamma)$-images of the antiholomorphic discrete components of
$\BB S(\HR^+)$ quasi-regular at the origin and infinity respectively.
\end{thm}

\pf
Changing the variables
$Y = \tilde e_3 (X+i)(X-i)^{-1}$, $Z = \tilde e_3 (Z'+i)(Z'-i)^{-1}$,
and using Lemma 10 with Proposition 11 from \cite{FL},
\begin{multline*}
\int_{Y \in \tilde H} \sgnx(Y) \cdot
\frac {(Y-Z)^{-1}}{N(Y-Z)} \cdot Dy \cdot f(Y)  \\
= \frac{(Z'-i)}{N(Z'-i)^{-1}} \cdot \int_{X \in SU(1,1)}
\frac {(X-Z')^{-1}}{N(X-Z')} \cdot Dx \cdot
\frac{(X-i)^{-1}}{N(X-i)} \cdot f\bigl(\tilde e_3 (X+i)(X-i)^{-1}\bigr)  \\
= -2\pi^2 \bigl( \P_0^{\mp}(f) - \P_{\infty}^{\mp}(f) \bigr)
\bigl(\tilde e_3 (Z'+i)(Z'-i)^{-1}\bigr)
= 2\pi^2 \bigl( \widetilde{\P}_{\infty}^{\mp}(f)
- \widetilde{\P}_0^{\mp}(f) \bigr)(Z)
\end{multline*}
by Theorem \ref{Hilbert-transform}.
\qed

\subsection{Basis of Harmonic Functions on $\BB M$}

Writing the Minkowski space $\BB M$ as
$$
\biggl\{ Y= y^0 \tilde e_0 + y^1 e_1 + y^2 e_2 + y^3 e_3 =
\begin{pmatrix} -iy^0-iy^3 & -iy^1-y^2 \\ -iy^1+y^2 & -iy^0+iy^3 \end{pmatrix}
;\: y^0, y^1, y^2, y^3 \in \BB R \biggr\},
$$
we can embed $\BB R^3$ into $\BB M$ so that
$(x^1,x^2,x^3) \leftrightsquigarrow x^1e_1+x^2e_2+x^3e_3$.
Note that, for $Y = t \tilde e_0 + x^1e_1 + x^2e_2 + x^3e_3 \in \BB M$,
$N(Y) = r^2-t^2$, where $r=\sqrt{(x^1)^2 + (x^2)^2 + (x^3)^2}$.
We parameterize the unit hyperboloid $\tilde H$ as
\begin{equation}  \label{hyper-param}
\begin{tabular}{lcc}
$x^1 = \cosh \rho \sin \theta \cos \phi$ & \qquad & $-\infty < \rho < \infty$\\
$x^2 = \cosh \rho \sin \theta \sin \phi$ & \qquad & $0 \le \phi < 2\pi$  \\
$x^3 = \cosh \rho \cos \theta$ & \qquad & $0 \le \theta < \pi$  \\
$t= \sinh \rho$
\end{tabular}
\end{equation}

\begin{lem}
We have the following $SO(3,1)$-invariant measure on the hyperboloid $\tilde H$:
$$
-\frac {dx^1 \wedge dx^2 \wedge dx^3}{t}
= \cosh^2 \rho \sin \theta \, d\rho \wedge d\phi \wedge d\theta
= i Dy \cdot Y^{-1} = \frac {dS}{\|Y\|}
$$
($\tilde H$ is oriented as the boundary of $\{Y \in \BB M ;\: N(Y)<1 \}$).
\end{lem}


The functions (\ref{Y-harm}) can be regarded as functions on
$\BB M^+ = \{ Y \in \BB M ;\: N(Y)>0 \}$, and they satisfy $\square_{3,1}\phi=0$.
From (\ref{Delta-spherical}) and (\ref{Y-egf}) by direct computation we obtain:

\begin{lem}
We have the following two families of solutions of $\square_{3,1}\phi=0$
on $\BB M^+$:
$$
(r^2-t^2)^l \cdot  r^{-l-1} \cdot Y_l^m(\theta,\phi),
\qquad l=0,1,2,\dots, \quad -l \le m \le l,
$$
of homogeneity degree $l-1$ and
$$
(r^2-t^2)^{-l-1} \cdot r^l \cdot Y_l^m(\theta,\phi),
\qquad l=0,1,2,\dots, \quad -l \le m \le l,
$$
of homogeneity degree $-l-2$.
\end{lem}

To construct a basis of solutions of $\square_{3,1}\phi=0$ we differentiate
these two families with respect to $t$:
\begin{equation}  \label{sln}
\frac {\partial^k}{\partial t^k} 
(r^2-t^2)^l \cdot  r^{-l-1} \cdot Y_l^m(\theta,\phi)
\qquad \text{and} \qquad
\frac {\partial^k}{\partial t^k} 
(r^2-t^2)^{-l-1} \cdot r^l \cdot Y_l^m(\theta,\phi),
\end{equation}
where $l=0,1,2,\dots$, $-l \le m \le l$, $0 \le k \le 2l$. We will see soon
that these are two orthogonal families of functions on $\BB M^+$.
Let $T = t/r$, then we can rewrite the first family in (\ref{sln}) as
$$
\frac{d^k}{dT^k}(1-T^2)^l \Bigr|_{T= \frac tr}
\cdot r^{l-k-1} \cdot Y_l^m(\theta,\phi).
$$
Then we rewrite this expression using the associated Legendre functions:
$$
(-2)^l \frac{k! \cdot l!}{(2l-k)!} \cdot r^{-1} \cdot
(r^2 - t^2)^{(l-k)/2} \cdot P^{(l-k)}_l (t/r)
\cdot Y_l^m(\theta,\phi).
$$
Same can be done with the second family in (\ref{sln}).

\begin{thm}  \label{M-basis-thm}
We have two families of solutions of $\square_{3,1}\phi=0$ on $\BB M^+$:
\begin{equation}  \label{basis1}
f_{l,m,n}^+(Y)=
\sqrt{ \frac {(2l+1)\cdot (l-n)!}{(l+n)!}}
\cdot r^{-1} \cdot (r^2 - t^2)^{n/2} \cdot P^{(n)}_l (t/r)
\cdot Y_l^m(\theta,\phi)
\end{equation}
and
\begin{equation}  \label{basis2}
f_{l,m,n}^-(Y)=
\sqrt{ \frac {(2l+1)\cdot (l-n)!}{(l+n)!}}
\cdot r^{-1} \cdot (r^2 - t^2)^{-n/2} \cdot P^{(n)}_l (t/r)
\cdot Y_l^m(\theta,\phi),
\end{equation}
where $l=0,1,2,\dots$, $-l \le m \le l$, $0 \le n \le l$.

These families of solutions satisfy the following orthogonality relations:
$$
\bigl\langle f_{l,m,n}^+(Y), f_{l',-m',n'}^{\pm}(Y) \bigr\rangle =
-\bigl\langle f_{l,m,n}^-(Y), f_{l',-m',n'}^{\pm}(Y) \bigr\rangle =
(-1)^m \frac2{\pi i}
\delta_{l,l'} \cdot \delta_{m,m'} \cdot \delta_{n,n'} \cdot (1-\delta_{n,0})
$$
with respect to the bilinear pairing
\begin{equation}  \label{M-pairing}
\langle \phi_1 , \phi_2 \rangle = \frac 1{2\pi^2i} \int_{Y \in \tilde H}
(\deg \phi_1)(Y) \cdot \phi_2(Y) \,\frac{dS}{\|Y\|}.
\end{equation}

Also, these two families of solutions are related by involutions
$\phi(Y) \mapsto N(Y)^{-1} \cdot \phi\bigl(Y/N(Y)\bigr)$ and
$\phi(Y) \mapsto \frac{(-1)^l}{N(Y)} \cdot \phi(Y^{-1})$.
\end{thm}

\pf
The functions (\ref{basis1}) are homogeneous of degree $n-1$ and
the functions (\ref{basis2}) are homogeneous of degree $-(n+1)$;
if $n=0$, then $\deg f_{l,m,n}^{\pm}(Y) \equiv 0$.
To prove the orthogonality relations for $f^{\pm}_{l,m,n}(Y)$'s, we integrate in
coordinates (\ref{hyper-param}) using orthogonality identities
(\ref{Y-orthogonality}), (\ref{P-orthogonality})
and reduce to the following calculation:
\begin{multline}  \label{orthog-calc}
\int_{-\infty}^{\infty} 
(\cosh\rho)^{-1} \cdot P^{(n)}_l (\tanh \rho)
\cdot (\cosh\rho)^{-1} \cdot P^{(n')}_l (\tanh \rho)
\cdot \cosh^2 \rho \,d\rho  \\
= \int_{-\infty}^{\infty} 
P^{(n)}_l (\tanh \rho) \cdot P^{(n')}_l (\tanh \rho) \,d\rho
= \int_{-1}^1 P^{(n)}_l (x) \cdot P^{(n')}_l (x) \, \frac{dx}{1-x^2}
= \frac {(l+n)!}{n (l-n)!} \cdot \delta_{n,n'}
\end{multline}
(unless $n=n'=0$).

The last statement follows immediately from (\ref{P(-x)}) and
(\ref{Y-P-relation}).
\qed

\begin{rem}  \label{M-discr-ser-lim-rem}
When $n=0$, (\ref{basis1}) and (\ref{basis2}) yield the same functions.
We call these functions $f_{l,m,0}^{\pm}(Y)$ the limits of the discrete series
in $\BB M$. Note that $P^{(n)}_l$ and $P^{(-n)}_l$ are proportional according to
(\ref{m-m}), and so only $n \ge 0$ appear.
\end{rem}

\subsection{Space ${\cal D}^0_{\BB M}$ and $\mathfrak{sl}(4,\BB C)$-action}

Let ${\cal D}^0_{\BB M}$ denote the linear span of the functions (\ref{basis1})
and (\ref{basis2}).
In this subsection we characterize the space ${\cal D}^0_{\BB M}$ as a subspace
of harmonic polynomials on $\BB M$, show that it is preserved by the
$\mathfrak{sl}(4,\BB C)$-action and identify this representation with the
minimal representation of $SO(3,3) \simeq SL(2,\HR)/\{\pm 1\}$.

First we identify the linear span of the basis (\ref{basis1}).
Then the linear span of the basis (\ref{basis2}) is obtained by applying an
involution $\phi(Y) \mapsto N(Y)^{-1} \cdot \phi\bigl(Y/N(Y)\bigr)$ or
$\phi(Y) \mapsto \frac{(-1)^l}{N(Y)} \cdot \phi(Y^{-1})$.
Recall that $r= \sqrt{(x^1)^2+(x^2)^2+(x^3)^2}$ and consider algebraic functions
$\phi \in \BB C [t,x^1,x^2,x^3, r^{-1}]$ satisfying $\square_{3,1}\phi =0$.
We denote by $\BB C [t,x^1,x^2,x^3, r^{-2}]^{\le 0}$ the subalgebra of polynomials
in $\BB C [t,x^1,x^2,x^3, r^{-2}]$ spanned by monomials of non-positive degree
of homogeneity.

\begin{prop}
The linear span of the basis (\ref{basis1}) can be characterized as
\begin{equation}  \label{M-discr-span}
\bigl\{ \phi \in r^{-1} \cdot \BB C [t,x^1,x^2,x^3, r^{-2}]^{\le 0} ;\:
\square_{3,1}\phi =0 \bigr\}.
\end{equation}
\end{prop}

\pf
From the definition of the associate Legendre functions (\ref{P-def}),
(\ref{sln}) and Lemma \ref{rY-poly} we see that the basis functions
(\ref{basis1}) indeed lie in $r^{-1} \cdot \BB C [t,x^1,x^2,x^3, r^{-2}]^{\le 0}$.
Also, if we fix $l$ and $m$, there are exactly $l+1$ such basis functions
which are distinguished by their homogeneity degrees $d$, $-l-1 \le d \le -1$.

Since the polynomials $r^l \cdot Y_l^m(\theta,\phi)$'s form a basis of
polynomial functions on the sphere $S^2$ lying in the linear span of 
$x^1,x^2,x^3$, every function in $\BB C [t,x^1,x^2,x^3, r^{-1}]$ can be uniquely
expressed as a finite linear combination of the monomials
$$
t^ar^bY_l^m(\theta,\phi),
\qquad a,l=0,1,2,\dots, \quad b \in \BB Z, \quad -l \le m \le l.
$$
From (\ref{Delta-spherical}) and (\ref{Y-egf}) we obtain
$$
\square_{3,1} \bigl( t^ar^bY_l^m(\theta,\phi) \bigr)
= \bigl( b(b+1)-l(l+1) \bigr) t^ar^{b-2}Y_l^m(\theta,\phi)
- a(a-1)t^{a-2}r^bY_l^m(\theta,\phi).
$$
Note that the coefficient of $t^ar^{b-2}Y_l^m$ is zero if and only if
$b=l$ or $b=-l-1$, and the coefficient of $t^{a-2}r^bY_l^m$ is zero if
and only if $a=0$ or $a=1$.
This implies that we can find a function $g \in \BB C [t, r^{\pm1}]$ such that
$\square_{3,1} \bigl( g \cdot Y_l^m \bigr) =0$ by starting
with a monomial $t^ar^lY_l^m$ (or $t^ar^{-l-1}Y_l^m$) and inductively finding
the coefficients of $t^{a-2}r^{l+2}Y_l^m$, $t^{a-4}r^{l+4}Y_l^m$,...
(or $t^{a-2}r^{-l+1}Y_l^m$, $t^{a-4}r^{-l+3}Y_l^m$,...) until we zig-zag
to $t^0r^*Y_l^m$ or $t^1r^*Y_l^m$.
We can summarize this observation as follows:

\begin{lem}
Every function in the space
$\{ \phi \in \BB C [t,x^1,x^2,x^3, r^{-1}] ;\: \square_{3,1}\phi=0 \}$
is a finite linear combination of functions of the type
$g \cdot Y_l^m(\theta,\phi)$, where $g \in \BB C [t, r^{\pm1}]$ is such that
$\square_{3,1} \bigl( g \cdot Y_l^m(\theta,\phi) \bigr) =0$.
Moreover, for a fixed $Y_l^m$ and $d \in \BB Z$,
$$
\dim \biggl\{ g \in \BB C [t, r^{\pm1}] ;\:
\begin{matrix} \text{$g$ is homogeneous of degree $d$} \\
\text{and $\square_{3,1} \bigl( g \cdot Y_l^m \bigr) =0$}\end{matrix} \biggr\} =
\begin{cases}
2 & \text{if $d \ge l$,}  \\
1 & \text{if $-l-1 \le d \le l-1$,} \\
0 &\text{if $d < -l-1$}.
\end{cases}
$$
\end{lem}

This proves that functions (\ref{basis1}) span all of (\ref{M-discr-span}).
\qed

\begin{prop}  \label{algebra-action}
The Lie algebra
$\mathfrak{sl}(4,\BB C) = \mathfrak{sl}(2,\HC) \subset \mathfrak{gl}(2,\HC)$
acts on ${\cal D}^0_{\BB M}$, i.e. it preserves the space of finite
linear combinations of functions from (\ref{basis1}) and (\ref{basis2}).
Moreover, ${\cal D}^0_{\BB M}$ is generated by
$f^+_{0,0,0}(Y)=f^-_{0,0,0}(Y) = r^{-1}$.
\end{prop}

\pf
It is sufficient to find the actions of
$\begin{pmatrix} 0 & B \\ 0 & 0 \end{pmatrix} \in \mathfrak{sl}(2,\HC)$,
$B \in \HC$, since together with their conjugates by
$\begin{pmatrix} 0 & 1 \\ 1 & 0 \end{pmatrix} \in GL(2,\HC)$
they generate all of $\mathfrak{sl}(2,\HC)$.
Thus we need to find the actions of
$\frac{\partial}{\partial t}$, $\frac{\partial}{\partial x^1}$, 
$\frac{\partial}{\partial x^2}$, $\frac{\partial}{\partial x^3}$
on the functions
$r^{-1} \cdot (r^2 - t^2)^{\pm n/2} \cdot P^{(n)}_l (t/r) \cdot Y_l^m(\theta,\phi)$.
First we consider $\frac{\partial}{\partial t}$,
then from (\ref{P-recursive-1}) and (\ref{P-recursive-2}) we have:
\begin{multline*}
\frac{\partial}{\partial t} \Bigl( r^{-1} \cdot (r^2 - t^2)^{\pm n/2} \cdot
P^{(n)}_l (t/r) \cdot Y_l^m(\theta,\phi) \Bigr)  \\
= r^{-1} \cdot (r^2 - t^2)^{\pm n/2} \cdot 
\Bigl( r^{-1} \cdot \bigl( P^{(n)}_l \bigr)' (t/r) \mp
\frac {nt}{r^2 - t^2} P^{(n)}_l (t/r) \Bigr) \cdot Y_l^m(\theta,\phi)  \\
= r^{-1} \cdot (r^2 - t^2)^{(\pm n - 1)/2} \cdot P^{(n \mp 1)}_l(t/r) \cdot
Y_l^m(\theta,\phi) \cdot
\begin{cases}
(l+n)(l-n+1) & \text{if ``$+n$'',}\\
-1 & \text{if ``$-n$''.}
\end{cases}
\end{multline*}
Hence
\begin{align}
\frac{\partial}{\partial t} : \quad
f^+_{l,m,n}(Y) \quad &\mapsto \quad \sqrt{(l+n)(l-n+1)} \cdot f^+_{l,m,n-1}(Y),
\label{t-deriv1} \\
f^-_{l,m,n}(Y) \quad &\mapsto \quad - \sqrt{(l+n+1)(l-n)} \cdot f^-_{l,m,n+1}(Y).
\label{t-deriv2}
\end{align}

Next we consider $\frac{\partial}{\partial x^i}$. It is sufficient to show that
$\frac{\partial}{\partial x^i}$ acts on the bases (\ref{sln}).
By (\ref{Y-P-relation}) $Y_l^m(\theta,\phi)$ is proportional to
$P^{(m)}_l(\cos \theta) \cdot e^{-im\phi}$.
Since $\frac{\partial}{\partial x^i}$ commute with
$\frac{\partial}{\partial t}$, it is enough to show that
$$
\frac {\partial}{\partial x^i} \Bigl(
(r^2-t^2)^l \cdot  r^{-l-1} \cdot P^{(m)}_l(\cos \theta) \cdot e^{-im\phi} \Bigr)
\qquad \text{and} \qquad
\frac {\partial}{\partial x^i} \Bigl(
(r^2-t^2)^{-l-1} \cdot  r^l \cdot P^{(m)}_l(\cos \theta) \cdot e^{-im\phi} \Bigr)
$$
are finite linear combinations of elements from (\ref{sln}).
Using $x^3/r=\cos\theta$, (\ref{id1}) and (\ref{id2}), we get:
\begin{multline*}
e^{im\phi} \cdot \frac {\partial}{\partial x^3} \Bigl(
(r^2-t^2)^l \cdot  r^{-l-1} \cdot P^{(m)}_l(\cos \theta) \cdot e^{-im\phi} \Bigr) \\
=
\Bigl( 2l (r^2-t^2)^{l-1} \cdot r^{-l} - (l+1) (r^2-t^2)^l \cdot r^{-l-2} \Bigr)
\cdot \cos\theta \cdot P^{(m)}_l(\cos \theta)
\\
+ (r^2-t^2)^l \cdot r^{-l-2} \cdot \Bigl( (l+m)P^{(m)}_{l-1}(\cos \theta)
-l \cos\theta \cdot P^{(m)}_l(\cos\theta) \Bigr)  \\
=
(l+m) (r^2-t^2)^l \cdot r^{-l-2} \cdot P^{(m)}_{l-1}(\cos \theta) \\
+
r^{-l} \cdot \Bigl(\frac{2l}{2l+1} (r^2-t^2)^{l-1} - (r^2-t^2)^l \cdot r^{-2} \Bigr)
\cdot \bigl( (l-m+1) P^{(m)}_{l+1}(\cos\theta)
+ (l+m) P^{(m)}_{l-1}(\cos\theta) \bigr)  \\
= \frac {2l(l+m)}{2l+1} (r^2-t^2)^{l-1} \cdot r^{-l} \cdot
P^{(m)}_{l-1}(\cos \theta)  \\
+
\frac{l-m+1}{2l+1} \Bigl( 2l r^2 (r^2-t^2)^{l-1} - (2l+1) (r^2-t^2)^l \Bigr)
\cdot r^{-l-2} \cdot P^{(m)}_{l+1}(\cos \theta)  \\
= \frac{2l(l+m)}{2l+1} (r^2-t^2)^{l-1} \cdot r^{-l} \cdot P^{(m)}_{l-1}(\cos \theta)
+ \frac{l-m+1}{2(2l+1)(l+1)} \frac{\partial^2}{\partial t^2}
\Bigl( (r^2-t^2)^{l+1} \cdot r^{-l-2} \cdot P^{(m)}_{l+1}(\cos \theta) \Bigr).
\end{multline*}

Similarly, we compute
\begin{multline*}
\biggl(\frac {\partial}{\partial x^1}
+i\frac {\partial}{\partial x^2} \biggr)\Bigl(
(r^2-t^2)^l \cdot  r^{-l-1} \cdot P^{(m)}_l(\cos \theta) \cdot e^{-im\phi} \Bigr) \\
=
\frac {(l-m+1)(l-m+2)}{2(2l+1)(l+1)} \frac{\partial^2}{\partial t^2}
\Bigl( (r^2-t^2)^{l+1} \cdot r^{-l-2} \cdot P^{(m-1)}_{l+1}(\cos \theta)
\cdot e^{-i(m-1)\phi} \Bigr) \\
- \frac {2l(l+m)(l+m+1)}{2l+1} (r^2-t^2)^{l-1} \cdot r^{-l}
\cdot P^{(m-1)}_{l-1}(\cos\theta) \cdot e^{-i(m-1)\phi}
\end{multline*}
and
\begin{multline*}
\biggl(\frac {\partial}{\partial x^1}
-i\frac {\partial}{\partial x^2} \biggr)\Bigl(
(r^2-t^2)^l \cdot  r^{-l-1} \cdot P^{(m)}_l(\cos \theta) \cdot e^{-im\phi} \Bigr) \\
=
\frac {2l}{2l+1} (r^2-t^2)^{l-1} \cdot r^{-l} \cdot P^{(m+1)}_{l-1}(\cos\theta)
\cdot e^{-i(m+1)\phi}  \\
- \frac {l}{(2l+1)(l+1)} \frac{\partial^2}{\partial t^2} \Bigl( (r^2-t^2)^{l+1}
\cdot r^{-l-2} \cdot P^{(m+1)}_{l+1}(\cos \theta) \cdot e^{-i(m+1)\phi}\Bigr).
\end{multline*}

The derivatives of
$(r^2-t^2)^{-l-1} \cdot  r^l \cdot P^{(m)}_l(\cos \theta) \cdot e^{-im\phi}$
are computed in the same way.
Finally, it follows from these expressions for the action of
$\mathfrak{sl}(2,\HC)$ that $f^{\pm}_{0,0,0}$ generates all of ${\cal D}^0_{\BB M}$.
\qed

\begin{thm}  \label{M-discrete}
As a representation of $\mathfrak{sl}(4,\BB C)$, ${\cal D}^0_{\BB M}$ is the
minimal representation $(\varpi_{\BB R^{2,2}}^{\text{min}}, {\cal H})$
in the notations of \cite{KobO}. In particular, it is irreducible.
\end{thm}

\pf
Consider a function
\begin{equation}  \label{phi_0}
\phi_0(X)= \Bigl( \bigl(1+N(X)\bigr)^2+4x_{12}x_{21} \Bigr)^{-1/2},
\qquad X \in \HR.
\end{equation}
It is easy to see that $\phi_0$ is a real analytic function on $\HR$
satisfying $\square_{2,2} \phi_0=0$.
In fact, $\phi_0$ is the rescaled generating function $f_0$ from \cite{KobO},
Part III, equation (5.4.1) for $p=q=3$.
Hence it generates the minimal representation
$(\varpi_{\BB R^{2,2}}^{\text{min}}, {\cal H})$ of
$SO(3,3) \simeq SL(2,\HR)/\{\pm 1\}$.
Note that $-N(\tilde e_3 Y -1)^2$ is a negative real number only when
$\re(e_3 Y)=0$ and
\begin{equation}  \label{sqrt-n-square}
\bigl(-N(\tilde e_3 Y -1)^2 \bigr)^{1/2}
= -i \sgnx(Y) \cdot N(\tilde e_3 Y -1).
\end{equation}
Using (\ref{Y-tilde}) we obtain:
$$
\bigl( \pi_l^0(\gamma)\phi_0 \bigr)(Y)
= 2i \sgnx(Y) \cdot \bigl( -4(y_{22}-y_{11})^2 - 16y_{12}y_{21} \bigr)^{-1/2}
= \sgnx(Y) \cdot \frac i{2r},
$$
which is proportional to $\sgnx(Y) \cdot f^{\pm}_{0,0,0}(Y)$.
As far as representations of Lie algebras are concerned, we can drop
the factor $\sgnx(Y)$. By the previous proposition,
this function generates the whole span, and the result follows.
\qed

The $K$-types of the minimal representation
$(\varpi_{\BB R^{2,2}}^{\text{min}}, {\cal H})$ of $\mathfrak{sl}(4,\BB C)$ with
respect to $K=SO(3) \times SO(3)$ were identified in \cite{KobO}, Part I, as
$$
\bigoplus_l V_l \boxtimes V_l, \qquad l=0,1,2,\dots,
$$
where $V_l$ denotes the irreducible representation of $SO(3)$ of dimension
$l+1$.

\begin{prop}
The basis functions (\ref{basis1}) and (\ref{basis2}) are the $K$-finite
vectors of the minimal representation ${\cal D}^0_{\BB M}$ of
$\mathfrak{sl}(4,\BB C)$ with respect to $K=SO(3) \times SO(3)$.
\end{prop}

Define a nondegenerate symmetric bilinear pairing on the linear span of
(\ref{basis1}) and (\ref{basis2}) by declaring
$$
\bigl\langle f_{l,m,n}^+(Y), f_{l',-m',n'}^-(Y) \bigr\rangle_{min} =
\bigl\langle f_{l,m,n}^-(Y), f_{l',-m',n'}^+(Y) \bigr\rangle_{min} =
(-1)^m \frac2{\pi i} \delta_{l,l'} \cdot \delta_{m,m'} \cdot \delta_{n,n'},
$$
$$
\bigl\langle f_{l,m,n}^+(Y), f_{l',m',n'}^+(Y) \bigr\rangle_{min} =
\bigl\langle f_{l,m,n}^-(Y), f_{l',m',n'}^-(Y) \bigr\rangle_{min} = 0.
$$
In the second line we exclude $n'=0$ because $f_{l',m',0}^+(Y)=f_{l',m',0}^-(Y)$.
By Theorem \ref{M-basis-thm}, this pairing partially agrees with the
bilinear form (\ref{M-pairing}) up to a sign.

\begin{prop}  \label{invar-pairing-M}
This bilinear pairing on ${\cal D}^0_{\BB M}$ is $\mathfrak{gl}(2,\HC)$-invariant:
$$
\bigl\langle \pi^0_l(Z) \phi_1, \phi_2 \bigr\rangle_{min}
+ \bigl\langle \phi_1, \pi^0_r(Z) \phi_2 \bigr\rangle_{min} =0,
\qquad \forall Z \in \mathfrak{gl}(2,\HC).
$$
\end{prop}

\pf
It is sufficient to show invariance of $\langle\:,\:\rangle_{min}$
under $\mathfrak{sl}(2,\BB C)$ embedded into $\mathfrak{gl}(2,\HC)$
as in (\ref{SL(2,C)-embedding}),
$\begin{pmatrix} 0 & 1 \\ 1 & 0 \end{pmatrix} \in GL(2, \HC)$,
$\begin{pmatrix} 0 & 1 \\ 0 & 0 \end{pmatrix} \in \mathfrak{gl}(2, \HC)$
and scalar matrices.
The invariance with respect to the scalar matrices and
$\begin{pmatrix} 0 & 1 \\ 1 & 0 \end{pmatrix}$ is clear,
and the invariance with respect to
$\begin{pmatrix} 0 & 1 \\ 0 & 0 \end{pmatrix}$ follows from
(\ref{t-deriv1})-(\ref{t-deriv2}).

To show the $\mathfrak{sl}(2,\BB C)$-invariance, we observe that the action of
$\mathfrak{sl}(2,\BB C)$ commutes with $\deg$.
Hence, for each fixed $n$, $\mathfrak{sl}(2,\BB C)$
preserves the linear spans of $\{ f_{l,m,n}^+(Y) \}$ and $\{ f_{l,m,n}^-(Y) \}$,
$l=n,n+1,n+2, \dots$, $-l \le m \le l$.
Then the $\mathfrak{sl}(2,\BB C)$-invariance follows from Theorem
\ref{M-basis-thm} for $n>0$, and the case $n=0$ has to be considered separately.
Heuristically, the action of $\mathfrak{sl}(2,\BB C)$ on the linear span
of $\{ f_{l,m,0}^+(Y) \}$ is obtained from the action on $\{ f_{l,m,n}^{\pm}(Y) \}$,
$n$ fixed, by setting the value of the parameter $n=0$.
Since this action is described by algebraic formulas and preserves
$\langle\:,\:\rangle_{min}$ for all $n>0$,
it preserves $\langle\:,\:\rangle_{min}$ for $n=0$ as well.
\qed

\subsection{Extension of  ${\cal H}(\HR^+)$}  \label{extension}

Let $\B{{\cal D}^0}$ denote the unitary representation of $SL(2,\HR)$ in the
Hilbert space generated by the function $\phi_0(X)$ defined by (\ref{phi_0}).
This is the minimal representation $(\varpi_{\BB R^{2,2}}^{\text{min}}, {\cal H})$
of $SO(3,3) \simeq SL(2,\HR)/\{\pm 1\}$ in the notations of \cite{KobO}.

\begin{lem}  \label{tilde-phi_0}
The function
$$
\tilde\phi_0(X)= \frac{ \bigl( N(X)+1 \bigr) \cdot N(X-i)}
{\Bigl( \bigl(1+N(X)\bigr)^2+4x_{12}x_{21} \Bigr)^{3/2}},
\qquad X \in \HR
$$
also generates the minimal representation $\B{{\cal D}^0}$.
\end{lem}

\pf
Using (\ref{Y-tilde}) and (\ref{sqrt-n-square}) we obtain:
$$
\bigl( \pi_l^0(\gamma)\tilde\phi_0 \bigr)(Y) = \sgnx(Y) \cdot
\frac{16i (y_{22}-y_{11})}{\bigl( -4(y_{22}-y_{11})^2 - 16y_{12}y_{21} \bigr)^{3/2}}
= -\sgnx(Y) \cdot \frac{\cos\theta}{2r^2},
$$
which is proportional to $\sgnx(Y) \cdot f^-_{1,0,1}(Y)$.
It follows from Proposition \ref {algebra-action} that $\tilde\phi_0$ generates
the same representation as $\phi_0$ does.
\qed

\begin{prop}  \label{H-cont}
The intersection $\B{{\cal D}^0} \cap {\cal H}(\HR^+)$ is precisely the
continuous series component of ${\cal H}(\HR^+)$.
\end{prop}

\pf
Clearly, $\tilde\phi_0(X) \in {\cal H}(\HR^+)$.
Then Lemma \ref{tilde-phi_0} implies that $\phi_0$ cannot have a discrete
series component and $\B{{\cal D}^0} \cap {\cal H}(\HR^+)$ is
contained in the continuous series component of ${\cal H}(\HR^+)$.
To prove the other inclusion, it is sufficient to show that
$$
\bigl\langle
t^l_{\frac12\,\underline{\frac12}}(X), \tilde\phi_0(X) \bigr\rangle_1 \ne 0
\qquad \text{and} \qquad
\bigl\langle t^l_{0\,\underline{0}}(X), \phi_0(X) \bigr\rangle_1 \ne 0
\qquad \text{for all $l = -\frac12 + i\lambda$ with $\lambda \in \BB R^{\times}$.}
$$
In coordinates (\ref{param}) the restrictions of $\phi_0$ and $\tilde\phi_0$
to $SU(1,1)$ become
$$
\phi_0(\phi,\tau,\psi) = \Bigl( 2\cosh \frac{\tau}2 \Bigr)^{-1}
\qquad \text{and} \qquad
\tilde\phi_0(\phi,\tau,\psi) =
\frac1{4i} \bigl( e^{i\frac{\phi+\psi}2} + e^{-i\frac{\phi+\psi}2} \bigr)
\cdot \Bigl( \cosh \frac{\tau}2 \Bigr)^{-2}.
$$
Substituting $v=\sinh^2\frac{\tau}2$ and using a special case of an integral
formula 7.512(10) from \cite{GR}
$$
\int_0^{\infty} (1+x)^{-r} \,_2F_1(a,b;1;-x)\,dx =
\frac{\Gamma(a+r-1)\Gamma(b+r-1)}{\Gamma(r)\Gamma(a+b+r-1)}
$$
valid when $\re(a+r-1)>0$, $\re(b+r-1)>0$, we obtain
\begin{multline*}
\bigl\langle t^l_{\frac12\,\underline{\frac12}}(X), \tilde\phi_0(X) \bigr\rangle_1 =
-\frac{2l+1}{8i} \int_0^{\infty} \sinh \tau \cdot
\Bigl( \cosh \frac{\tau}2 \Bigr)^{-2}
\cdot \mathfrak{P}^l_{\frac12\,\frac12}(\cosh \tau)\,d\tau  \\
=
-\frac{2l+1}{4i} \int_0^{\infty} \sinh \frac{\tau}2 \cdot
\,_2F_1\Bigl(l+3/2,-l+1/2;1; -\sinh^2\frac{\tau}2 \Bigr)\,d\tau  \\
=
-\frac{2l+1}{4i} \int_0^{\infty}
(1+v)^{-\frac12} \cdot \,_2F_1(l+3/2,-l+1/2;1;-v)\,dv  \\
=
-\frac{2l+1}{4i} \cdot
\frac{\Gamma(l+1)\Gamma(-l)}{\Gamma(1/2)\Gamma(3/2)}
= \frac{2l+1}{2i \sin(\pi l)}
= - \frac{2l+1}{2i \cosh(\pi \im l)} \ne 0.
\end{multline*}
In the other case we need to deal with convergence issues, so we observe that
$$
\mathfrak{P}^l_{0\,0}(\cosh\tau) = \lim_{s \to 0^+} \mathfrak{P}^l_{0\,0}(\cosh\tau)
\cdot \Bigl( \cosh\frac{\tau}2 \Bigr)^{-2s}
$$
as a distribution in $\im l$. Substituting $v=\sinh^2\frac{\tau}2$ we obtain
\begin{multline*}
\bigl\langle t^l_{0\,\underline{0}}(X), \phi_0(X) \bigr\rangle_1 =
\lim_{s \to 0^+} -\frac{2l+1}4 \int_0^{\infty} \sinh \tau \cdot
\Bigl( \cosh \frac{\tau}2 \Bigr)^{-1-2s}
\cdot \mathfrak{P}^l_{0\,0}(\cosh \tau)\,d\tau  \\
=
\lim_{s \to 0^+} -\frac{2l+1}2 \int_0^{\infty} \sinh \frac{\tau}2 \cdot
\Bigl( \cosh \frac{\tau}2 \Bigr)^{-2s} \cdot
\,_2F_1\Bigl(l+1,-l;1; -\sinh^2\frac{\tau}2 \Bigr)\,d\tau  \\
=
\lim_{s \to 0^+} -\frac{2l+1}2 \int_0^{\infty}
(1+v)^{-\frac12-s} \cdot \,_2F_1(l+1,-l;1;-v)\,dv  \\
=
\lim_{s \to 0^+} -\frac{2l+1}2 \cdot
\frac{\Gamma(l+s+1/2)\Gamma(s-l-1/2)}{\bigl(\Gamma(s+1/2)\bigr)^2}  \\
=
-\frac{2l+1}2 \cdot
\frac{\Gamma(l+1/2)\Gamma(-l-1/2)}{\bigl(\Gamma(1/2)\bigr)^2}
= \frac1{\cos(\pi l)} =  \frac1{i\sinh(\pi \im l)} \ne 0.
\end{multline*}
\qed

\begin{rem}  \label{switch}
We can summarize the observation made at the beginning of
Subsection \ref{M-cont}, Theorem \ref{M-discrete} and Proposition \ref{H-cont}
as follows: The Cayley transform switches the discrete and continuous series
components of the spaces of harmonic functions on $\HR$ and $\BB M$.
\end{rem}

The Schwartz space ${\cal H}(\HR^+)$ ``almost'' contains the representations
$\B{{\cal D}^0}$, ${\cal D}^-$ and ${\cal D}^+$, meaning that ${\cal H}(\HR^+)$
is missing the limits of the discrete series and the functions on $\HR$
corresponding to the limits of the discrete series on $\BB M$
(see Remark \ref{M-discr-ser-lim-rem}).
From the point of view of representation theory it is natural to add these
missing functions to ${\cal H}(\HR^+)$. Then the resulting space is just a
direct sum of three irreducible representations of $SL(2,\HR)$.
Expansions of $\frac1{N(X-W)}$ contain these missing functions, which is
another reason for extending ${\cal H}(\HR^+)$.
We define an extended space
$$
\widehat{\cal H}(\HR^+) =
\B{{\cal D}^0} \oplus \B{{\cal D}^-} \oplus \B{{\cal D}^+}
\quad \supsetneq \quad {\cal H}(\HR^+),
$$
where $\B{{\cal D}^{\mp}}$ denote the Hilbert space completions of
${\cal D}^{\mp}$ with respect to an $\mathfrak{su}(2,2)$-invariant
inner product (which is unique up to scaling).
Similarly we can define
\begin{equation}  \label{H(M)}
{\cal H}(\BB M^+) = \pi_l^0(\gamma) \bigl( \B{{\cal D}^-} \bigr) \oplus
\pi_l^0(\gamma) \bigl( \B{{\cal D}^+} \bigr) \oplus
\bigl\{ \sgnx(Y) \cdot (\pi_l^0(\gamma)\phi)(Y);\:
\phi \in \B{{\cal D}^0} \bigr\}.
\end{equation}

Next we define a nondegenerate $\mathfrak{gl}(2,\HC)$-invariant
symmetric bilinear pairing $\langle\:,\:\rangle_{\widehat{\cal H}(\HR^+)}$ on
$\widehat{\cal H}(\HR^+)$. We declare
$$
\langle \B{{\cal D}^0}, \B{{\cal D}^-} \rangle_{\widehat{\cal H}(\HR^+)}
= \langle \B{{\cal D}^0}, \B{{\cal D}^+} \rangle_{\widehat{\cal H}(\HR^+)}
= \langle \B{{\cal D}^-},\B{{\cal D}^-} \rangle_{\widehat{\cal H}(\HR^+)}
= \langle \B{{\cal D}^+}, \B{{\cal D}^+} \rangle_{\widehat{\cal H}(\HR^+)} =0.
$$
Then we extend the pairing $\langle\:,\:\rangle_{{\cal D}^- \times {\cal D}^+}$
described in Proposition \ref{invar-pairing} by continuity to
$\B{{\cal D}^-} \oplus \B{{\cal D}^+}$.
Finally, we define the pairing on $\B{{\cal D}^0}$ by
$$
\bigl\langle \phi_1(X),\phi_2(X) \bigr\rangle_{\widehat{\cal H}(\HR^+)} =
\bigl\langle \sgnx(Y) \cdot (\pi_l^0(\gamma)\phi_1)(Y),
\sgnx(Y) \cdot (\pi_l^0(\gamma)\phi_2)(Y) \bigr\rangle_{min},
$$
where $\langle \:,\:\rangle_{min}$ is the pairing described in
Proposition \ref{invar-pairing-M} and extended by continuity to
the last summand of (\ref{H(M)}).

For example, in this light the matrix coefficient expansions given in
Proposition \ref{expansion2} mean that $-\frac 1{N(Z-W)}$ is the
reproducing kernel for $\B{{\cal D}^-}$ and $\B{{\cal D}^+}$ when,
respectively, $W \in \Gamma^-$ and $\Gamma^+$.
In other words, for each $W \in \Gamma^-$, the function $\frac 1{N(Z-W)}$
lies in the extended space $\widehat{\cal H}(\HR^+)$ and
$$
\phi(Z) \quad \mapsto \quad - \Bigl\langle \phi(Z) , \frac 1{N(Z-W)}
\Bigr\rangle_{\widehat{\cal H}(\HR^+)}
$$
is the projection of $\phi$ onto its holomorphic discrete series component.
Similarly, for each $W \in \Gamma^+$, the function
$\frac 1{N(Z-W)} \in \widehat{\cal H}(\HR^+)$ and
$$
\phi(Z) \quad \mapsto \quad - \Bigl\langle \phi(Z) , \frac 1{N(Z-W)}
\Bigr\rangle_{\widehat{\cal H}(\HR^+)}
$$
is the projection of $\phi$ onto its antiholomorphic discrete series component.

\subsection{Expansion of $\frac 1{N(Z-W)}$}

\begin{lem}
Let $Z,W \in \BB M^+$, write $Z = t \tilde e_0 + x^1 e_1 + x^2 e_2 + x^3 e_3$,
$W= \tilde t \tilde e_0 + \tilde x^1 e_1 + \tilde x^2 e_2 + \tilde x^3 e_3$,
then we have the following expansions:
\begin{equation}  \label{P(x)_1}
\sum_{m=-l}^l (-1)^m Y_l^{-m}(\theta,\phi) \cdot Y_l^m(\tilde\theta,\tilde\phi)
=
P_l \Bigl(\frac{x^1 \tilde x^1+x^2 \tilde x^2+x^3 \tilde x^3}{r \tilde r}\Bigr),
\end{equation}
\begin{equation}  \label{P(x)_2}
\sum_{n=-l}^l \frac{(l-n)!}{(l+n)!} \cdot
\frac {(r^2 - t^2)^{n/2}}{(\tilde r^2 - \tilde t^2)^{n/2}} \cdot
P_l^{(n)}(t/r) \cdot P_l^{(n)}(\tilde t / \tilde r)
=
P_l \Bigl( \frac {t \tilde t}{r \tilde r} +
\frac {r^2-t^2 + \tilde r^2- \tilde t^2}{2 r \tilde r} \Bigr),
\end{equation}
\begin{equation}  \label{Q(x)_2}
\sum_{n=-l}^l \frac{(l-n)!}{(l+n)!} \cdot
\frac {(r^2 - t^2)^{n/2}}{(\tilde r^2 - \tilde t^2)^{n/2}} \cdot
\tilde Q_l^{(n)}(t/r) \cdot P_l^{(n)}(\tilde t / \tilde r)
=
\tilde Q_l \Bigl( \frac {t \tilde t}{r \tilde r} +
\frac {r^2-t^2 + \tilde r^2- \tilde t^2}{2 r \tilde r} \Bigr).
\end{equation}
(For the last expansion we assume
$(2t \tilde t +r^2-t^2 + \tilde r^2- \tilde t^2)/ 2 r \tilde r$
lies in the domain of $\tilde Q_l$.)
\end{lem}

\pf
Since $\cos\theta = x^3/r$, $\cos\tilde\theta = \tilde x^3/\tilde r$ and
$$
\sin \theta \sin \tilde\theta \cos (\phi-\tilde\phi)
=
\sin \theta \cos \phi \sin \tilde\theta \cos \tilde\phi
+
\sin \theta \sin \phi \sin \tilde\theta \sin \tilde\phi
=
\frac {x^1}r \cdot \frac {\tilde x^1}{\tilde r} +
\frac {x^2}r \cdot \frac {\tilde x^2}{\tilde r},
$$
from (\ref{Y-P-relation}), (\ref{P-sum}) and (\ref{m-m}) we get (\ref{P(x)_1}).

Applying (\ref{P-sum}) with $s \in \BB C$ such that
$e^{is} = \sqrt{ \frac {r^2 - t^2}{\tilde r^2 - \tilde t^2}}$,
$\cos s = (e^{is}+e^{-is})/2$,
$\cos \theta = t/r$ and $\cos \tilde\theta = \tilde t / \tilde r$ we get:
\begin{multline*}
\sum_{n=-l}^l \frac{(l-n)!}{(l+n)!} \cdot
\frac {(r^2 - t^2)^{n/2}}{(\tilde r^2 - \tilde t^2)^{n/2}} \cdot 
P_l^{(n)}(t/r) \cdot P_l^{(n)}(\tilde t / \tilde r)
=
P_l (\cos \theta \cos \tilde\theta + \sin \theta \sin \tilde\theta \cos \phi)\\
=
P_l \Biggl(\frac tr \frac {\tilde t}{\tilde r}
+ \frac 12 \sqrt{1 - \frac {t^2}{r^2}} \sqrt{1 - \frac {\tilde t^2}{\tilde r^2}}
\biggl( \sqrt{ \frac {r^2 - t^2}{\tilde r^2 - \tilde t^2}} +
\sqrt{ \frac {\tilde r^2 - \tilde t^2}{r^2 - t^2}} \biggr) \Biggr) \\
= P_l \Bigl( \frac {t \tilde t}{r \tilde r} +
\frac {|r^2-t^2| + |\tilde r^2- \tilde t^2|}{2 r \tilde r} \Bigr)
= P_l \Bigl( \frac {t \tilde t}{r \tilde r} +
\frac {r^2-t^2 + \tilde r^2- \tilde t^2}{2 r \tilde r} \Bigr),
\end{multline*}
which proves (\ref{P(x)_2}).

Using (\ref{Q-sum}) instead of (\ref{P-sum}) we prove the last expansion
(\ref{Q(x)_2}). Strictly speaking, (\ref{Q-sum}) applies only if
$\sqrt{ \frac {r^2 - t^2}{\tilde r^2 - \tilde t^2}} = e^{is}$ with $s \in \BB R$,
but both sides depend on $s$ analytically.
\qed

\begin{prop}
Let $W \in \BB M^+$ and
$Z \in \{\BB M^+ +ae_0; \: a\in \BB R,\: a \ne 0 \} \subset \HC$.
Write $Z = t \tilde e_0 + x^1 e_1 + x^2 e_2 + x^3 e_3$ with $t \in \BB C$,
$W= \tilde t \tilde e_0 + \tilde x^1 e_1 + \tilde x^2 e_2 + \tilde x^3 e_3$,
then we have the following expansions:
\begin{multline*}
\sum_{l=0}^{\infty} \frac{2l+1}{2r\tilde r} \Biggl( \sum_{m,n=-l}^l
(-1)^m (i \cdot \sgn \im t)^n \frac{(l-n)!}{(l+n)!} \cdot
\frac{(r^2-t^2)^{n/2}}{(\tilde r^2 - \tilde t^2)^{n/2}} \\
\times Q_l^{(n)}(t/r) \cdot Y_l^{-m}(\theta,\phi) \cdot
P_l^{(n)} (\tilde t / \tilde r) \cdot Y_l^m(\tilde\theta,\tilde\phi) \Biggr)
= \frac 1{N(Z-W)},
\end{multline*}
provided that $\tilde t > \re t$.
Similarly, if $Z = t \tilde e_0 + x^1 e_1 + x^2 e_2 + x^3 e_3\in \BB M^+$ and
$W = \tilde t \tilde e_0 + \tilde x^1 e_1 + \tilde x^2 e_2 + \tilde x^3 e_3
\in \{\BB M^+ +ae_0; \: a\in \BB R,\: a \ne 0 \}$ with $\tilde t \in \BB C$,
then
\begin{multline*}
\sum_{l=0}^{\infty} \frac{2l+1}{2r\tilde r} \Biggl( \sum_{m,n=-l}^l
(-1)^m (i \cdot \sgn \im \tilde t)^n \frac{(l-n)!}{(l+n)!} \cdot
\frac{(r^2 - t^2)^{n/2}}{(\tilde r^2 - \tilde t^2)^{n/2}} \\
\times P_l^{(n)}(t/r) \cdot Y_l^{-m}(\theta,\phi) \cdot
Q_l^{(n)}(\tilde t/\tilde r) \cdot Y_l^m(\tilde\theta,\tilde\phi) \Biggr)
= \frac 1{N(Z-W)},
\end{multline*}
provided that $t > \re \tilde t$.
\end{prop}

\pf
Letting $t_0=\re t$, $a=\im t$ and using (\ref{P(x)_2}), (\ref{Q(x)_2}),
(\ref{Q-tilde}), the first sum reduces to
\begin{multline*}
\sum_{n=-l}^l (i \cdot \sgn a)^n \frac{(l-n)!}{(l+n)!} \cdot
\frac{(r^2-t^2)^{n/2}}{(\tilde r^2-\tilde t^2)^{n/2}} \cdot
Q_l^{(n)}\Bigl(\frac{t_0+ia}r\Bigr) \cdot P_l^{(n)}(\tilde t / \tilde r) \\
=
\sum_{n=-l}^l \frac{(l-n)!}{(l+n)!} \cdot
\frac {(r^2-t^2)^{n/2}}{(\tilde r^2-\tilde t^2)^{n/2}}
\cdot \biggl( \tilde Q_l^{(n)}\Bigl(\frac{t_0+ia}r\Bigr)
- \sgn a \cdot \frac{\pi i}2 P_l^{(n)}\Bigl(\frac{t_0+ia}r\Bigr)\biggr)
\cdot P_l^{(n)}(\tilde t / \tilde r)  \\
=
\tilde Q_l \Bigl( \frac {t \tilde t}{r \tilde r} +
\frac {r^2-t^2 + \tilde r^2- \tilde t^2}{2 r \tilde r} \Bigr)
- \sgn a \cdot \frac{\pi i}2 P_l \Bigl( \frac {t \tilde t}{r \tilde r} +
\frac {r^2-t^2 + \tilde r^2- \tilde t^2}{2 r \tilde r} \Bigr)  \\
=
Q_l \Bigl( \frac {t \tilde t}{r \tilde r} +
\frac {r^2-t^2 + \tilde r^2- \tilde t^2}{2 r \tilde r} \Bigr),
\end{multline*}
provided that 
$$
\im \Bigl( \frac {t \tilde t}{r \tilde r} +
\frac {r^2-t^2 + \tilde r^2- \tilde t^2}{2 r \tilde r} \Bigr)
= a \frac{\tilde t-t_0}{r\tilde r}
$$
has the same sign as $a$.
Using 
$$
N(Z-W) = 
2t \tilde t + r^2-t^2 + \tilde r^2- \tilde t^2
- 2(x^1 \tilde x^1+x^2 \tilde x^2+x^3 \tilde x^3),
$$
(\ref{PQ-sum}) and (\ref{P(x)_1}) we see that
$$
\sum_{l=0}^{\infty} \frac{2l+1}{2r\tilde r}
P_l \Bigl(\frac {x^1 \tilde x^1+x^2 \tilde x^2+x^3 \tilde x^3}{r \tilde r} \Bigr)
\cdot Q_l \Bigl( \frac {t \tilde t}{r \tilde r} +
\frac {r^2-t^2 + \tilde r^2- \tilde t^2 +\epsilon^2}{2 r \tilde r} \Bigr)
= \frac 1{N(Z-W)}.
$$
The other expansion is proved by switching the variables $Z$ and $W$.
\qed

Let $Z, W \in \BB M^+$ and write $Z = t \tilde e_0 + x^1 e_1 + x^2 e_2 + x^3 e_3$,
$W= \tilde t \tilde e_0 + \tilde x^1 e_1 + \tilde x^2 e_2 + \tilde x^3 e_3$.
Then
$$
\frac 1{N(Z-W \pm i\epsilon \tilde e_0)}
= \frac 1{N(Z-W) + \epsilon^2 \mp 2i\epsilon(t-\tilde t)}.
$$
Letting $\epsilon \to 0^+$, from the above expansions and (\ref{Q-Q=Pn})
we obtain
\begin{multline*}
\pi i \sum_{l=0}^{\infty} \frac{2l+1}{2r\tilde r} \Biggl( \sum_{m,n=-l}^l
(-1)^m \frac{(l-n)!}{(l+n)!} \cdot
\frac{(r^2 - t^2)^{n/2}}{(\tilde r^2 - \tilde t^2)^{n/2}} \cdot 
P_l^{(n)}(t/r) \cdot Y_l^{-m}(\theta,\phi) \cdot
P_l^{(n)}(\tilde t / \tilde r) \cdot Y_l^m(\tilde\theta,\tilde\phi) \Biggr)  \\
=
\lim_{\epsilon \to 0^+} \biggl( \frac 1{N(Z-W)+\epsilon^2-2i\epsilon|t-\tilde t|}
- \frac 1{N(Z-W)+\epsilon^2+2i\epsilon|t-\tilde t|} \biggr).
\end{multline*}
As distributions on the hyperboloids $\{ N(Z)=const \}$, the limits
$$
\lim_{\epsilon \to 0^+} \frac1{N(Z-W)+\epsilon^2 \pm 2i\epsilon|t-\tilde t|}
\qquad \text{are equivalent to} \qquad \frac1{N(Z-W) \pm i0}.
$$
Thus we can formally write
\begin{multline}  \label{1/N-formal-expansion}
\pi i \sum_{l=0}^{\infty}  \frac{2l+1}{2r\tilde r}
\Biggl( \sum_{m,n=-l}^l (-1)^m \frac{(l-n)!}{(l+n)!} \cdot
\frac{(r^2 - t^2)^{n/2}}{(\tilde r^2 - \tilde t^2)^{n/2}} \cdot 
P_l^{(n)}(t/r) \cdot Y_l^{-m}(\theta,\phi) \cdot
P_l^{(n)}(\tilde t / \tilde r) \cdot Y_l^m(\tilde\theta,\tilde\phi) \Biggr)  \\
= \frac1{N(Z-W)-i0} - \frac1{N(Z-W)+i0}.
\end{multline}

\subsection{The Discrete Series Projector on $\BB M$}

Using expansions of $\frac1{N(Z-W) \pm i0}$ obtained in previous subsection,
we obtain projectors onto the discrete series component on $\BB M$.
Recall that the space of harmonic functions ${\cal H}(\BB M^+)$ was defined in
(\ref{H(M)}). Let $\tilde H_R = \{ Y \in \BB M ;\: N(Y)=R^2 \}$, $R>0$, and
define an operator on ${\cal H}(\BB M^+)$ by
\begin{multline}  \label{M-discr-oper}
\bigl( \S_R^{\BB M} \phi \bigr)(W) =
\lim_{\epsilon \to 0^+} \frac 1{2\pi^2i} \int_{Y \in \tilde H_R}
\biggl( \frac 1{N(Y-W) + i\epsilon} - \frac 1{N(Y-W) - i\epsilon} \biggr)
\cdot \deg\phi(Y) \,\frac{dS}{\|Y\|} \\
= \lim_{\epsilon \to 0^+} \frac 1{2\pi^2i}
\int_{Y \in \tilde H_R} \biggl( \frac {\deg\phi(Y)}{N(Y-W) + i\epsilon}
- \frac {\deg\phi(Y)}{N(Y-W) - i\epsilon} \biggr)
\cosh^2 \rho \sin \theta \, d\rho d\theta d\phi,
\end{multline}
$W \in \BB M^+$.

\begin{thm}  \label{M-discrete-proj}
The operator $\S_R^{\BB M}$ annihilates the continuous series component
of ${\cal H}(\BB M^+)$ and acts on the discrete series component of
${\cal H}(\BB M^+)$ by sending
\begin{align*}
f^+_{l,m,n}(Y) \quad &\mapsto \quad
- \bigl( f^+_{l,m,n}(W) + R^{2n} \cdot f^-_{l,m,n}(W) \bigr), \\
f^-_{l,m,n}(Y) \quad &\mapsto \quad
R^{-2n} \cdot f^+_{l,m,n}(W) + f^-_{l,m,n}(W)
\end{align*}
if $n \ne 0$ and annihilates $f^{\pm}_{l,m,0}(Y)$.
\end{thm}

\begin{rem}  \label{M-loc-rem}
We can ``localize'' the difference of the integrals at the set where $N(Y-W)=0$.
Thus the operator (\ref{M-discr-oper})
depends only on the values of $\deg\phi(Y)$, where $Y$ ranges over a
two-dimensional hyperboloid (or cone if $W\in \tilde H_R$)
$$
\{Y \in \BB M;\: N(Y)=R^2, \: N(Y-W)=0\}.
$$
\end{rem}

\pf
The discrete series on $\BB M^+$ is spanned by (\ref{basis1}) and
(\ref{basis2}). Hence it is enough to check the statement for the functions
$\phi(Y) =
r^{-1} \cdot (r^2 - t^2)^{\pm n/2} \cdot P^{(n)}_l (t/r) \cdot Y_l^m(\theta,\phi)$,
where $l=1,2,3,\dots$, $-l \le m \le l$, $1 \le n \le l$.
First we expand using (\ref{PQ-sum}) and (\ref{P(x)_1}):
\begin{multline*}
\frac 1{N(Y-W) \pm i\epsilon}
=
\sum_{l=0}^{\infty} \frac{2l+1}{2r\tilde r}
Q_l \Bigl( \frac {t \tilde t}{r \tilde r} +
\frac {r^2-t^2 + \tilde r^2- \tilde t^2}{2 r \tilde r}
\pm \frac{i\epsilon}{2r\tilde r} \Bigr) \cdot
P_l \Bigl(\frac {x^1 \tilde x^1+x^2 \tilde x^2+x^3 \tilde x^3}{r \tilde r} \Bigr)
\\
=
\sum_{l=0}^{\infty} \frac{2l+1}{2r\tilde r}
Q_l \Bigl( \frac {t \tilde t}{r \tilde r} +
\frac {r^2-t^2 + \tilde r^2- \tilde t^2}{2 r \tilde r}
\pm \frac{i\epsilon}{2r\tilde r} \Bigr)
\biggl( \sum_{m=-l}^l (-1)^m Y_l^{-m}(\theta,\phi) \cdot
Y_l^m(\tilde\theta,\tilde\phi) \biggr).
\end{multline*}
We have $\deg \phi(Y) = \pm n \phi(Y)$.
By the orthogonality relation (\ref{Y-orthogonality}),
after we integrate out $\phi$ and $\theta$ we are left with
integrating over $\rho$
\begin{multline*}
\pm \frac{2\pi n}{r^2 \tilde r} \biggl( Q_l \Bigl( \frac{t\tilde t}{r\tilde r}
+ \frac{r^2-t^2 + \tilde r^2- \tilde t^2}{2 r \tilde r}
+ \frac{i\epsilon}{r \tilde r} \Bigr)
- Q_l \Bigl( \frac{t \tilde t}{r \tilde r} +
\frac {r^2-t^2 + \tilde r^2- \tilde t^2}{2 r \tilde r}
- \frac{i\epsilon}{r \tilde r} \Bigr) \biggr)  \\
\times
(r^2 - t^2)^{\pm n/2} \cdot P^{(n)}_l (t/r) \cdot Y_l^m(\tilde\theta,\tilde\phi).
\end{multline*}
As $\epsilon \to 0^+$, by (\ref{Q-Q=Pn}) we get
$$
\mp 2\pi^2i \cdot \frac{nR^{\pm n}}{r^2 \tilde r}
P_l\Bigl(\frac{t \tilde t}{r \tilde r}
+ \frac {r^2-t^2 + \tilde r^2- \tilde t^2}{2r\tilde r}\Bigr)
\cdot P^{(n)}_l (t/r) \cdot Y_l^m(\tilde\theta,\tilde\phi).
$$
Then using (\ref{m-m}), (\ref{P(x)_2}) and our earlier calculation
(\ref{orthog-calc}) we can integrate over $\rho$ to get
$$
\mp 2\pi^2 i \cdot \tilde r^{-1} \cdot R^{\pm n} \cdot P^{(n)}_l(\tilde t/\tilde r)
\cdot Y_l^m(\tilde\theta,\tilde\phi) \cdot \bigl(
R^{-n}(\tilde r^2 - \tilde t^2)^{n/2} + R^n(\tilde r^2 - \tilde t^2)^{- n/2} \bigr).
$$

Finally, the fact that the continuous series on $\BB M$ gets annihilated
follows from Corollary \ref{M-cont-cor} or, alternatively, since
our expansion (\ref{1/N-formal-expansion}) contains the discrete series only.
\qed

We can rewrite the expression $\frac 1{N(Y-W)+i0} - \frac 1{N(Y-W)-i0}$
geometrically as follows.

\begin{lem}
Fix $W \in \BB M^+$ and let $Y$ range over the hyperboloid $\tilde H_R$.
Fix a number $\xi \in \BB C \setminus \BB R$ and let
$W'=(1+\epsilon\xi)W \in \HC$, $\epsilon>0$, then
$N(Y-W')$ is never zero for all $\epsilon$ sufficiently small.
\end{lem}

\pf
Let $A=\tr(YW^+)$, then
$$
N(Y-W') = N(Y) + (1+\epsilon\xi)^2 N(W) - (1+\epsilon\xi) A.
$$
Write $\xi=a+ib$, $a,b \in \BB R$, then the imaginary part of $N(Y-W')$
is zero if and only if
$$
2b (\epsilon^2 a + \epsilon) N(W) - \epsilon b A =0
\qquad\Longleftrightarrow\qquad
A = 2(\epsilon a + 1) N(W).
$$
Then the real part of $N(Y-W')$ is
$$
N(Y) + \bigl( 1+2\epsilon a + \epsilon^2(a^2-b^2) \bigr) N(W)
- 2(1+\epsilon a)^2 N(W)
= R^2 - ( 1 + 2\epsilon a + \epsilon^2 |\xi|^2) N(W) \ne 0
$$
for all $\epsilon \ll 1$.
\qed

When $N(Y-W)=0$,
$N(Y-W')=\epsilon\xi \bigl( N(W)-N(Y) \bigr) + \epsilon^2 \xi^2 N(W)$.
Choose $\xi_1, \xi_2 \in \BB C$ with $\im(\xi_1)>0$ and $\im(\xi_2)<0$, then
\begin{multline*}
\lim_{\epsilon \to 0^+} \sgn\bigl( N(W)-R^2 \bigr) \cdot \int_{Y \in \tilde H_R}
\biggl( \frac 1{N(Y-W) + i\epsilon} - \frac 1{N(Y-W) - i\epsilon} \biggr)
\cdot \deg\phi(Y) \,\frac{dS}{\|Y\|} \\
= \lim_{\epsilon \to 0^+}
\int_{Y \in \tilde H_R} \biggl( \frac 1{N\bigl( Y-(1+\epsilon\xi_1)W \bigr)}
- \frac 1{N\bigl( Y-(1+\epsilon\xi_2)W \bigr)} \biggr)
\cdot \deg\phi(Y) \,\frac{dS}{\|Y\|}.
\end{multline*}

We finish this subsection with an analogue of Theorem \ref{M-discrete-proj}
for left-regular functions.
Of course, similar result holds for right-regular functions.
For convergence reasons we need to exclude the functions with components
containing $f_{l,m,0}^{\pm}(Y)$. Thus we define $\BB S(\BB M^+)$ to be the space
of $\BB S$-valued left-regular functions on $\BB M^+$ with both components
orthogonal to all $f_{l,m,0}^{\pm}(Y)$'s with respect to
$\langle\:,\:\rangle_{min}$. Let $\P^{discr(\BB M)}_0$ and $\P^{discr(\BB M)}_{\infty}$
be the projections onto the discrete series components of $\BB S(\BB M^+)$
quasi-regular at $0$ and $\infty$ respectively.

\begin{thm}  \label{M-discr-proj-reg}
Let $f \in \BB S(\BB M^+)$ a left-regular function and $W \in \BB M^+$, then
\begin{multline*}
\bigl( \P^{discr(\BB M)}_0(f) - \P^{discr(\BB M)}_{\infty}(f) \bigr)(W)  \\
=\lim_{\epsilon \to 0^+} \frac i{2\pi^2} \int_{Y \in \tilde H_R}
\biggl( \frac{(Y-W)^+}{(N(Y-W) +i\epsilon)^2}
- \frac{(Y-W)^+}{(N(Y-W) -i\epsilon)^2} \biggr) \cdot Dy \cdot f(Y).
\end{multline*}
\end{thm}

\pf
Suppose that $f: \BB M^+ \to \BB S$ is left-regular and homogeneous of degree
$n-2$ (i.e. $(\deg + 1)f= nf$) with $n \ne 0$.
Let $\phi(Y)=Yf(Y): \BB M^+ \to \BB S$ be an $\BB S$-valued function,
then $\deg \phi = n \phi$, each coordinate component of $\phi$ is a harmonic
function and ``restricting'' (\ref{deg-nabla}) to $\BB M$ we obtain
\begin{equation}  \label{deg-nabla-M}
\nabla_{\BB M} \phi = (Y^+ \nabla^+_{\BB M} + \nabla_{\BB M} Y) f
= 2(\deg + 1)f = 2nf \ne 0.
\end{equation}
For concreteness, let us assume that $n>0$.
(If $f$ is quasi-regular at infinity, the signs
get reversed when we apply Theorem \ref{M-discrete-proj}.) Then
\begin{multline*}
\phi(W) + R^{2n} \cdot N(W)^{-1} \cdot \phi\bigl(W/N(W)\bigr) \\
= \lim_{\epsilon \to 0^+} \frac i{2\pi^2} \int_{Y \in \tilde H_R}
\biggl( \frac 1{N(Y-W) +i\epsilon} - \frac 1{N(Y-W) -i\epsilon} \biggr)
\cdot (\deg_Y \phi)(Y) \,\frac{dS}{\|Y\|}  \\
= \lim_{\epsilon \to 0^+} \frac i{2\pi^2} \int_{Y \in \tilde H_R}
\biggl( \frac 1{N(Y-W) +i\epsilon} - \frac 1{N(Y-W) -i\epsilon} \biggr)
\cdot Dy \cdot (\deg_Y +1) (Y^{-1}\phi)(Y)  \\
= \lim_{\epsilon \to 0^+} \frac i{4\pi^2} \int_{Y \in \tilde H_R}
\biggl( \frac 1{N(Y-W) +i\epsilon} - \frac 1{N(Y-W) -i\epsilon} \biggr)
\cdot Dy \cdot \nabla_{\BB M} \phi(Y).
\end{multline*}
Differentiating with respect to $W$ we get
\begin{multline*}
\nabla_{\BB M} \Bigl(
\phi(W) + R^{2n} \cdot N(W)^{-1} \cdot \phi\bigl(W/N(W)\bigr) \Bigr) \\
= \lim_{\epsilon \to 0^+} \frac i{2\pi^2} \int_{Y \in \tilde H_R} \biggl(
\frac{(Y-W)^+}{(N(Y-W) +i\epsilon)^2} - \frac{(Y-W)^+}{(N(Y-W) -i\epsilon)^2}
\biggr) \cdot Dy \cdot \nabla_{\BB M} \phi(Y).
\end{multline*}
Finally, by (\ref{deg-nabla-M}),
\begin{multline*}
\nabla_{\BB M} \Bigl( N(W)^{-1} \cdot \phi\bigl(W/N(W)\bigr) \Bigr) \\
= \frac1{N(W)^2} \cdot (\nabla_{\BB M} \phi)\Bigl(\frac{W}{N(W)}\Bigr)
- 2 \frac{W^{-1}}{N(W)} \cdot \phi\Bigl(\frac{W}{N(W)}\Bigr)
- 2 \frac{W^{-1}}{N(W)^2} \cdot (\operatorname{deg} \phi)
\Bigl(\frac{W}{N(W)}\Bigr)  \\
= \frac1{N(W)^2} \bigl( \nabla_{\BB M} \phi - 2f - 2\deg f \bigr)
\bigl(W/N(W)\bigr) =0.
\end{multline*}
\qed

\section{Separation of the Series and the Plancherel Measure of $SU(1,1)$}

\subsection{Convergence and Equivariant Properties of Poisson Integrals}

In this section we regularize the Poisson and Cauchy-Fueter integrals by
replacing $N(X-W)$ in the denominator with $N(X-W) \pm i\epsilon$ and
letting $\epsilon \to 0^+$. First we prove that the resulting limits converge.
Recall that $H_R = \{ X \in \HR ;\: N(X) = R^2 \} \subset \HR^+$, where $R>0$.

\begin{prop}  \label{main-prop}
Fix an $R>0$ and a positive integer $n$, then, for any Schwartz function
$\psi$ on $H_R$ and any $W \in \HR^+$ with $N(W) \ne R^2$, the limit
\begin{equation}  \label{tilde-psi-n}
\tilde \psi(W) = \lim_{\epsilon \to 0^+} \int_{X \in H_R}
\frac{\psi(X)}{\bigl( N(X-W) + i\epsilon \bigr)^n} \,\frac{dS}{\|X\|}
\end{equation}
determines a real-analytic function of $W$.
When $n=1$, $\square_{2,2} \tilde\psi = 0$.
Similar statements hold for $\epsilon \to 0^-$ as well.
\end{prop}

\begin{rem}
The function $\psi$ itself need not satisfy $\square_{2,2} \psi = 0$.
The limits as $\epsilon \to 0^+$ and $\epsilon \to 0^-$ can yield
different functions.
We will see later that the difference between these two limits corresponds
to the continuous series component of $\psi$.
\end{rem}

\pf
This proposition would be trivial if
$\supp \psi \cap \{ X \in H_R ;\: N(X-W)=0\} =\varnothing$.
On the other hand, we will rewrite the integral so that the only
``problematic'' points are the critical points of $N(X-W)$ restricted
to $H_R$ (as a function of $X$).
It is easy to see that the only such critical points are
$X_{crit}^{\pm} = \pm \frac R{\sqrt{N(W)}} W$,
and $N(X_{crit}^{\pm}-W) \ne 0$.
Then we apply a partition of unity argument to break up the function
$\psi$ into $\psi_0 + \psi_c$ so that $\psi_0$ is supported away from
$\{ X \in H_R ;\: N(X-W)=0\}$ and $\psi_c$ is supported away from $X_{crit}^{\pm}$.

Let $\eta : \BB R \to \BB R$ be a smooth compactly supported function such that
$\eta(t)=1$ for $t$ near $0$ and $\eta(t)=0$ for $t$ near
$N(X_{crit}^{\pm}-W)$.
Then the composition function $\eta \bigl(N(X-W) \bigr)$ is one
in a neighborhood of $\{ X \in H_R ;\: N(X-W)=0\}$ and zero in a
neighborhood of $X_{crit}^{\pm}$. Define
$$
\psi_0 (X) = \bigl( 1- \eta\bigl(N(X-W)\bigr) \bigr) \cdot \psi(X),  \qquad
\psi_c (X) = \eta\bigl(N(X-W)\bigr) \cdot \psi (X),
$$
so that $\psi = \psi_0 + \psi_c$, and set
$$
\psi_{0,\epsilon}(W) =
\int_{X \in  H_R} \frac{\psi_0(X)}{\bigl( N(X-W) + i\epsilon \bigr)^n}
\,\frac{dS}{\|X\|},  \qquad
\psi_{c,\epsilon}(W) =
\int_{X \in H_R} \frac{\psi_c(X)}{\bigl( N(X-W) + i\epsilon \bigr)^n}
\,\frac{dS}{\|X\|}.
$$
Clearly, the limit
$$
\tilde\psi_0(W) = \lim_{\epsilon \to 0^+} \psi_{0,\epsilon}(W)
= \int_{X \in H_R} \frac {\psi_0(X)}{N(X-W)^n} \,\frac {dS}{\|X\|}
$$
is a well defined real-analytic function of $W$.

Next we analyze $\psi_{c,\epsilon}(W)$.
Consider the case $n=1$ first and observe that
$$
\frac 1{N(X-W) + i\epsilon} = -i \int_0^{\infty}
\exp \bigl( it \bigr( N(X-W) + i\epsilon \bigr) \bigr) \,dt,
$$
and the integral converges absolutely since the real
part of the exponent is $-\epsilon t <0$.
Thus we can rewrite $\psi_{c,\epsilon}(W)$ as
$$
\psi_{c,\epsilon}(W) = -i \int_0^{\infty} \biggl( \int_{X \in H_R}
\exp \bigl( it \bigl( N(X-W) + i\epsilon \bigr) \bigr) \cdot \psi_c(X)
\,\frac {dS}{\|X\|} \biggr) \,dt,
$$
where we are allowed to switch the order of integration because
$\bigl| \exp \bigl( it \bigl( N(X-W) + i\epsilon \bigr) \bigr) \bigr| \le 1$
and the function $ \psi_c$ is absolutely integrable on $H_R$ with respect
to the measure $\frac {dS}{\|X\|}$.
Then the limit
$$
\tilde\psi_c(W) = \lim_{\epsilon \to 0^+} \psi_{c,\epsilon}(W)
= -i \int_0^{\infty} \biggl( \int_{X \in H_R}
\exp \bigl( it N(X-W) \bigr) \cdot \psi_c(X) \,\frac {dS}{\|X\|} \biggr) \,dt
$$
is a well defined real-analytic function of $W$.
(The integral converges by a standard integration by parts argument;
this is where the property $X_{crit}^{\pm} \notin \supp \psi_c$ is used.)

If $n>1$, then we can apply the same argument with
$$
\frac 1{\bigl( N(X-W) + i\epsilon \bigr)^n}
= (-i)^n \int_0^{\infty} \dots \int_0^{\infty}
\exp \bigl( i(t_1+ \dots +t_n) \bigr( N(X-W) + i\epsilon \bigr) \bigr)
\,dt_1 \dots dt_n.
$$
Differentiating under the integral sign we see that
$\square_{2,2} \tilde\psi =0$ when $n=1$.
\qed

\begin{lem}
The limits
\begin{equation}  \label{2limits}
\lim_{\epsilon \to 0^+}
\int_{X \in H_R} \frac{\psi(X)}{N(X-W) + i\epsilon} \,\frac {dS}{\|X\|}
\qquad \text{and} \qquad
\lim_{\epsilon \to 0^+}
\int_{X \in H_R} \frac{\psi(X)}{N(X-W) - i\epsilon} \,\frac {dS}{\|X\|}
\end{equation}
determine continuous functions of $W$ on $\HR^+$ (even if $N(W)=R^2$).
\end{lem}

\pf
Note that when $N(W)=R^2$, i.e. $W \in H_R$, the above argument fails
because $X_{crit}^{\pm} = \pm W$ and $N(X_{crit}^+ - W)=0$.
However, if we remove a small neighborhood $V$ of $W$ in $H_R$, the above
argument still shows that the limits
$$
\lim_{\epsilon \to 0^+} \int_{X \in H_R \setminus V}
\frac{\psi(X)}{N(X-W) \pm i\epsilon}\,\frac {dS}{\|X\|}
$$
converge. On the other hand, as explained in \cite{GS}, the limits
$$
\lim_{\epsilon \to 0^+}
\int_{X \in V} \frac{\psi(X)}{N(X-W) \pm i\epsilon} \,\frac {dS}{\|X\|}
$$
converge as well (which is not hard to see directly).
Thus we conclude that for $R>0$, $W \in H_R$ and a Schwartz function
$\psi$ on $H_R$ the limits (\ref{2limits})
converge.
\qed

Recall that $\deg$ is the degree operator plus identity:
$\deg = 1 + \sum_{i=0}^3 x^i \frac{\partial}{\partial x^i}$.
For each $R>0$, we define operators $I_R^+$ and $I_R^-$ from
${\cal H}(\HR^+)$ into functions on $\HR^+$
\begin{align*}
(I^{+}_R \phi)(W) &= \lim_{\epsilon \to 0^+} \frac 1{2\pi^2} \int_{X \in H_R}
\frac {(\deg_X \phi)(X)}{N(X-W) +i\epsilon} \,\frac{dS}{\|X\|},  \\
(I^{-}_R \phi)(W) &= \lim_{\epsilon \to 0^+} \frac 1{2\pi^2} \int_{X \in H_R}
\frac {(\deg_X \phi)(X)}{N(X-W) -i\epsilon} \,\frac{dS}{\|X\|},
\qquad \phi \in {\cal H}(\HR^+), \quad W \in \HR^+.
\end{align*}
The functions $(I^{\pm}_R \phi)(W)$ are real-analytic and harmonic for
$N(W) \ne R^2$, continuous when $N(W)=R^2$.
By analogy with Theorem 34 in \cite{FL} and similar results of this article,
such as Theorems \ref{discrete_ser_proj} and \ref{discr_proj-harm}, we call
these integrals ``regularized Poisson integrals''.
However, this name is somewhat misleading since these integrals
do not always reproduce functions from ${\cal H}(\HR^+)$.
We finish this subsection by proving that these regularized Poisson integrals
are equivariant with respect to the $\mathfrak{g}(H_R)$-action.


\begin{prop}
The operators $I_R^+$ and $I_R^-$ are equivariant with respect to the
$\pi_l^0$-action of $\mathfrak{g}(H_R)$, where
$\mathfrak{g}(H_R) \simeq \mathfrak{so}(3,2)$ is the Lie algebra
introduced in Corollary \ref{G(H_R)}.
\end{prop}

\pf
The operators $I^{\pm}_R$ are equivariant with respect to the $\pi_l^0$ action of
$SU(1,1) \times SU(1,1)$. Thus it is sufficient to show that $I^{\pm}_R$ are
equivariant with respect to the one-dimensional algebra $\mathfrak{g}(H_R)'$
introduced in the proof of Proposition \ref{invar-pairing-prop}.
Let $g=\exp\begin{pmatrix} 0 & Rt \\ R^{-1}t & 0 \end{pmatrix} \in G(H_R)'$.
We want to compare $(I^{\pm}_R \tilde\phi)(W)$ with
$(\widetilde {I^{\pm}_R \phi})(W)$, where $\tilde\phi = \pi^0_l(g) \phi$
and $(\widetilde {I^{\pm}_R \phi}) = \pi^0_l(g) (I^{\pm}_R (\phi))$.
For $t \to 0$ and modulo terms of order $t^2$, we have (\ref{X-tilde-1}),
(\ref{X-tilde-2}) and (\ref{X-tilde-3}).
Using Lemma \ref{deg-x-tilde}, we get
\begin{multline*}
2\pi^2 \cdot (I^{\pm}_R \tilde\phi)(W)
= \lim_{\epsilon \to 0^+} \int_{X \in H_R}
\frac 1{N(X-W) \pm i\epsilon} \cdot \deg_X \biggl(
\frac{\phi(\tilde X)}{N(-R^{-1} \sinh tX + \cosh t)} \biggr)
\,\frac{dS}{\|X\|}  \\
= \lim_{\epsilon \to 0^+} \int_{X \in H_R}
\frac {(1 + 4t \re X /R) \cdot (\deg_X \phi)(\tilde X)}{N(X-W) \pm i\epsilon}
\,\frac{dS}{\|X\|}.
\end{multline*}
On the other hand, using Lemma 10 from \cite{FL} and Lemma \ref{Jacobian-lemma},
\begin{multline*}
2\pi^2 \cdot (\widetilde {I^{\pm}_R \phi})(W)
= \lim_{\epsilon \to 0^+} \frac 1{N(-R^{-1}\sinh tW + \cosh t)}
\cdot \int_{\tilde X \in H_R} \frac {(\deg_X \phi)(\tilde X)}
{N(\tilde X - \tilde W) \pm i\epsilon} \,\frac{dS}{\|\tilde X\|} \\
= \lim_{\epsilon \to 0^+} \int_{X \in H_R}
\frac {(1 + 4t \re X /R) \cdot (\deg_X \phi)(\tilde X)}
{N(X-W) \pm i\epsilon (1 - 2t \re X /R)(1 - 2t \re W/R)} \,\frac{dS}{\|X\|}  \\
= \lim_{\epsilon \to 0^+} \int_{X \in H_R}
\frac {(1 + 4t \re X /R) \cdot (\deg_X \phi)(\tilde X)}
{N(X-W) \pm i\epsilon} \,\frac{dS}{\|X\|}.
\end{multline*}
This proves that the maps $I^{\pm}_R$ are $\mathfrak{g}(H_R)$-equivariant.
\qed

\subsection{Regularized Integrals and the Discrete Series Component on $\HR$}

In this subsection we study the effect of the operators $I^+_R$ and $I^-_R$
on the discrete series component of ${\cal H}(\HR^+)$.
We prove an analogue of Poisson formula for the discrete series.

\begin{thm}  \label{Poisson}
The operators $I^+_R$ and $I^-_R$ coincide on the discrete series component
of ${\cal H}(\HR^+)$. Moreover,
\begin{align*}
I^{\pm}_R: \qquad
N(Z)^{-2l-1} \cdot t^l_{n\,\underline{m}}(Z) &\mapsto
\begin{cases}
- N(W)^{-2l-1} \cdot t^l_{n\,\underline{m}}(W) & \text{if $N(W) \le R^2$;} \\
- R^{-2(2l+1)} \cdot t^l_{n\,\underline{m}}(W) & \text{if $N(W) \ge R^2$;}
\end{cases}  \\
t^l_{n\,\underline{m}}(Z) &\mapsto
\begin{cases}
R^{2(2l+1)} \cdot N(W)^{-2l-1} \cdot t^l_{n\,\underline{m}}(W) &
\text{if $N(W) \le R^2$;}\\
t^l_{n\,\underline{m}}(W) & \text{if $N(W) \ge R^2$.}
\end{cases}
\end{align*}
(Compare with Theorem \ref{discrete_ser_proj}.)
\end{thm}

\pf
The proof is based on two lemmas.

\begin{lem}  \label{proportionality}
Let $\phi$ be a discrete series matrix coefficient
$$
t^l_{l\,\underline{l}}(Z), \qquad N(Z)^{-2l-1} \cdot t^l_{l\,\underline{l}}(Z), \qquad
t^l_{-l\,\underline{-l}}(Z), \qquad N(Z)^{-2l-1} \cdot t^l_{-l\,\underline{-l}}(Z), \qquad
l= -1, -\frac32, -2 \dots.
$$
Then, for any $R'>0$, $R'\ne R$, the functions on $H_{R'}$
$$
(I^+_R \phi) \bigr|_{H_{R'}}, \qquad (I^-_R \phi) \bigr|_{H_{R'}}
\qquad \text{and} \qquad \phi \bigr|_{H_{R'}}
\qquad \text{are proportional}.
$$
\end{lem}

\pf
The maps $\phi \mapsto (I^{\pm}_R \phi) \bigr|_{H_{R'}}$ are
$\bigl( SU(1,1) \times SU(1,1) \bigr)$-equivariant.
Let $k_{\phi}=
\begin{pmatrix} e^{i\frac{\phi}2} & 0 \\ 0 & e^{-i\frac{\phi}2} \end{pmatrix}$,
then by (\ref{t-german}) we have
$$
t^l_{n\,\underline{m}}(k_{\phi} \cdot Z \cdot k_{\psi}) =
e^{-i(n\phi+m\psi)} \cdot t^l_{n\,\underline{m}}(Z),
\qquad N(k_{\phi} \cdot Z \cdot k_{\psi}) = N(Z).
$$
Let
$$
H = \begin{pmatrix} 1 & 0 \\ 0 & -1 \end{pmatrix}, \qquad
E = \begin{pmatrix} 0 & 1 \\ 0 & 0 \end{pmatrix}, \qquad
F = \begin{pmatrix} 0 & 0 \\ 1 & 0 \end{pmatrix}
$$
be the standard generators of $\mathfrak{sl}(2, \BB C)$.
From (\ref{t=0}) and the Lie algebra action formulas given in the proof of
Theorem \ref{D-component} we see that
$$
\pi_l^0 \begin{pmatrix} E & 0 \\ 0 & 0 \end{pmatrix}
= \pi_r^0 \begin{pmatrix} E & 0 \\ 0 & 0 \end{pmatrix}, \quad
\pi_l^0 \begin{pmatrix} 0 & 0 \\ 0 & F \end{pmatrix}
= \pi_r^0 \begin{pmatrix} 0 & 0 \\ 0 & F \end{pmatrix}
\quad \text{annihilate} \quad
\begin{matrix} t^l_{l\,\underline{l}}(Z) \quad \text{and} \\
N(Z)^{-2l-1} \cdot t^l_{l\,\underline{l}}(Z);\end{matrix}
$$
$$
\pi_l^0 \begin{pmatrix} F & 0 \\ 0 & 0 \end{pmatrix}
= \pi_r^0 \begin{pmatrix} F & 0 \\ 0 & 0 \end{pmatrix}, \quad
\pi_l^0 \begin{pmatrix} 0 & 0 \\ 0 & E \end{pmatrix}
= \pi_r^0 \begin{pmatrix} 0 & 0 \\ 0 & E \end{pmatrix}
\quad \text{annihilate} \quad
\begin{matrix} t^l_{-l\,\underline{-l}}(Z) \quad \text{and} \\
N(Z)^{-2l-1} \cdot t^l_{-l\,\underline{-l}}(Z).\end{matrix}
$$
At a generic point $g \in SU(1,1)$, vector fields
$$
\pi_l^0 \begin{pmatrix} H & 0 \\ 0 & 0 \end{pmatrix}, \qquad
\pi_l^0 \begin{pmatrix} 0 & 0 \\ 0 & H \end{pmatrix}, \qquad
\pi_l^0 \begin{pmatrix} E & 0 \\ 0 & 0 \end{pmatrix}, \qquad
\pi_l^0 \begin{pmatrix} 0 & 0 \\ 0 & F \end{pmatrix}
$$
span the entire complexified tangent space $T^{\BB C}_g SU(1,1)$ at $g$.
Since $(I^{\pm}_R \phi) \bigr|_{H_{R'}}$ with $\phi= t^l_{l\,\underline{l}}(Z)$
or $N(Z)^{-2l-1} \cdot t^l_{l\,\underline{l}}(Z)$ are analytic, it follows from the
$\bigl( SU(1,1) \times SU(1,1) \bigr)$-equivariance that
$(I^{\pm}_R \phi) \bigr|_{H_{R'}}$ are proportional to $\phi\bigr|_{H_{R'}}$.
Similarly, $(I^{\pm}_R \phi) \bigr|_{H_{R'}}$ are proportional to
$\phi\bigr|_{H_{R'}}$ when $\phi = t^l_{-l\,\underline{-l}}(Z)$ or
$N(Z)^{-2l-1} \cdot t^l_{-l\,\underline{-l}}(Z)$.
\qed

\begin{lem} \label{computation}
Let $\phi(Z)$ be $t^l_{l\,\underline{l}}(Z) = (z_{22})^{2l}$ or
$t^l_{-l\,\underline{-l}}(Z) = (z_{11})^{2l}$, and let
$W= \begin{pmatrix} \lambda & 0 \\ 0 & \lambda \end{pmatrix}$ be diagonal
with $\lambda >0$.
Then
$$
(I^{\pm}_R \, t^l_{l\,\underline{l}})(W) = (I^{\pm}_R \, t^l_{-l\,\underline{-l}})(W) =
\begin{cases}
R^{2(2l+1)} \cdot \lambda^{-2(l+1)} & \text{if $\lambda < R$;}\\
\lambda^{2l} & \text{if $\lambda > R$.}
\end{cases}
$$
\end{lem}

The proof is done by direct computation, and we postpone it until next
subsection. Together Lemmas \ref{proportionality} and \ref{computation}
imply that the theorem holds for functions
$$
t^l_{l\,\underline{l}}(Z), \qquad t^l_{-l\,\underline{-l}}(Z), \qquad
N(Z)^{-2l-1} \cdot t^l_{l\,\underline{l}}(Z), \qquad
N(Z)^{-2l-1} \cdot t^l_{-l\,\underline{-l}}(Z)
$$
for $l \le -1$.
Then the theorem follows from the
$\bigl( SU(1,1) \times SU(1,1) \bigr)$-equivariance of $I^{\pm}_R$.
\qed

\subsection{Proof of Lemma \ref{computation}}

In this subsection we give a proof of Lemma \ref{computation}.

\pf
We use integral formulas 2.117 and 2.662 from \cite{GR}
$$
\int \frac {dx}{x^n(a+bx)} =
\sum_{k=1}^{n-1} \frac {(-1)^k b^{k-1}}{(n-k)a^kx^{n-k}}
+ (-1)^n \frac {b^{n-1}}{a^n} \log \frac{a+bx}{x} \quad +C,
$$
$$
\int \frac {dx}{x^n(a-bx)} =
-\sum_{k=1}^{n-1} \frac {b^{k-1}}{(n-k)a^kx^{n-k}}
- \frac {b^{n-1}}{a^n} \log \frac{bx-a}{x} \quad +C,
$$
$$
\int e^{ax} (\cos x)^{2m} \,dx =
\begin{pmatrix} 2m \\ m \end{pmatrix} \frac{e^{ax}}{2^{2m}a} +
\frac{e^{ax}}{2^{2m-1}} \sum_{k=1}^m \begin{pmatrix} 2m \\ m-k \end{pmatrix}
\frac{a \cos (2kx) +2k \sin (2kx)}{a^2 + 4k^2} \quad +C,
$$
$$
\int e^{ax} (\cos x)^{2m+1} \,dx =
\frac{e^{ax}}{2^{2m}} \sum_{k=0}^m \begin{pmatrix} 2m+1 \\ m-k \end{pmatrix}
\frac{a \cos (2k+1)x +(2k+1) \sin (2k+1)x}{a^2 + (2k+1)^2} \quad +C.
$$
We parameterize $H_R$ as in (\ref{param}) so that
$$
X(\phi,\tau,\psi) = R
\begin{pmatrix} \cosh \frac{\tau}2 \cdot e^{i\frac{\phi+\psi}2} &
\sinh \frac{\tau}2 \cdot e^{i\frac{\phi-\psi}2} \\
\sinh \frac{\tau}2 \cdot e^{i\frac{\psi-\phi}2} &
\cosh \frac{\tau}2 \cdot e^{-i\frac{\phi+\psi}2} \end{pmatrix},
\qquad
\begin{matrix}
0 \le \phi < 2\pi, \\ 0 < \tau < \infty, \\ -2\pi \le \psi < 2\pi.
\end{matrix}
$$
Then $\frac {dS}{\|X\|} = \frac {R^2}8 \sinh\tau \,d\phi d\tau d\psi$ and
$$
N(X-W) = R^2+\lambda^2-2R\lambda\cosh\frac{\tau}2 \cdot \cos\frac{\phi+\psi}2.
$$
Changing variables $\theta =  \frac{\phi+\psi}2$,
$\theta' =  \frac{\phi-\psi}2$, $x= \cosh\frac{\tau}2$ and letting
$a(\epsilon)= R^2+\lambda^2 + i \epsilon$, $b(\theta)=2R\lambda\cos\theta$,
we obtain
\begin{multline*}
(I^{\pm}_R \, t^l_{l\,\underline{l}})(W) = (I^{\pm}_R \, t^l_{-l\,\underline{-l}})(W)  \\
=
\frac{2l+1}{16\pi^2} \lim_{\epsilon \to 0^{\pm}}
\int_{\psi=-2\pi}^{2\pi} \int_{\tau=0}^{\infty} \int_{\phi=0}^{2\pi}
\frac {R^{2l+2} \cdot (\cosh\frac{\tau}2)^{2l} \cdot e^{\pm il(\phi+\psi)} \cdot
\sinh \tau} {N(W-X)+i\epsilon} \,d\phi d\tau d\psi\\
=
\frac{2l+1}{4\pi^2} \lim_{\epsilon \to 0^{\pm}} \int_{\theta=0}^{2\pi}
\int_{\tau=0}^{\infty} \int_{\theta'=0}^{2\pi} \frac {R^{2l+2} \cdot
(\cosh\frac{\tau}2)^{2l+1} \cdot e^{\pm i2l\theta} \cdot \sinh\frac{\tau}2}
{R^2+\lambda^2 - 2R\lambda \cos\theta \cdot \cosh\frac{\tau}2 +i\epsilon}
\,d\theta' d\tau d\theta\\
=
\frac{2l+1}{\pi} \lim_{\epsilon \to 0^{\pm}} \int_{\theta=-\pi/2}^{\pi/2}
R^{2l+2} \cdot e^{\pm i2l\theta} \int_{x=1}^{\infty}
\biggl( \frac{x^{2l+1}}{a(\epsilon) - b(\theta)x}
+ \frac{(-1)^{2l} x^{2l+1}}{a(\epsilon) + b(\theta)x} \biggr)\,dx d\theta\\
=
\frac{2l+1}{\pi} \lim_{\epsilon \to 0^{\pm}} \int_{\theta=-\pi/2}^{\pi/2} R^{2l+2} \cdot
e^{\pm i2l\theta} \biggl( \frac{a(\epsilon)^{2l+1}}{b(\theta)^{2l+2}} \Bigl(
\log \Bigl( \frac{a(\epsilon)-b(\theta)}{b(\theta)}\Bigr)
+ \log \Bigl( \frac{a(\epsilon)+b(\theta)}{b(\theta)}\Bigr) \mp \pi i \Bigr) \\
\hspace{2in}
- \sum_{k=1}^{-2l-2} \frac {\bigl( 1-(-1)^{2l+k} \bigr) \cdot b(\theta)^{k-1}}
{(2l+k+1)a(\epsilon)^k} \biggr)\,d\theta\\
=
\frac{2l+1}{4\pi} \cdot \frac{(R^2+\lambda^2)^{2l+1}}{2^{2l} \lambda^{2l+2}}
\lim_{\epsilon \to 0^{\pm}} \int_{\theta=-\pi/2}^{\pi/2}
\frac{e^{\pm i2l\theta}}{(\cos \theta)^{2l+2}} \cdot \biggl(
\log \Bigl( \frac{a(\epsilon)-b(\theta)}{b(\theta)}\Bigr)
+ \log \Bigl( \frac{a(\epsilon)+b(\theta)}{b(\theta)}\Bigr)
\biggr)\,d\theta.
\end{multline*}
From now on we assume $2l$ is odd; the other case is similar.
Integrating by parts we get
$$
(I^{\pm}_R \, t^l_{l\,\underline{l}})(W) = (I^{\pm}_R \, t^l_{-l\,\underline{-l}})(W)
= \frac{2l+1}{2\pi} \cdot \frac{(R^2+\lambda^2)^{2l+1}}{\lambda^{2l+2}}
\sum_{k=0}^{-l-\frac32} \begin{pmatrix} -2l-2 \\ -l-k- \frac32 \end{pmatrix}
\frac{I(l,k)}{4l^2 - (2k+1)^2},
$$
where integrals $I(l,k)$ can be computed using residues:
\begin{multline*}
I(l,k) = \int_{\theta=-\pi/2}^{\pi/2} e^{\pm i2l\theta}
\bigl((2k+1) \sin(2k+1)\theta \pm i2l \cos(2k+1)\theta \bigr) \cdot \biggl(
\log \Bigl( \frac{(R^2+\lambda^2)^2-b^2(\theta)}{b^2(\theta)}\Bigr)\biggr)'
\,d\theta\\
=
2(R^2+\lambda^2)^2 \int_{\theta=-\pi/2}^{\pi/2} e^{\pm i2l\theta} \cdot \tan\theta \cdot
\frac{(2k+1) \sin (2k+1)\theta \pm i2l \cos (2k+1)\theta}
{(R^2+\lambda^2)^2-4R^2\lambda^2\cos^2\theta} \,d\theta\\
=
(R^2+\lambda^2)^2 \int_{\theta=-\pi}^{\pi} \tan\theta \cdot
\frac{(2k+1) \cos(2l\theta) \cdot \sin (2k+1)\theta
- 2l \sin(2l\theta) \cdot \cos (2k+1)\theta}
{(R^2+\lambda^2)^2-4R^2\lambda^2\cos^2\theta} \,d\theta = \\
\frac{i (R^2+\lambda^2)^2}4 \oint \frac{s^2-1}{s^2+1} \cdot
\frac{(2k+1)(s^{2l}+s^{-2l})(s^{2k+1}-s^{-2k-1}) - 2l(s^{2l}-s^{-2l})(s^{2k+1}+s^{-2k-1})}
{(R^2+\lambda^2)^2-R^2\lambda^2 (s+s^{-1})^2} \,\frac{ds}s.
\end{multline*}
From the factorization
$$
(R^2+\lambda^2)^2-R^2\lambda^2 (s+s^{-1})^2 = -\frac {R^2\lambda^2}{s^2}
(s-R/\lambda) (s+R/\lambda) (s-\lambda/R) (s+\lambda/R)
$$
and decomposition
$$
\frac{(R^2+\lambda^2)^2 \cdot (s^2-1)}{(R^2+\lambda^2)^2-R^2\lambda^2 (s+s^{-1})^2}
=
\frac{R^2}{R^2-\lambda^2 s^2} + \frac{\lambda^2}{\lambda^2-R^2s^2}
- \frac2{1+s^2}
$$
we see that the residues are
\begin{multline*}
\frac 1{2i} \Bigl( (2k+1-2l) \bigl(
(\lambda/R)^{2l+2k+1} - (R/\lambda)^{2l+2k+1} \bigr)  \\
+ (2l+2k+1) \bigl( (\lambda/R)^{2k+1-2l} - (R/\lambda)^{2k+1-2l} \bigr)
\Bigr) \cdot \begin{cases}
+1 & \text{at $\pm \lambda/R$,}\\
-1 & \text{at $\pm R/\lambda$}
\end{cases}
\end{multline*}
and
\begin{multline*}
i(2k+1-2l) \bigl(
(\lambda/R)^{2l+2k+1} + (R/\lambda)^{2l+2k+1} -2 \cdot i^{2l+2k+1} \bigr)  \\
- i(2l+2k+1) \bigl( (\lambda/R)^{2l-2k-1}
+ (R/\lambda)^{2l-2k-1} -2 \cdot i^{2l-2k-1} \bigr)
\qquad \text{at $0$}.
\end{multline*}
Letting
$r= \begin{cases} \lambda/R & \text{if $\lambda<R$,}\\
R/\lambda & \text{if $\lambda>R$,} \end{cases}$ we get:
\begin{multline*}
(I^{\pm}_R \, t^l_{l\,\underline{l}})(W) = (I^{\pm}_R \, t^l_{-l\,\underline{-l}})(W) = \\
-2(2l+1) \cdot \frac{(R^2+\lambda^2)^{2l+1}}{\lambda^{2l+2}}
\sum_{k=0}^{-l-\frac32} \begin{pmatrix} -2l-2 \\ -l-k- \frac32 \end{pmatrix}
\biggl( \frac{ r^{-2l-2k-1} - i^{-2l-2k-1}}{-2l-2k-1}
+ \frac{r^{-2l+2k+1} - i^{-2l+2k+1}}{-2l+2k+1} \biggr) \\
=
-2(2l+1) \cdot \frac{(R^2+\lambda^2)^{2l+1}}{\lambda^{2l+2}}
\sum_{k=0}^{-l-\frac32} \begin{pmatrix} -2l-2 \\ -l-k- \frac32 \end{pmatrix}
\int_i^r \bigl( x^{-2l-2k-2} + x^{-2l+2k} \bigr) \,dx \\
=
-\frac{(R^2+\lambda^2)^{2l+1}}{\lambda^{2l+2}}
\int_i^r 2(2l+1)x(1+x^2)^{-2l-2} \,dx  \\
=
\frac{(R^2+\lambda^2)^{2l+1}}{\lambda^{2l+2}} \cdot \bigl(1+r^2\bigr)^{-2l-1}
=
\begin{cases}
R^{2(2l+1)} \cdot \lambda^{-2(l+1)} & \text{if $\lambda <R$;}\\
\lambda^{2l} & \text{if $\lambda >R$.}
\end{cases}
\end{multline*}
\qed

\subsection{Regularized Integrals and the Continuous Series Component on $\HR$}

In this subsection we determine the actions of the operators $I^{\pm}_R$
on the continuous series component of ${\cal H}(\HR^+)$.
Combining this with Theorem \ref{Poisson} we get a complete
description of the actions of $I^{\pm}_R$ on ${\cal H}(\HR^+)$.
We state the main result of this section:

\begin{thm}  \label{cont_ser_proj}
The operator $I^+_R - I^-_R$ annihilates the discrete series component of
${\cal H}(\HR^+)$. Its action on the continuous series component of
${\cal H}(\HR^+)$ is given by
$$
\bigl( (I^+_R - I^-_R)t^l_{n\,\underline{m}} \bigr)(W) =
\bigl(1 + R^{4 \im l} \cdot N(W)^{-2l-1} \bigr) \cdot t^l_{n\,\underline{m}}(W) \cdot
\begin{cases}
\coth(\pi\im l) & \text{if $m,n \in \BB Z$;}  \\
\tanh(\pi\im l) & \text{if $m,n \in \BB Z +\frac12$,}
\end{cases}
$$
where $\re l =-\frac12$ and $W \in \HR^+$.
\end{thm}

The proof of this theorem will be given in the next subsection.

\begin{rem}
Strictly speaking, the matrix coefficient functions of the continuous series
representations $t^l_{n\,\underline{m}}(X)$, $l=-1/2+i\lambda$ with
$\lambda \in \BB R$, do not belong to ${\cal H}(\HR^+)$.
So the theorem should be restated as follows: If
$$
\phi(X) = \int_{\lambda \in \BB R} t^{-\frac12 + i\lambda}_{n\,\underline{m}}(X) \cdot
\eta(\lambda) \,d\lambda,
$$
where $\eta(\lambda)$ is a smooth compactly supported function on
$\BB R \setminus \{0\}$, then $\bigl( (I^+_R - I^-_R)\phi \bigr)(W)$ is
$$
\begin{cases}
\int_{\lambda \in \BB R} \coth(\pi\lambda) \cdot
\bigl(1 + R^{4 \lambda} \cdot N(W)^{-2\lambda} \bigr) \cdot
t^{-\frac12 + i\lambda}_{n\,\underline{m}}(W) \cdot \eta(\lambda) \,d\lambda
& \text{if $m,n \in \BB Z$;}  \\
\int_{\lambda \in \BB R} \tanh(\pi\lambda) \cdot
\bigl(1 + R^{4 \lambda} \cdot N(W)^{-2\lambda} \bigr) \cdot
t^{-\frac12 + i\lambda}_{n\,\underline{m}}(W) \cdot \eta(\lambda) \,d\lambda
& \text{if $m,n \in \BB Z +\frac12$.}
\end{cases}
$$
\end{rem}

Note that the expression
$$
\begin{cases}
\frac{\coth(\pi\im l)}{\im l} & \text{if $m,n \in \BB Z$;}  \\
\frac{\tanh(\pi\im l)}{\im l} & \text{if $m,n \in \BB Z +\frac12$}
\end{cases}
$$
is precisely the inverse of the Plancherel measure of $SU(1,1)$!
Since $\deg t^l_{n\,\underline{m}} = 2\im l \cdot t^l_{n\,\underline{m}}$,
we can reformulate Theorem \ref{cont_ser_proj} as follows:

\begin{thm}  \label{cont_ser_proj-2}
The operator
$$
\phi(X) \mapsto (\operatorname{Pl}_R \phi)(W) =
\lim_{\epsilon \to 0^+} \frac1{\pi^2} \int_{X \in H_R}
\biggl( \frac1{N(X-W) +i\epsilon} - \frac1{N(X-W) -i\epsilon} \biggr)
\phi(X) \,\frac{dS}{\|X\|}
$$
annihilates the discrete series component of ${\cal H}(\HR^+)$.
If $\re l =-\frac12$ and $W \in \HR^+$,
$$
(\operatorname{Pl}_R t^l_{n\,\underline{m}})(W) =
\bigl(1 + R^{4 \im l} \cdot N(W)^{-2l-1} \bigr) \cdot t^l_{n\,\underline{m}}(W) \cdot
\begin{cases}
\frac{\coth(\pi\im l)}{\im l} & \text{if $m,n \in \BB Z$;}  \\
\frac{\tanh(\pi\im l)}{\im l} & \text{if $m,n \in \BB Z +\frac12$.}
\end{cases}
$$
\end{thm}

\begin{rem}
As in the Minkowski case (see Remark \ref{M-loc-rem}), for each $W \in \HR^+$,
integrals $\bigl( (I^+_R - I^-_R)\phi \bigr)(W)$ and
$(\operatorname{Pl}_R \phi)(W)$ depend only on the
values of $(\deg\phi)(X)$, where  $X$ ranges over a two-dimensional hyperboloid
(or cone if $W \in H_R$)
$$
\{X \in \HR;\: N(X)=R^2, \: N(X-W)=0\}.
$$
\end{rem}


We can rewrite the expression $\frac 1{N(X-W)+i0} - \frac 1{N(X-W)-i0}$
geometrically as follows. Let us consider $W'=(1+\epsilon\xi)W \in \HC$,
$\epsilon>0$, $\xi \in \BB C \setminus \BB R$.
Then, for $X \in SU(1,1)$, $W \in \HR^+$,
$$
N(X-W') = N(W) \cdot N\bigl( X \cdot W^{-1}-(1+\epsilon\xi) \bigr) \ne 0
$$
for all $\epsilon \ll 1$ because $X \cdot W^{-1} \in N(W)^{-1/2} \cdot SU(1,1)$
cannot have an eigenvalue $1+\epsilon\xi$.
When $N(X-W)=0$,
$N(X-W')=\epsilon\xi \bigl( N(W)- N(X) \bigr) + \epsilon^2 \xi^2 N(W)$.
Choose $\xi_1, \xi_2 \in \BB C$ with $\im(\xi_1)>0$ and $\im(\xi_2)<0$, then
\begin{multline*}
\lim_{\epsilon \to 0^+} \sgn\bigl( N(W)-1 \bigr) \cdot \int_{X \in SU(1,1)}
\biggl( \frac1{N(X-W) +i\epsilon} - \frac1{N(X-W) -i\epsilon} \biggr)
\phi(X) \,\frac{dS}{\|X\|} \\
 = \lim_{\epsilon \to 0^+} \int_{X \in SU(1,1)}
\biggl( \frac1{N\bigl( X-(1+\epsilon\xi_1)W \bigr)} -
\frac1{N\bigl( X-(1+\epsilon\xi_2)W \bigr)} \biggr) \phi(X) \,\frac{dS}{\|X\|}.
\end{multline*}

\begin{prop}
The operator $I^-_R + I^+_R$ annihilates the continuous series component of
${\cal H}(\HR^+)$.
\end{prop}

\pf
Let $\phi$ be in the continuous series component of ${\cal H}(\HR^+)$.
Since $(I^-_R + I^+_R)\phi$ is analytic, it is enough to show that
$\deg^k \bigl( (I^-_R + I^+_R)\phi \bigr)$ vanishes on $H_R$ for $k=0,1,2,\dots$.
By Theorem \ref{discr_proj-harm},
$(\S_{R,\sigma^{-1}}^-\phi)(W)$ and $(\S_{R,\sigma}^+\phi)(W)$
vanish identically for $W$ in the open set
$$
R \cdot \begin{pmatrix} \sigma^{-1} & 0 \\ 0 & \sigma \end{pmatrix} \cdot
\Gamma^- \cap R \cdot \begin{pmatrix} \sigma & 0 \\ 0 & \sigma^{-1} \end{pmatrix}
\cdot \Gamma^+, \qquad \forall \sigma >1.
$$
For $W \in H_R$ we have:
\begin{multline*}
- \lim_{\sigma \to 1^+} \bigl( \S_{R,\sigma^{-1}}^-\phi + \S_{R,\sigma}^+\phi \bigr)(W) \\
= \lim_{\epsilon \to 0^+} \frac1{2\pi^2}\int_{X \in H_R}
\biggl( \frac {(\deg \phi)(X)}{N(X-W) +i\epsilon}
+ \frac {(\deg \phi)(X)}{N(X-W) -i\epsilon} \biggr)\,\frac{dS}{\|X\|}
= (I^-_R\phi + I^+_R\phi)(W).
\end{multline*}
Hence $(I^-_R + I^+_R)\phi$ vanishes on $H_R$ along with all derivatives
$\deg^k \bigl( (I^-_R + I^+_R)\phi \bigr)$, $k=0,1,2,\dots$.
\qed

This proposition combined with Theorems \ref{Poisson} and \ref{cont_ser_proj}
gives us a complete description of the actions of $I^-_R$ and $I^+_R$ on the
discrete and continuous series components of ${\cal H}(\HR^+)$.

We conclude this subsection with a result for left-regular functions.
(Of course, similar result holds for right-regular functions.)

\begin{thm}
The operator sending left-regular functions $f(X) \in \BB S(\HR^+)$ into
$$
f(X) \quad \mapsto \quad \lim_{\epsilon \to 0^+} \frac1{2\pi^2} \int_{X \in H_R}
\biggl( \frac{(X-W)^+}{(N(X-W) +i\epsilon)^2}
- \frac{(X-W)^+}{(N(X-W) -i\epsilon)^2} \biggr) \cdot Dx \cdot f(X)
$$
annihilates the discrete series component. If $\re l =-\frac12$,
$$
\begin{pmatrix} (l-m+ \frac 12) t^l_{n \, \underline{m+ \frac 12}}(X)  \\
(l+m+ \frac 12) t^l_{n \, \underline{m- \frac 12}}(X) \end{pmatrix} \quad \mapsto \quad
\begin{pmatrix} (l-m+ \frac 12) t^l_{n \, \underline{m+ \frac 12}}(W)  \\
(l+m+ \frac 12) t^l_{n \, \underline{m- \frac 12}}(W) \end{pmatrix} \cdot
\begin{cases}
\coth(\pi\im l) & \text{if $m,n \in \BB Z$;}  \\
\tanh(\pi\im l) & \text{if $m,n \in \BB Z +\frac12$.}
\end{cases}
$$
\end{thm}

This theorem follows from Theorem \ref{cont_ser_proj} in exactly the same way
Theorem \ref{M-discr-proj-reg} follows from Theorem \ref{M-discrete-proj}
in the Minkowski case.

\subsection{Proof of Theorem \ref{cont_ser_proj}}

In this subsection we give a proof of Theorem \ref{cont_ser_proj}.

\pf
The operators $I^{\pm}_R$ are equivariant with respect to the subgroup
$SU(1,1) \times SU(1,1)$ of $SL(2,\HC)$.
Also, by Theorem \ref{Poisson} the operators $I^{\pm}_R$ coincide on the
discrete series component of ${\cal H}(\HR^+)$.
Therefore, the operator $I^+_R - I^-_R$ annihilates the discrete series
and sends the continuous series component into itself.
To determine the action of $I^+_R - I^-_R$ on the continuous series component
we use a couple of reductions.
First, it is sufficient to prove the theorem for $R=1$ only.
Indeed, substituting $X' = R^{-1} \cdot X \in SU(1,1)$ and
$W' =  R^{-1} \cdot W$ we get
$$
\bigl( I^{\pm}_R t^l_{n\,\underline{m}} \bigr)(W)
=
\lim_{\epsilon \to 0^+} \frac{R^{2l}}{2\pi^2} \int_{X' \in H_1}
\frac{(\deg_{X'} t^l_{n\,\underline{m}})(X')}{N(X'-W') \pm i\epsilon}
\,\frac{dS}{\|X'\|}
=
R^{2l} \cdot \bigl( I^{\pm}_1 t^l_{n\,\underline{m}} \bigr)(R^{-1} \cdot W).
$$
Secondly, we can reduce the proof of the theorem to the case $N(W)=1$:

\begin{prop}  \label{cont_ser_proj_1}
When $N(W)=1$, and $\re l =-\frac12$
$$
\bigl( (I^+_1 - I^-_1)t^l_{n\,\underline{m}} \bigr)(W) =
2 t^l_{n\,\underline{m}}(W) \cdot
\begin{cases}
\coth(\pi\im l) & \text{if $m,n \in \BB Z$;}  \\
\tanh(\pi\im l) & \text{if $m,n \in \BB Z +\frac12$.}
\end{cases}
$$
\end{prop}

We postpone the proof of this proposition until the end of this subsection.
By Lemma 10 in \cite{FL} and direct computation we have:
$$
\pi^0_l(\gamma)_X \circ \pi^0_l(\gamma)_W
\biggl( \frac 1{N(X-W) \pm i\epsilon} \biggr) =
\frac1{N(X'-Y) \pm \frac {i\epsilon}4 \cdot
N(\tilde e_3X'-1) \cdot N(\tilde e_3Y-1)},
$$
where $X'= \pi_l(\gamma^{-1})(X)$ and $Y= \pi_l(\gamma^{-1})(W)$.
When $N(X)=N(W)=1$,
$$
N(\tilde e_3X' -1) \cdot N(\tilde e_3Y-1)
= \tr(\tilde e_3 X') \cdot \tr(\tilde e_3 Y) = -4 \re(e_3 X') \cdot \re(e_3 Y).
$$
Therefore, by Theorem \ref{inner_product_pullback},
\begin{multline}  \label{H-M-switch}
\sgnx(Y) \cdot \pi^0_l(\gamma) \bigl( (I^+_1 - I^-_1) \phi \bigr) (Y) \\
=
\lim_{\epsilon \to 0^+} \frac i{2\pi^2} \int_{X' \in \tilde H}
\biggl( \frac 1{N(X'-Y) +i\epsilon} - \frac 1{N(X'-Y) -i\epsilon} \biggr)
\cdot \bigl(\deg_{X'} (\pi_l^0(\gamma)\phi) \bigr)(X') \,\frac{dS}{\|X'\|}  \\
=
- \S_1^{\BB M} (\pi_l^0(\gamma)\phi)(Y).
\end{multline}
Then by Theorem \ref{M-discrete-proj} and Lemma \ref{degrees} we have
$\deg \bigl( (I^+_1 - I^-_1) \phi \bigr)(W)=0$ whenever $W \in SU(1,1)$.
Now Theorem \ref{cont_ser_proj} follows from Proposition \ref{cont_ser_proj_1}
and the uniqueness part of the theorem given on page 343 in \cite{St}.
That theorem essentially states that the Cauchy problem for
$\square_{2,2}\phi=0$ with boundary data on $SU(1,1)$ in the continuous series
component of $L^2(SU(1,1))$ has a unique analytic solution.
(Caution: there can be other non-analytic solutions.)
\qed

The rest of this subsection will be devoted to the proof of
Proposition \ref{cont_ser_proj_1}.

\pf
We start with (\ref{H-M-switch}) and track down carefully which functions
in the continuous part of ${\cal H}(\HR^+)$ get sent where in
${\cal H}(\BB M^+)$. The operator $\phi \mapsto (I^+_1 - I^-_1) \phi$ is
$SU(1,1) \times SU(1,1)$-equivariant.
This means that if $\phi$ is a continuous series function, then
$(I^+_1 - I^-_1) \phi \bigr|_{SU(1,1)}$ is proportional to $\phi$ with a
coefficient of proportionality possibly depending on $l$ and whether $m$, $n$
are integers or half-integers. We essentially compare the expansions of
$\deg_Y \bigl( \pi_l^0(\gamma)\phi \bigr)(Y)$ and
$\sgnx(Y) \cdot \bigl( \pi_l^0(\gamma)\phi \bigr)(Y)$
in terms of basis functions (\ref{basis1}) or (\ref{basis2}).
In order to find this coefficient of proportionality we pick a basis function
$f^{\pm}_{l,m,n}$ and find the ratio between the coefficients in those expansions.

Let
$Y = \begin{pmatrix} y_{11} & y_{12} \\ y_{21} & y_{22} \end{pmatrix}
\in \tilde H \subset \BB M^+$
and
$W = \begin{pmatrix} w_{11} & w_{12} \\ w_{21} & w_{22} \end{pmatrix}
= \pi_l(\gamma)(Y) \in SU(1,1)$. By (\ref{Y-tilde}),
$$
W = \frac{-i}{y_{22}-y_{11}}
\begin{pmatrix} 2+y_{11}+y_{22} & 2y_{12} \\ -2y_{21} & 2-y_{11}-y_{22}\end{pmatrix}.
$$
From (\ref{t-german})
$$
t^l_{n\,\underline{m}}(\psi_1, \tau, \psi_2)
= e^{-i(n+m) \frac{\psi_1 + \psi_2}2} \cdot e^{-i(n-m) \frac{\psi_1 - \psi_2}2}
\cdot \mathfrak{P}^l_{n\,m}(\cosh \tau).
$$
Since
$$
\cosh^2 \frac{\tau}2 = x_{11}x_{22} = -\frac {4-(y_{11}+y_{22})^2}{(y_{22}-y_{11})^2}
= \frac {1+t^2}{(x^3)^2} = \frac 1{\cos^2 \theta},
\qquad
\cosh \tau = \frac 2{\cos^2\theta} -1,
$$
$$
e^{i(\psi_1 + \psi_2)} = \frac{x_{11}}{x_{22}}
= \frac{2+y_{11}+y_{22}}{2-(y_{11}+y_{22})} = \frac{1-it}{1+it}
= \frac{1-i\sinh\rho}{1+i\sinh\rho},
\qquad
e^{-i\frac{\psi_1+\psi_2}2} = \pm \sqrt{\frac{1+i\sinh\rho}{1-i\sinh\rho}}
$$
(which is ambiguous), and
$$
e^{i\frac{\psi_1-\psi_2}2} = \frac{x_{12}}{\sinh \frac{\tau}2}
= - \frac {2i}{\sinh \frac{\tau}2} \frac{y_{12}}{y_{22}-y_{11}}
= i \frac{\tan\theta}{\sinh \frac{\tau}2} e^{-i\phi}
= i \sgn(\cos\theta) \cdot e^{-i\phi},
$$
we obtain
$$
t^l_{n\,\underline{m}}(Y) = \bigl(i \sgn(\cos\theta) \bigr)^{m-n} \cdot
\mathfrak{P}^l_{n\,m}\Bigl(\frac 2{\cos^2\theta} -1\Bigr) \cdot
\biggl( \pm \sqrt{\frac{1+i\sinh\rho}{1-i\sinh\rho}} \biggr)^{n+m}
\cdot e^{i(n-m)\phi}.
$$
To avoid the sign ambiguity we simply set $m+n=0$.
Then, by Lemma \ref{degrees},
\begin{align*}
\sgnx(Y) \cdot \bigl(\pi_l^0(\gamma) t^l_{-m\,\underline{m}}\bigr)(Y)
&= i \frac{\bigl( i\sgn(\cos\theta) \bigr)^{2m}}
{\cosh\rho \cdot |\cos\theta|} \cdot
\mathfrak{P}^l_{-m\,m}\Bigl(\frac 2{\cos^2\theta} -1\Bigr) \cdot e^{-2im\phi}, \\
\deg_Y \bigl( \pi_l^0(\gamma) t^l_{-m\,\underline{m}} \bigr)(Y)
&= (2l+1) \frac{\bigl( i\sgn(\cos\theta) \bigr)^{2m}}
{\cosh^2\rho \cdot \cos^2\theta} \cdot
\mathfrak{P}^l_{-m\,m}\Bigl(\frac 2{\cos^2\theta} -1\Bigr) \cdot e^{-2im\phi}.
\end{align*}
On the other hand, when $N(Y)=1$, the functions $f^{\pm}_{\lambda,\mu,\nu}(Y)$
are proportional to
$$
r^{-1} \cdot P^{(\nu)}_{\lambda} (t/r) \cdot P_{\lambda}^{(\mu)}(\cos\theta)
\cdot e^{-i\mu\phi}
=
(\cosh\rho)^{-1} \cdot P^{(\nu)}_{\lambda} (\tanh\rho) \cdot
P_{\lambda}^{(\mu)}(\cos\theta) \cdot e^{-i\mu\phi}.
$$
Thus $\mu=2m$ and we need to compare integrals
$$
i \biggl( \int_0^{\pi} \bigl( \sgn(\cos\theta) \bigr)^{2m} \cdot
\mathfrak{P}^l_{-m\,m}\Bigl(\frac 2{\cos^2\theta} -1\Bigr)
\cdot P_{\lambda}^{(2m)}(\cos\theta) \cdot |\tan\theta| \,d\theta \biggr)
\cdot
\biggl( \int_{-\infty}^{\infty} P^{(\nu)}_{\lambda} (\tanh\rho) \,d\rho \biggr)
$$
and
$$
(2l+1)\biggl( \int_0^{\pi} \bigl( \sgn(\cos\theta) \bigr)^{2m} \cdot
\mathfrak{P}^l_{-m\,m}\Bigl(\frac 2{\cos^2\theta} -1\Bigr) \cdot
P_{\lambda}^{(2m)}(\cos\theta) \cdot \frac{\tan\theta}{\cos\theta}\,d\theta \biggr)
\cdot
\biggl( \int_{-\infty}^{\infty} \frac{P^{(\nu)}_{\lambda} (\tanh\rho)}{\cosh\rho}
\,d\rho \biggr).
$$
Substituting $u=\tanh\rho$, we find
$$
\int_{-\infty}^{\infty} P^{(\nu)}_{\lambda} (\tanh\rho) \,d\rho
= \int_{-1}^1 \frac{P^{(\nu)}_{\lambda}(u)}{1-u^2} \,du,
\qquad
\int_{-\infty}^{\infty} \frac{P^{(\nu)}_{\lambda} (\tanh\rho)}{\cosh\rho} \,d\rho
= \int_{-1}^1 \frac{P^{(\nu)}_{\lambda}(u)}{\sqrt{1-u^2}} \,du.
$$
If we take $\nu=\lambda$, $P^{(\nu)}_{\lambda}(u)$ is proportional to
$(1-u^2)^{\lambda/2}$, and
$$
\int_{-1}^1 (1-u^2)^{a-1} \,du = \frac{\sqrt{\pi} \, \Gamma(a)}{\Gamma(a+1/2)}.
$$
Thus,
$$
\frac{\int_{-\infty}^{\infty}\frac{P^{(\lambda)}_{\lambda}(\tanh\rho)}{\cosh\rho}\,d\rho}
{\int_{-\infty}^{\infty}P^{(\lambda)}_{\lambda}(\tanh\rho)\,d\rho}
= \frac{\bigl(\Gamma(\lambda/2+1/2)\bigr)^2}
{\Gamma(\lambda/2) \cdot \Gamma(\lambda/2+1)}.
$$

It remains to compare the integrals over $\theta$.
When $m \le 0$, by (\ref{P-hypergeometric}),
$$
\mathfrak{P}^l_{-m\,m}\Bigl(\frac 2{\cos^2\theta} -1\Bigr)=
\frac{\Gamma(l-m+1)}{\Gamma(l+m+1)(-2m)!} \cdot
|\tan\theta|^{-2m} \cdot
\,_2F_1\bigl( l-m+1,-l-m;-2m+1; -\tan^2\theta \bigr).
$$
If $m$ is a half-integer, we need $P_{\lambda}^{(2m)}(x)$ to be odd,
so we can take $\lambda =1-2m$; then $P_{\lambda}^{(2m)}(\cos\theta)$
is proportional to $\cos\theta \cdot \sin^{-2m}\theta$.
Substituting $v=\tan^2\theta$ and using a special case of an integral
formula 7.512(10) from \cite{GR}
$$
\int_0^{\infty} x^{c-1}(1+x)^{-r} \,_2F_1(a,b;c;-x)\,dx =
\frac{\Gamma(c)\Gamma(a+r-c)\Gamma(b+r-c)}{\Gamma(r)\Gamma(a+b+r-c)}
$$
valid when $\re c >0$, $\re(a+r-c)>0$, $\re(b+r-c)>0$, we obtain
\begin{multline*}
2 \frac{\Gamma(l+m+1)(-2m)!}{\Gamma(l-m+1)}
\int_0^{\pi/2} \mathfrak{P}^l_{-m\,m}\Bigl(\frac 2{\cos^2\theta} -1\Bigr)
\cdot \cos\theta \cdot \sin^{-2m}\theta \cdot \tan\theta \,d\theta \\
= \int_0^{\infty} v^{-2m} \cdot (1+v)^{m-3/2} \cdot
\,_2F_1\bigl( l-m+1,-l-m;-2m+1; -v \bigr)\,dv  \\
= \frac{\Gamma(1-2m)\Gamma(3/2+l)\Gamma(1/2-l)}{\bigl(\Gamma(3/2-m)\bigr)^2}
= \frac{\pi}{2 \cos(\pi l)}
\frac{(2l+1)\Gamma(1-2m)}{\bigl(\Gamma(3/2-m)\bigr)^2},
\end{multline*}
\begin{multline}  \label{coeff-calc}
2 \frac{\Gamma(l+m+1)(-2m)!}{\Gamma(l-m+1)}
\int_0^{\pi/2} \mathfrak{P}^l_{-m\,m}\Bigl(\frac 2{\cos^2\theta} -1\Bigr)
\cdot \sin^{-2m}\theta \cdot \tan\theta \,d\theta \\
= \int_0^{\infty} v^{-2m} \cdot (1+v)^{m-1} \cdot
\,_2F_1\bigl( l-m+1,-l-m;-2m+1; -v \bigr)\,dv  \\
= \frac{\Gamma(1-2m)\Gamma(1+l)\Gamma(-l)}{\bigl(\Gamma(1-m)\bigr)^2}
= - \frac{\pi}{\sin(\pi l)} \frac{\Gamma(1-2m)}{\bigl(\Gamma(1-m)\bigr)^2}.
\end{multline}
Hence the ratio of the coefficients of $f^{\pm}_{\lambda,1-\lambda,\lambda}(Y)$
in the expansions of $\deg_Y \bigl( \pi_l^0(\gamma)t^l_{-m\,\underline{m}} \bigr)(Y)$
and $\sgnx(Y) \cdot \bigl(\pi_l^0(\gamma) t^l_{-m\,\underline{m}}\bigr)(Y)$ is
$$
2i \cot(\pi l) \cdot \frac{\bigl(\Gamma(3/2-m)\bigr)^2}
{\bigl(\Gamma(1-m)\bigr)^2}
\cdot \frac{\bigl(\Gamma(\lambda/2+1/2)\bigr)^2}
{\Gamma(\lambda/2) \cdot \Gamma(\lambda/2+1)}
= \lambda \tanh (\pi\im l).
$$

If $m$ is an integer, we need $P_{\lambda}^{(2m)}(x)$ to be even, so we can take
$\lambda =-2m$; then $P_{\lambda}^{(2m)}(\cos\theta)$ is proportional to
$\sin^{-2m}\theta$.
As was already computed in (\ref{coeff-calc}),
$$
2 \frac{\Gamma(l+m+1)(-2m)!}{\Gamma(l-m+1)}
\int_0^{\pi/2} \mathfrak{P}^l_{-m\,m}\Bigl(\frac 2{\cos^2\theta} -1\Bigr)
\cdot \sin^{-2m}\theta \cdot \tan\theta \,d\theta
= \frac{- \pi}{\sin(\pi l)} \frac{\Gamma(1-2m)}{\bigl(\Gamma(1-m)\bigr)^2}.
$$
In the other case we need to deal with convergence issues, so we observe that
$$
\mathfrak{P}^l_{-m\,m}\Bigl(\frac 2{\cos^2\theta} -1\Bigr) =
\lim_{r \to 0^+} \mathfrak{P}^l_{-m\,m}\Bigl(\frac 2{\cos^2\theta} -1\Bigr)
\cdot \cos^{2r} \theta
$$
as a distribution in $\im l$.
Then substituting $v=\tan^2\theta$,
\begin{multline*}
2 \frac{\Gamma(l+m+1)(-2m)!}{\Gamma(l-m+1)}
\int_0^{\pi/2} \mathfrak{P}^l_{-m\,m}\Bigl(\frac 2{\cos^2\theta} -1\Bigr)
\cdot \cos^{2r-1}\theta \cdot \sin^{-2m}\theta \cdot \tan\theta \,d\theta \\
= \int_0^{\infty} v^{-2m} \cdot (1+v)^{m-1/2-r} \cdot
\,_2F_1\bigl( l-m+1,-l-m;-2m+1; -v \bigr)\,dv  \\
= \frac{\Gamma(1-2m)\Gamma(1/2+l+r)\Gamma(-1/2-l+r)}
{\bigl(\Gamma(1/2-m+r)\bigr)^2}
\quad \xrightarrow{\text{as $r\to 0^+$}} \quad
\frac{-2\pi}{\cos(\pi l)}
\frac{\Gamma(1-2m)}{(2l+1) \bigl(\Gamma(1/2-m)\bigr)^2}
\end{multline*}
as a distribution in $\im l$.
Hence the ratio of the coefficients of $f^{\pm}_{\lambda,-\lambda,\lambda}(Y)$
in the expansions of $\deg_Y \bigl( \pi_l^0(\gamma)t^l_{-m\,\underline{m}} \bigr)(Y)$
and $\sgnx(Y) \cdot \bigl(\pi_l^0(\gamma) t^l_{-m\,\underline{m}}\bigr)(Y)$ is
$$
-2i \tan(\pi l) \cdot
\frac{\bigl(\Gamma(1-m)\bigr)^2}{\bigl(\Gamma(1/2-m)\bigr)^2}
\cdot \frac{\bigl(\Gamma(\lambda/2+1/2)\bigr)^2}
{\Gamma(\lambda/2) \cdot \Gamma(\lambda/2+1)}
= \lambda \coth (\pi\im l).
$$
\qed

\section{Appendix: Review of Special Functions}

In this appendix we review special functions used in this paper:
the spherical harmonics $Y_l^m(\theta,\phi)$,
the associated Legendre functions $P_l^{(m)}(x)$,
the associated Legendre functions of the second kind $Q_l^{(m)}(x)$
and some of their identities.
References: \cite{Ba}, \cite{GR}, \cite{V}.

\subsection{Spherical Harmonics}

Let $\BB R^3$ be the three-dimensional space with coordinates $(x^1,x^2,x^3)$.
We denote the Laplacian by
$$
\Delta_3 = \frac{\partial^2}{(\partial x^1)^2} +
\frac{\partial^2}{(\partial x^2)^2} + \frac{\partial^2}{(\partial x^3)^2}.
$$
We use the spherical coordinates
\begin{center}
\begin{tabular}{lcc}
$x^1 = r \sin \theta \cos \phi$ & \qquad & $r \ge 0$  \\
$x^2 = r \sin \theta \sin \phi$ & \qquad & $0 \le \phi < 2\pi$  \\
$x^3 = r \cos \theta$ & \qquad & $0 \le \theta < \pi$.
\end{tabular}
\end{center}
In these coordinates
\begin{equation}  \label{Delta-spherical}
\Delta_3 f = \frac 1{r^2} \cdot \frac{\partial}{\partial r}
\biggl( r^2 \frac{\partial f}{\partial r} \biggr)
+ \frac 1{r^2 \sin \theta} \cdot \frac{\partial}{\partial \theta}
\biggl( \sin \theta \frac{\partial f}{\partial \theta} \biggr)
+ \frac 1{r^2 \sin^2 \theta} \cdot \frac{\partial^2 f}{\partial \phi^2}.
\end{equation}

Let ${\bf t}^l_{m\,n}$ denote the matrix coefficients of $SU(2)$:
\begin{multline*}
{\bf t}^l_{m\,n} (X) =
\frac 1{2\pi i} \sqrt{\frac{(l-m)!(l+m)!}{(l-n)!(l+n)!}}
\oint (sx_{11}+x_{21})^{l-n} (sx_{12}+x_{22})^{l+n} s^{-l+m} \,\frac{ds}s,  \\
X=\begin{pmatrix} x_{11} & x_{12} \\ x_{21} & x_{22} \end{pmatrix} \in SU(2),
\qquad \begin{matrix} l=0, \frac 12, 1, \frac 32, 2, \dots, \\
m,n = -l, -l+1, \dots, l, \end{matrix}
\end{multline*}
where the integral is taken over a loop in $\BB C$ going once around the origin
in the counterclockwise direction (see \cite{V}).
We define the {\em spherical harmonics} on $\BB R^3$ by setting
$$
Y_l^m(\theta,\phi) = {\bf t}^l_{m\,0}(g),
\qquad l=0, 1, 2, \dots, \quad m = -l, -l+1, \dots, l,
$$
where
\begin{equation}  \label{g}
g= \begin{pmatrix} \cos \frac{\theta}2 \cdot e^{i\frac{\phi+\psi}2} &
i \sin \frac{\theta}2 \cdot e^{i\frac{\phi-\psi}2} \\
i \sin \frac{\theta}2 \cdot e^{i\frac{\psi-\phi}2} &
\cos \frac{\theta}2 \cdot e^{-i\frac{\phi+\psi}2} \end{pmatrix}
=
\begin{pmatrix} e^{i\frac{\phi}2} & 0 \\ 0 & e^{-i\frac{\phi}2} \end{pmatrix}
\begin{pmatrix} \cos \frac{\theta}2 & i \sin \frac{\theta}2 \\
i \sin \frac{\theta}2 & \cos \frac{\theta}2 \end{pmatrix}
\begin{pmatrix} e^{i\frac{\psi}2} & 0 \\ 0 & e^{-i\frac{\psi}2} \end{pmatrix}
\end{equation}
with $0 \le \phi < 2\pi$, $0 < \theta < \pi$, $-2\pi \le \psi < 2\pi$.
Using induction on $l$, relation between the spherical harmonics and the
associate Legendre functions (\ref{Y-P-relation}) and recursive relation
(\ref{id1}), one can prove:

\begin{lem}  \label{rY-poly}
The functions $r^l \cdot Y_l^m(\theta,\phi)$ are polynomial in
$x^1,x^2,x^3$ of degree $l$.
\end{lem}


The functions $Y_l^m(\theta,\phi)$'s satisfy
\begin{equation}  \label{Y-egf}
\Delta_3 Y_l^m(\theta,\phi) = -\frac{l(l+1)}{r^2} \cdot Y_l^m(\theta,\phi)
\end{equation}
and form a complete orthogonal basis of functions on the unit sphere $S^2$.
They also satisfy an orthogonality relation
\begin{equation}  \label{Y-orthogonality}
\int_{\theta=0}^{\pi} \int_{\phi=0}^{2\pi} 
Y_l^m(\theta,\phi) \cdot Y_{l'}^{-m'}(\theta,\phi) \cdot
\sin \theta \,d\phi d\theta
= (-1)^m \frac{4\pi}{2l+1} \cdot \delta_{l,l'} \cdot \delta_{m,m'},
\end{equation}
which follows from (\ref{Y-P-relation}), (\ref{P-orthogonality_2}) and
(\ref{m-m}) using substitution $x=\cos \theta$.
From (\ref{Delta-spherical}) and (\ref{Y-egf}) we obtain two classes of
solutions of the harmonic equation $\Delta_3f=0$:
\begin{equation}  \label{Y-harm}
r^l \cdot Y_l^m(\theta,\phi) \qquad \text{and} \qquad
r^{-l-1} \cdot Y_l^m(\theta,\phi), \qquad l=0,1,2,\dots,\quad -l \le m \le l.
\end{equation}

\subsection{Legendre Functions}

The {associated Legendre functions} $P_l^{(m)}(x)$ and
{the associated Legendre functions of the second kind} $Q_l^{(m)}(x)$
are two linearly independent solutions of the second order differential
equation
$$
(1-x^2) \frac {d^2 f}{dx^2} - 2x \frac {df}{dx} +
\Bigl( l(l+1) - \frac {m^2}{1-x^2} \Bigr) f =0.
$$
While the parameters $l$ and $m$ can be arbitrary complex numbers,
we will only be interested in $l=0,1,2,\dots$, $m \in \BB Z$, $-l \le m \le l$.
The {\em associated Legendre functions (of the first kind)} are functions
defined by the formula
\begin{equation}  \label{P-def}
P_l^{(m)}(x) = \frac {(-1)^{l+m}}{2^l l!} (1-x^2)^{m/2}
\frac{d^{l+m}}{dx^{l+m}} (1-x^2)^l.
\end{equation}
If we set $m=0$ we get {\em Legendre polynomials}:
$$
P_l(x) = P_l^{(0)}(x)
= \frac {(-1)^l}{2^l l!} \frac{d^l}{dx^l} (1-x^2)^l,
\qquad x \in \BB C.
$$
Then, for $m \ge 0$,
$$
P_l^{(m)}(x) = (-1)^m (1-x^2)^{m/2} \frac{d^m}{dx^m} P_l(x).
$$
Functions $P_l^{(m)}$ and $P_l^{(-m)}$ are related as follows:
\begin{equation}  \label{m-m}
P_l^{(-m)}(x) = (-1)^m \frac{(l-m)!}{(l+m)!} \cdot P_l^{(m)}(x).
\end{equation}

When $m$ is even, $P_l^{(m)}(x)$ is a polynomial of degree $l$ and defined
on the entire complex plane. On the other hand, when $m$ is odd, $P_l^{(m)}(x)$
is defined as long as $\sqrt{1-x^2}$ is single-valued.
There are different conventions for choosing the domain of $P_l^{(m)}(x)$,
ours is that $P_l^{(m)}(x)$ is defined on the complex plane away from the two
cuts along the real line: along $(-\infty,-1]$ and $[1,\infty)$.
(But, for example, in \cite{Ba} the associated Legendre functions are defined
on the complex plane with the interval $[-1,1]$ removed from the real line.)

The {\em associated Legendre functions of the second kind} are defined
on the complex plane
with the interval $[-1,1]$ removed from the real line.
One way to define these functions is by declaring
$$
Q_0(z) = Q_0^{(0)}(z) = \frac 12 \log \biggl( \frac{z+1}{z-1} \biggr),
\qquad
Q_1(z) = Q_1^{(0)}(z) = \frac z2 \log \biggl( \frac{z+1}{z-1} \biggr) -1,
$$
and then by recursive relations
$$
(2l+1)z \cdot Q^{(m)}_l(z) =
(l-m+1) \cdot Q^{(m)}_{l+1}(z) + (l+m) \cdot Q^{(m)}_{l-1}(z),
$$
which is formally identical to (\ref{id1}),
$$
Q_l^{(m)}(z) = (z^2-1)^{m/2} \frac{d^m}{dz^m} Q_l(z),
\qquad m \ge 0,
$$
and
$$
Q_l^{(-m)}(z) = \frac{(l-m)!}{(l+m)!} \cdot Q_l^{(m)}(z).
$$

Finally, we introduce functions $\tilde Q_l^{(m)}(z)$ defined on the complex
plane away from the two cuts along the real line: along $(-\infty,-1]$ and
$[1,\infty)$ (same domain as for $P_l^{(m)}(z)$'s).
We set
$$
\tilde Q_l^{(m)}(x) = \frac {i^m}2 Q_l^{(m)}(x+i0) + \frac {i^{-m}}2 Q_l^{(m)}(x-i0),
\qquad x \in (-1,1) \subset \BB R.
$$
From the identity
$$
Q^{(m)}_l(x \pm i0) =
i^{\mp m} \Bigl( \tilde Q^{(m)}_l(x) \mp \frac{\pi i}2 P^{(m)}_l(x) \Bigr),
\qquad x \in (-1,1),
$$
we see that $\tilde Q_l^{(m)}(x)$ extends analytically to the upper and lower
half-planes by
\begin{equation}  \label{Q-tilde}
\tilde Q^{(m)}_l(z) = \begin{cases}
i^{m} Q^{(m)}_l(z) + \frac {\pi i}2 P^{(m)}_l(z) & \text{if $\im z>0$;}\\
i^{-m} Q^{(m)}_l(z) - \frac {\pi i}2 P^{(m)}_l(z) & \text{if $\im z<0$.}
\end{cases}
\end{equation}

\subsection{Identities}



The associated Legendre functions satisfy two kinds of orthogonality relations:
\begin{equation} \label{P-orthogonality_2}
\int_{-1}^1 P_k^{(m)}(x) \cdot P_l^{(m)}(x) \,dx = 
\frac {2(l+m)!}{(2l+1)(l-m)!} \cdot \delta_{k,l},
\qquad 0 \le m \le l,
\end{equation}
and
\begin{equation} \label{P-orthogonality}
\int_{-1}^1 P_l^{(m)}(x) \cdot P_l^{(n)}(x) \,\frac{dx}{1-x^2} =
\begin{cases}
0 & \text{if $m \ne n$;} \\
\frac {(l+m)!}{m(l-m)!} & \text{if $m=n \ne 0$;} \\
\infty & \text{if $m=n=0$;}
\end{cases} \qquad 0 \le m,n \le l.
\end{equation}
There is a recursive relation for the associated Legendre functions:
\begin{equation}  \label{P-recursive-1}
\frac d{dx} P^{(m)}_l(x) + \frac {mx}{1-x^2} \cdot P^{(m)}_l(x)
= - \frac 1{\sqrt{1-x^2}} \cdot P^{(m+1)}_l(x).
\end{equation}
Combining it with (\ref{m-m}) we obtain:
\begin{equation}  \label{P-recursive-2}
\frac d{dx} P^{(m)}_l(x) - \frac {mx}{1-x^2} \cdot P^{(m)}_l(x)
= \frac {(l+m)(l-m+1)}{\sqrt{1-x^2}} \cdot P^{(m-1)}_l(x).
\end{equation}
We will also use the following identities:
\begin{equation}  \label{P(-x)}
P^{(m)}_l(-x) = (-1)^{l+m} \cdot P^{(m)}_l(x),
\end{equation}
\begin{equation}  \label{id1}
(2l+1)x \cdot P^{(m)}_l(x) =
(l-m+1) \cdot P^{(m)}_{l+1}(x) + (l+m) \cdot P^{(m)}_{l-1}(x),
\end{equation}
\begin{equation}  \label{id2}
(1-x^2) \cdot \frac d{dx} P^{(m)}_l(x) =
- lx \cdot P^{(m)}_l(x) + (l+m) \cdot P^{(m)}_{l-1}(x),
\end{equation}
\begin{equation}  \label{Qm-m}
\tilde Q_l^{(-m)}(x) = (-1)^m \frac{(l-m)!}{(l+m)!} \cdot \tilde Q_l^{(m)}(x),
\qquad x \in (-1,1),
\end{equation}
\begin{equation}   \label{Q-Q=Pn}
i^m Q_l^{(m)}(x+i0) - i^{-m} Q_l^{(m)}(x-i0) = - \pi i P_l^{(m)}(x),
\qquad x \in (-1,1).
\end{equation}

The associated Legendre functions and the spherical harmonics on $\BB R^3$
are related as follows:
\begin{equation}  \label{Y-P-relation}
Y_l^m(\theta,\phi) = {\bf t}^l_{m\,0}(g) =
(-i)^m \cdot \sqrt{\frac{(l-m)!}{(l+m)!}}
\cdot P_l^{(m)}(\cos \theta) \cdot e^{-im\phi},
\end{equation}
where $g$ is as in (\ref{g}).
From the multiplicativity property for the matrix coefficients
${\bf t}^l_{m\,n}$'s and (\ref{Y-P-relation}) one can obtain
\begin{equation}  \label{P-sum}
\sum_{m=-l}^l \frac{(l-m)!}{(l+m)!} e^{im\phi} \cdot
P_l^{(m)}(\cos \theta) \cdot P_l^{(m)}(\cos \tilde\theta)
= P_l(\cos \gamma),
\qquad \phi \in \BB C,
\end{equation}
where $\cos \gamma =
\cos \theta \cos \tilde\theta + \sin \theta \sin \tilde\theta \cos \phi$.
There also is a similar formula involving the associated Legendre functions
of the second kind:
\begin{equation}  \label{Q-sum}
\sum_{m=-l}^l \frac{(l-m)!}{(l+m)!} e^{im\phi} \cdot
P_l^{(m)}(\cos \theta) \cdot \tilde Q_l^{(m)}(\cos \tilde\theta)
= \tilde Q_l(\cos \gamma),
\qquad \phi \in \BB R,
\end{equation}
which follows from (\ref{Qm-m}) and
$$
\tilde Q_l(\cos\gamma) = P_l(\cos \theta) \cdot \tilde Q_l(\cos\tilde\theta)
+ 2 \sum_{m=1}^l \frac{(l-m)!}{(l+m)!}
P_l^{(m)}(\cos\theta) \cdot \tilde Q_l^{(m)}(\cos\tilde\theta) \cdot \cos(m\phi).
$$
We will also use the following relation:
\begin{equation}  \label{PQ-sum}
\sum_{l=0}^{\infty} (2l+1) P_l(x) Q_l(y) = \frac 1{y-x},
\end{equation}
the sum converges uniformly when $x$ ranges over compact subsets lying
inside the ellipse passing through $y$ and foci at the points $\pm 1$.

\separate

\separate

\noindent
{\em Department of Mathematics, Yale University,
P.O. Box 208283, New Haven, CT 06520-8283}\\
{\em Department of Mathematics, Indiana University,
Rawles Hall, 831 East 3rd St, Bloomington, IN 47405}

\end{document}